\documentclass[12pt]{article}

\usepackage{amssymb}
\usepackage{ecltree}
\usepackage{epic}
\usepackage{fancybox}
\usepackage{subfigure}
\newenvironment{Bundle}[1]{\begin{bundle}{\vspace{-2mm}%
\rule[-1mm]{0mm}{1mm}#1}}{\end{bundle}}
\newcommand{\E}{{\cal E}}

\newcommand{\ebeta}{{\cal E}_{\beta}}
\newcommand{\eu}{{\cal E}^u}

\newtheorem{defn}{Definition}[section]

\newtheorem{cor}[defn]{Corollary}
\newtheorem{thm}[defn]{Theorem}
\newtheorem{prop}[defn]{Proposition}
\newtheorem{lemma}[defn]{Lemma}

\newtheorem{remark}[defn]{Remark}
\newtheorem{example}[defn]{Example}
\newtheorem{hypotheses}[defn]{Hypotheses}
\newtheorem{proc}[defn]{Procedure}

\newcommand{\brmk}{\begin{remark}\per\begin{em}}
\newcommand{\ermk}{\end{em}\end{remark}}
\newcommand{\bexa}{\begin{example}\per\begin{em}}
\newcommand{\eexa}{\end{em}\end{example}}
\newcommand{\bhyp}{\begin{hypotheses}\per\begin{em}}
\newcommand{\ehyp}{\end{em}\end{hypotheses}}
\newcommand{\bproc}{\begin{proc}\begin{em}}
\newcommand{\eproc}{\end{em}\end{proc}}

\def\theequation{\arabic{section}.\arabic{equation}}
\def\thedefn{\arabic{section}.\arabic{defn}}
\newcommand{\beginsec}{\setcounter{equation}{0}}

\newcounter{bean}
\newcommand{\benuma}{\setlength{\labelwidth}{.25in}
\begin{list}
{(\alph{bean})}{\usecounter{bean}}}
\newcommand{\eenuma}{\end{list}}

\newcommand{\be}{\begin{equation}}
\newcommand{\ee}{\end{equation}}

\newcommand{\bea}{\begin{eqnarray}}
\newcommand{\eea}{\end{eqnarray}}

\newcommand{\beas}{\begin{eqnarray*}}
\newcommand{\eeas}{\end{eqnarray*}}

\newcommand{\goto}{\rightarrow}

\newcommand{\per}{\hspace{-.072in}{\bf .  }}

\newcommand{\lan}{\langle}
\newcommand{\ran}{\rangle}

\newcommand{\noi}{\noindent}

\newcommand{\skp}{\vspace{\baselineskip}}

\newcommand{\ds}{\displaystyle}

\newcommand{\ink}{\rule{.5\baselineskip}{.55\baselineskip}}

\newcommand{\R}{I\!\!R}

\newcommand{\ZZ}{{Z\!\!\!Z}}

\newcommand{\N}{{I\!\!N}}

\newcommand{\cp}{{\cal P}}

\newcommand{\cx}{{\cal X}}

\newcommand{\nin}{n \in \N}

\newcommand{\vphi}{\varphi}

\newcommand{\ngi}{n \goto \infty}

\newcommand{\can}{\mbox{Can}}
\newcommand{\ccan}{\mbox{{\em Can}}}
\newcommand{\micro}{\mbox{Micro}}
\newcommand{\mmicro}{\mbox{{\em Micro}}}

\newcommand{\bone}{\beta^1}
\newcommand{\btwo}{\beta^2}
\newcommand{\thone}{\tilde{H}^1}
\newcommand{\thtwo}{\tilde{H}^2}
\newcommand{\hnone}{H_{n}^1}
\newcommand{\hntwo}{H_{n}^2}
\newcommand{\uone}{u^1}
\newcommand{\utwo}{u^2}
\newcommand{\jbone}{J_{\bone}}
\newcommand{\sbone}{s_{\bone}}

\newcommand{\sutwo}{s^{\utwo}}
\newcommand{\vphibone}{\vphi_{\bone}}
\newcommand{\ebonebtwo}{{\cal E}_{\bone,\btwo}}
\newcommand{\eboneutwo}{{\cal E}_{\bone}^{\utwo}}
\newcommand{\vphiutwo}{\vphi^{\utwo}}
\newcommand{\euoneutwo}{{\cal E}^{\uone,\utwo}}

\newcommand{\iuoneutwo}{I^{\uone,\utwo}}

\newcommand{\ital}{\textit}
\newcommand{\trm}{\textrm}

\newcommand{\LL}{\ensuremath{\mathcal{L}}}

\newcommand{\ve}{\ensuremath{\varepsilon}}

\newcommand{\X}{\ensuremath{\mathcal{X}}}

\newcommand{\Y}{\ensuremath{\mathcal{Y}}}
\newcommand{\Z}{\ensuremath{\mathcal{Z}}}

\begin{document}

\title{{\bf
Large Deviation Principles and Complete \\ Equivalence and Nonequivalence
Results \\ for Pure and Mixed Ensembles}} 

\author{Richard S. Ellis\thanks{This research
was supported by a
grant from the Department of Energy (DE-FG02-99ER25376)
and by a grant from the National
Science Foundation (NSF-DMS-9700852).}, 
Kyle Haven\thanks{This research
was supported by a
grant from the Department of Energy (DE-FG02-99ER25376).}, and 
Bruce Turkington\thanks{This research
was supported by a grant
from the Department of Energy (DE-FG02-99ER25376) 
and by a grant from the National
Science Foundation (NSF-DMS-9971204).}\\ 
\small{rsellis@math.umass.edu, haven@math.umass.edu, turk@math.umass.edu}\\
Department of Mathematics and Statistics\\ 
University of Massachusetts\\ 
Amherst, MA  01003} 

\maketitle 

\thispagestyle{empty}

\begin{abstract}

We consider a general class of statistical mechanical
models of coherent structures in turbulence, which 
includes models of two-dimensional fluid motion, 
quasi-geostrophic flows, and dispersive waves.  First, 
large deviation principles are proved for the canonical 
ensemble and the microcanonical ensemble.  For each 
ensemble the set of equilibrium macrostates is defined
as the set on which the corresponding rate function 
attains its minimum of 0.  We then present complete 
equivalence and nonequivalence results at the level of 
equilibrium macrostates for the two ensembles.

Microcanonical equilibrium macrostates are characterized 
as the solutions of a certain constrained minimization 
problem, while canonical equilibrium macrostates are
characterized as the solutions of an unconstrained 
minimization problem in which the constraint in the 
first problem is replaced by a Lagrange multiplier.  
The analysis of equivalence and nonequivalence of 
ensembles reduces to the following question in global 
optimization.  What are the relationships between the 
set of solutions of the constrained minimization problem 
that characterizes microcanonical equilibrium macrostates 
and the set of solutions of the unconstrained minimization 
problem that characterizes canonical equilibrium macrostates?

In general terms, our main result is that a necessary and 
sufficient condition for equivalence of ensembles to hold 
at the level of equilibrium macrostates is that it holds 
at the level of thermodynamic functions, which is the case 
if and only if the microcanonical entropy is concave.  The 
necessity of this condition is new and has the following 
striking formulation.  If the microcanonical entropy is 
not concave at some value of its argument, then the 
ensembles are nonequivalent in the sense that the 
corresponding set of microcanonical equilibrium macrostates 
is disjoint from any set of canonical equilibrium 
macrostates.  We point out a number of models of physical 
interest in which nonconcave microcanonical entropies arise.

We also introduce a new class of ensembles called mixed 
ensembles, obtained by treating a subset of the dynamical 
invariants canonically and the complementary set 
microcanonically.  Such ensembles arise naturally in 
applications where there are several independent dynamical 
invariants, including models of dispersive waves for the 
nonlinear Schr\"{o}dinger equation.  Complete equivalence 
and nonequivalence results are presented at the level of 
equilibrium macrostates for the pure canonical, the pure 
microcanonical, and the mixed ensembles.
\end{abstract}

\noi
{\it American Mathematical Society 1991 subject classifications.}  Primary 60F10, Secondary 82B99

\noi

\tableofcontents

\renewcommand{\theequation}{\arabic{section}.\arabic{subsection}.\arabic{equation}}
\renewcommand{\thedefn}{\arabic{section}.\arabic{subsection}.\arabic{defn}}

\newpage
\section{Introduction}
\beginsec

\subsection{Overview}
\beginsec

A wide variety of complex physical systems described by nonlinear
partial differential equations exhibit asymptotic
phenomena that are much too complicated
to study by purely analytic methods.  
In order to gain a fuller understanding of such phenomena,
analytic methods are supplemented by numerical simulations or the systems
are modeled via the formalism of statistical mechanics, which often yields 
uncannily accurate predictions concerning the system's asymptotic behavior.

An important class of complex physical systems for which the formalism of
statistical mechanics provides accurate predictions arises in the study of
turbulence; e.g., 
two-dimensional fluid motions,
quasi-geostrophic flows, two-dimensional magnetofluids, 
plasmas, and dispersive waves.  In each
case important features of the asymptotic behavior of 
the underlying nonlinear partial differential equation---the
two-dimensional Euler equations, the quasi-geostrophic potential vorticity
equation, the magnetohydrodynamic equations, the Vlasov-Poisson equation,
and the nonlinear
Schr\"{o}dinger equation---can be effectively
captured in a statistical mechanical model.
A distinguishing feature
of such systems is that a free evolution from a generic initial condition
exhibits a separation-of-scales behavior: 
coherent structures are formed on 
large scales---e.g., vortices and shears
in the case of fluid motion or solitons 
in the case of dispersive waves---while random fluctuations 
are generated on small scales.  A major goal of any
description of the system, whether analytic, numeric, or statistical,
is to predict the formation, interaction, and persistence of such
coherent structures.

The purpose of the present paper is to provide the theoretical basis for 
statistical mechanical studies of specific models of turbulence that are 
analyzed elsewhere.  These include 
two-dimensional fluids \cite{BouEllTur},
quasi-geostrophic flows \cite{EllHavTur}, and 
dispersive waves \cite{EllJorTur}.  
In each case the model is defined on a fixed flow domain
in terms of a sequence of finite-dimensional
systems indexed by $\nin$.  Coherent structures are studied 
in the continuum limit,
obtained by sending $n \goto \infty$. 
They are characterized by variational principles, the solutions of which
define equilibrium macrostates.
In contrast to the detailed description required
by the associated nonlinear partial
differential equation and by the finite-dimensional systems that 
discretize them, these equilibrium macrostates provide a vastly contracted description.  The variational principles are derived and
analyzed via the theory of large deviations 
and duality theory for concave functions.  

In these models
the sequence of finite-dimensional systems
is defined on a fixed domain in terms of a
long-range interaction with a local mean-field scaling.  In order
to obtain a nontrivial limit, one must scale the inverse temperature by a parameter tending to infinity.  By altering the scaling and making other
superficial changes, our results
can also be applied to classical lattice models
such as the Ising model
of a ferromagnet.  Such models are typically 
defined in terms of the thermodynamic limit
of a sequence of finite-dimensional systems having a finite-range or
summable interaction.  In such applications a basic stochastic process
that arises in the 
large deviation analysis is the empirical field, which has been
studied by a number of authors including \cite{DeuStrZes,FolOre,Geo,Oll}.
Other papers that investigate the equivalence of ensembles in the traditional
thermodynamic or bulk limit include \cite{AizGolLeb} and \cite{SchWeg}.

There is a large literature on the equivalence of ensembles for classical
lattice systems and related models.  
It is reviewed in part in the introduction to \cite{LewPfiSul}, to which the reader is referred for references.  In particular, a number of papers including \cite{DeuStrZes,Geo,LewPfiSul2,RoeZes} 
investigate the equivalence of ensembles using the theory
of large deviations.  Of these papers, \cite{LewPfiSul2} considers
the problem in the greatest generality, obtaining a criterion for the 
equivalence of ensembles in terms of 
the vanishing of the specific information gain
of a sequence of conditioned measures with respect
to a sequence of tilted measures.
However, despite the mathematical
sophistication of these and other studies, none of them explicitly addresses
the general issue of the nonequivalence of ensembles, which seems
to be the typical behavior for the models of turbulence that the
present paper analyzes.  In \cite[\S7.3]{LewPfiSul2} and \cite[\S7]{LewPfiSul} 
there is a discussion of the nonequivalence of ensembles for the 
simplest mean-field model in statistical mechanics; namely, the
Curie-Weiss model of a ferromagnet.  For a general class of local 
mean-field models
of turbulence, the present paper addresses this and related issues.

In much of the classical literature on statistical mechanical approaches to
two-dimensional turbulence, it is tacitly assumed that the microcanonical
and canonical ensembles give equivalent results \cite{Kra,MilWeiCro}.
Recently, however, in the context of the point vortex and related models, this
tacit assumption has been directly addressed.  
Questions concerning the equivalence and nonequivalence
of ensembles for these models 
have been investigated by a number of authors, including
\cite{CagLioMarPul2,EyiSpo,Kie,KieLeb}.  The present paper,
inspired in part by \cite{EyiSpo}, is the first to present complete and definitive results for a general class of models,
with a particular emphasis upon the nonequivalence of ensembles.

An unexpected connection of our work in this paper is to
dynamic stability analysis.  To date, all studies of the nonlinear stability
of two-dimensional flows have been carried out using the Lyapunov
functionals introduced by Arnold \cite{Arn,ArnKhe,MarPul}.
When these deterministic results are reformulated in the setting
of statistical mechanical models, they can be expressed in terms of 
the second-order conditions
satisfied by canonical equilibrium macrostates.  In the cases when
the microcanonical entropy is not concave 
and thus the ensembles are nonequivalent, the Arnold sufficient conditions
for nonlinear stability are not satisfied by the microcanonical
equilibrium macrostates.  Nevertheless, 
the second-order conditions
satisfied by these macrostates allow us to refine
the classical Arnold theorems by proving the nonlinear stability of 
a new class of two-dimensional flows.  In \cite{EllHavTur} 
these ideas are developed
for the quasi-geostrophic potential vorticity equation, which describes
the dynamics of rotating, shallow water systems in nearly geostrophic 
balance.  The work in that paper has possible applications to the stability
of planetary flows; specifically, to the stability of zonal shear flows and
embedded vortices in Jovian-type atmospheres.

In the next two subsections we present an overview of 
the main results in this paper, stripped of all technicalities.  
This is done in the context of a well-known statistical mechanical model 
of the two-dimensional Euler equations known
as the Miller-Robert model.  Results formulated in great
generality to apply to this and other models of turbulence
are given in Sections 2-5 of this paper.
We start by presenting large deviation principles with respect to the canonical ensemble and the microcanonical ensemble.  For each ensemble
we then define the set of equilibrium
macrostates as the set on which the associated
rate function attains its minimum of 0.  A fundamental question
arises.  Are the two ensembles equivalent at
the level of equilibrium macrostates?  That is, does
each equilibrium macrostate
with respect to one ensemble correspond to an equilibrium macrostate
with respect to the other ensemble?  In Section 4, 
definitive and sharp results on the equivalence
and nonequivalence of the ensembles are presented.
 
In general terms, our main result is that a necessary
and sufficient condition for the equivalence of ensembles
to hold at the level of equilibrium macrostates is that 
it holds at the level of thermodynamic functions.  
In proving this, we go beyond
the important work in \cite{LewPfiSul2},
which proves that for a general 
class of models including the classical lattice gas
thermodynamic equivalence of ensembles is a sufficient condition for macrostate
equivalence of ensembles.  
Our proof that thermodynamic equivalence is also a necessary
condition for macrostate equivalence 
is perhaps the most striking discovery of our work.  
Specifically, we show that 
whenever a quantity known as the microcanonical entropy is not concave, 
the ensembles are nonequivalent in the sense that the set
of microcanonical equilibrium macrostates is richer than the set of
canonical equilibrium macrostates.  In fact, the latter set contains
none of the microcanonical
equilibrium macrostates corresponding to nonconcave portions
of the entropy [see Thm.\ \ref{prop:z}(b)].  
Useful, but less concrete, connections between the nonconcavity of the microcanonical entropy
and nonequivalence of ensembles can also be deduced from 
the abstract results in \cite{LewPfiSul2} [see their \S5 and \S6].
On the other hand, our results are formulated in order to 
apply directly to 
statistical mechanical models of turbulence for which nonconcave
microcanonical entropies frequently and 
naturally arise, particularly in physically
interesting regions corresponding to a range of negative temperatures.  
Several such examples are mentioned in Section 1.4.

Besides the results on equivalence and nonequivalence of ensembles, we also
prove that for the Miller-Robert model
and other models microcanonical equilibrium macrostates have an equivalent
characterization in terms of constrained
maximum entropy principles (see Remark \ref{rmk:equiv}).
Our approach to this question 
seems simpler and more intuitive than the approach taken 
in \cite{MicRob,Rob1,Rob2}.  The 
derivation of constrained maximum entropy principles based on the microcanonical
ensemble brings to fruition the work begun in  \cite{BouEllTur}, where
unconstrained maximum entropy principles 
based on the canonical ensemble are derived.  
Our proof that microcanonical equilibrium macrostates are
characterized as solutions of constrained maximum 
entropy principles is an important contribution because
such principles are the basis for numerical computations of equilibrium
macrostates and coherent structures for the Miller-Robert model and other models
\cite{DibMajTur,TurWhi,WhiTur}.  

In systems having multiple conserved quantities,
one also has the option of studying mixed ensembles.  These are defined by
treating a subset of the conserved quantities canonically 
and the complementary subset of conserved quantities microcanonically.  
In Section 5 we derive large deviation principles with respect to such 
ensembles and give complete results on their equivalence and nonequivalence,
at the level of equilibrium macrostates, with the microcanonical ensemble
and the canonical ensemble.
Although mixed ensembles arise naturally in a number of applications, 
they have not
been studied in a general setting in the statistical mechanical literature.

An important
application of mixed ensembles is to the study of dispersive
waves and soliton turbulence 
for the nonlinear Schr\"{o}dinger equation \cite{EllJorTur}.
This equation 
has two conserved quantities, the Hamiltonian and the particle number.  In the associated statistical mechanical model, 
the canonical ensemble cannot be defined because the
partition function does not converge.  Instead, one must consider either a microcanonical ensemble or a mixed ensemble in which the Hamiltonian is
treated canonically and the particle number microcanonically.  
By applying to the mixed ensemble
a large deviation result for Gaussian processes derived in
\cite{EllRos}, in \cite{EllJorTur} 
we are able to justify rigorously a mean-field
theoretic approach
to soliton turbulence presented in \cite{JorTurZir}.
The agreement
between the predictions of the statistical mechanical model and long-time
simulations of the microscopic dynamics is excellent \cite{JorJos}.

\subsection{Ensembles and Large Deviation Principles}
\beginsec
\setcounter{defn}{0}

The Euler equations describe the time evolution of the velocity field
of an inviscid, incompressible fluid in a spatial domain, 
which for simplicity 
we take to be the unit torus $T^2$ with periodic boundary conditions.  
At time $t >0$ the
velocity field at a position $x = (x_1,x_2) \in T^2 $ is denoted
$(v_1(x,t), v_2(x,t))$.  The Euler equations can be cast in the form
of an infinite-dimensional
Hamiltonian system having a family of other conserved quantities called 
generalized enstophies.  A central goal of theoretical, 
numerical, and statistical studies is to relate the asymptotic behavior
of the vorticity $\omega(x,t) \doteq v_{2,x_1}(x,t) - v_{1,x_2}(x,t)$ 
to the formation, interaction, and persistence 
of coherent structures of the fluid motion.

A model that can be used to carry this out was proposed independently by Miller
et.\ al.\ \cite{Mil,MilWeiCro} and Robert et.\ al.\ \cite{Rob2,RobSom} 
and is known as the Miller-Robert model.  
In order to define it, one first discretizes 
the continuum dynamics described by the Euler equations, and then 
in terms of the discretized dynamics one
defines a sequence of statistical equilibrium models on 
suitable finite lattices ${\cal L}_n$ of $T^2$.  
Details are given in part (b) of Example \ref{exa:models}.  
These lattice models describe
the joint probability distributions of certain vorticity random variables
$\zeta(s)$ defined for each site $s \in {\cal L}_n$.   
We denote by $\zeta$
the configuration or microstate
$\{\zeta(s), s \in \LL_n\}$; by $a_n$ the number of sites
in $\LL_n$; by $\Y$ the common
range of $\zeta(s)$; by $H_n(\zeta)$ the Hamiltonian for $\zeta$, which is
a certain quadratic function of the $\zeta(s)$ that approximates the 
continuum Hamiltonian; by $A_n(\zeta)$ the generalized enstrophy of $\zeta$, which approximates the continuum generalized enstrophy; and by $P_n$ the prior distribution of $\zeta$,
which is a certain product measure on the configuration space
$\Y^{a_n}$.  In order to simplify
the present description, 
we absorb $A_n$ in $P_n$; in \cite{EllHavTur} 
a physical justification is given, in the context of a related model, for absorbing the generalized enstrophy $A_n$ in the prior distribution $P_n$.   
Thus for the purpose of this introduction,
the Miller-Robert model is defined in terms of a single conserved quantity,
the Hamiltonian.  As in many other models of turbulence, the Hamiltonian
in the Miller-Robert model has a long-range interaction and incorporates 
a local mean-field scaling. 

For other models of turbulence having the Hamiltonian as the only conserved quantity, much of the following discussion is valid with minimal changes in notation; in particular, the forms of the large deviation principles in the present subsection and the results 
on equivalence and nonequivalence of ensembles 
in the next subsection.  For models having multiple conserved quantities, the following discussion is easily adapted by replacing certain scalars with 
vectors.  The general class of models considered in this paper is defined in terms of the quantities in Hypotheses \ref{hyp:prob}.  In order for a large
deviation analysis of the model to be feasible, these quantities must satisfy
Hypotheses \ref{hyp:general}.  

We begin our overview of the main results in this paper by appealing to 
the formalism of 
equilibrium statistical mechanics, which provides two joint probability distributions for microstates $\zeta \in \Y^{a_n}$.   
The physically fundamental distribution known
as the microcanonical ensemble models the fact that
the Hamiltonian is a constant of the Euler dynamics.  
Probabilistically, this is expressed by conditioning $P_n$
on the energy shell $\{\zeta \in \Y^{a_n} \! : \! H_n(\zeta) = u\}$, 
where $u \in \R$ is determined
by the initial conditions.  However, in order to avoid problems with
the existence of regular conditional probability distributions,
we shall condition $P_n$ on the thickened energy shell
$\{H_n(\zeta) \in [u-r,u+r]\}$, where $r>0$.   
Thus, the microcanonical ensemble is the measure
defined for Borel subsets $B$ of $\Y^{a_n}$ by 
\[
P_n^{u,r}\{B\}= P_n\{B \, | \, H_n \in [u-r,u+r]\} =
\frac{P_n\{B \cap \{H_n \in [u-r,u+r]\}\}}{P_n\{H_n \in [u-r,u+r]\}};
\]
this is well defined provided the denominator in the last expression is 
positive.  The letter $u$ is used in
the definition of the microcanonical ensemble 
rather than the more usual letter $E$
because this is a special case of a general theory that applies to 
models having multiple conserved quantities; for such models $u \in \R$ is replaced by a vector $u$ representing a fixed value of the vector
of conserved quantities.

A mathematically more
tractable joint probability distribution is the canonical ensemble, defined 
for Borel subsets $B$ of $\Y^{a_n}$ by
\[
P_{n,\beta}\{B\}\doteq \frac{1}{Z(n,\beta)} \cdot
 \int_B \exp[-\beta H_n] \, dP_n.
\]
Here $\beta$ is a real number denoting the inverse temperature and
$Z(n,\beta)$ is the partition function 
$\int_{\Y^{a_n}} \exp[-\beta H_n] \, dP_n$.  This is a normalization
constant that makes $P_{n,\beta}$ a
probability measure.  

The main mathematical tool that we shall use to predict the formation of
coherent structures is the theory of large deviations.  
In the case of the Miller-Robert model, a crucial innovation 
implemented in \cite{BouEllTur} for the canonical ensemble is to study the asymptotic behavior of a random probability measure $Y_n(\zeta)$
that is closely related to a 
certain coarse graining of the random vorticity field
(see part (b) of Example \ref{exa:models}).
This coarse graining is defined in terms of the empirical measures
of $\zeta(s)$ for $s$ in certain macrocells of the lattice
$\LL_n$.  $Y_n$ takes values in a certain subset $\X$ of the
space of probability measures on $T^2 \times \Y$.  
Elements $\mu$ of $\X$ are called macrostates.  
While $Y_n$ is basic to analyzing the asymptotic behavior 
of the model, its definition is far from obvious.
For that reason we call $Y_n$ a hidden process and $\X$ a hidden space 
for the Miller-Robert model.  

The hidden process $Y_n$ has two properties that make a large
deviation analysis of the Miller-Robert model possible.  For details,
the reader is referred to \cite{BouEllTur}.  First, an application
of Sanov's Theorem shows that with respect
to the a priori distribution $P_n$, $Y_n$ satisfies the large deviation
principle on $\X$ with rate function $I(\mu)$ given by the relative entropy
of $\mu \in \X$ with respect to a certain base measure.  We record this fact
by the formal notation
\be
\label{eqn:sanov}
P_n\{Y_n \in B(\mu,\alpha)\} 
\approx \exp[-a_n I(\mu)] \ \mbox{ as } n \goto \infty, \alpha \goto 0.
\ee
In this formula $B(\mu,\alpha)$ denotes the open
ball with center $\mu$ and radius $\alpha$ with respect to an appropriate
metric on $\X$.  
Second, there exists a bounded continuous function $\tilde{H}$
mapping $\X$ into $\R$ with the property that uniformly over microstates 
the Hamiltonian $H_n(\zeta)$ is asymptotic to
$\tilde{H}(Y_n(\zeta))$ as $n \goto \infty$; in symbols,
\be
\label{eqn:tildehapprox}
\lim_{n \goto \infty} \sup_{\zeta \in \Y^{a_n}}
|H_n(\zeta) - \tilde{H}(Y_n(\zeta))| = 0.
\ee
$\tilde{H}$ is called the Hamiltonian representation function.

Using (\ref{eqn:tildehapprox}), one
derives from the large deviation principle for the $P_n$-distributions of $Y_n$
the asymptotic behavior of $Y_n$ 
with respect to the two ensembles $P_n^{u,r}$ and $P_{n,a_n \beta}$.
For appropriate values of $u$ and $\beta$
these are expressed by the formal notation 
\be
\label{eqn:microasymp}
P_n^{u,r}\{Y_n \in B(\mu,\alpha)\} 
\approx \exp[-a_n I^u(\mu)] \ \mbox{ as } n \goto \infty, 
r \goto 0, \alpha \goto 0
\ee
and 
\be
\label{eqn:canasymp}
P_{n,a_n \beta}\{Y_n \in B(\mu,\alpha)\} \approx 
\exp[-a_n I_\beta(\mu)] \ \mbox{ as } n \goto \infty, \alpha \goto 0.
\ee
In these formulas 
$I^u$ and $I_\beta$ are rate functions that map $\X$ into $[0,\infty]$
and are defined in terms of the 
relative entropy $I$ appearing in (\ref{eqn:sanov}).  Because the Miller-Robert model is defined in terms of a long-range interaction having 
a local mean-field
scaling, in order to obtain a nontrivial 
asymptotic theory $\beta$ must be scaled
by $a_n$ in the definition of the canonical ensemble $P_{n,\beta}$ 
\cite[\S3]{BouEllTur}.  
For the general formulation of (\ref{eqn:microasymp}) and (\ref{eqn:canasymp})
as large deviation principles for a general class of models, the reader is referred to Theorem \ref{thm:microcan} and Theorem \ref{thm:canon}, respectively.

It is not difficult to motivate the forms of $I^u$ and $I_\beta$.  In order to
do so, we introduce two basic thermodynamic functions, one associated with 
each ensemble.
Since the groundbreaking work
of Lanford on equilibrium macrostates in classical statistical mechanics
\cite{Lan}, it has been recognized that 
the basic thermodynamic
function associated with the microcanonical ensemble is the 
microcanonical entropy $s$.  In terms of the
distribution $P_n\{H_n \in \cdot\}$, this quantity measures
the multiplicity of microstates $\zeta \in \Y^{a_n}$ consistent
with a given energy value $u$.  It is defined by 
\be
\label{eqn:microsz}
s(u) \doteq \lim_{r \goto 0} \, \lim_{n \goto \infty} \frac{1}{a_n} \log
P_n\{H_n \in [u-r,u+r]\}.
\ee
For appropriate values of $u$, 
the limit exists and is given by (\ref{eqn:s(z)}), which
is a variational formula over macrostates $\mu$.
For $\beta \in \R$
the basic thermodynamic function associated with the canonical
ensemble is the canonical free energy
\be
\label{eqn:canonvphi}
\varphi(\beta) \doteq -\lim_{n \goto \infty} \frac{1}{a_n} 
\log Z(n,a_n \beta).
\ee
The limit exists and is given by (\ref{eqn:vphivar}), which is
also a variational formula over macrostates.

We first motivate the form of $I_\beta$.  If $Y_n \in B(\mu,\alpha)$,
then for all sufficiently small $\alpha$ 
and all sufficiently large $n$ (\ref{eqn:tildehapprox}) implies that
\[
H_n(\zeta) \approx \tilde{H}(Y_n(\zeta)) \approx \tilde{H}(\mu).
\]
Hence for all sufficiently small $\alpha$ and all sufficiently large $n$, 
the asymptotic formula (\ref{eqn:sanov}) and the definition of
$\varphi$ yield
\begin{eqnarray*}
P_{n,a_n \beta}\{Y_n \in B(\mu,\alpha)\} & \doteq & \frac{1}{Z(n,\beta)} 
\int_{\{Y_n \in B(\mu,\alpha)\}} \exp[-a_n \beta H_n] \, dP_n \\
& \approx & \frac{1}{Z(n,\beta)} \exp[-a_n \beta \tilde{H}(\mu)] \, P_n\{Y_n \in 
B(\mu,\alpha)\} \\ 
& \approx & \exp[-a_n (I(\mu) + \beta \tilde{H}(\mu) - \varphi(\beta))].
\end{eqnarray*}
Comparing this with the desired asymptotic form (\ref{eqn:canasymp}) 
motivates the formula
\be
\label{eqn:motivateibeta}
I_\beta(\mu) = I(\mu) + \beta \tilde{H}(\mu) - \varphi(\beta).
\ee
The actual proof of the large deviation principle
for the $P_{n,a_n \beta}$-distributions of $Y_n$ with this rate function
follows the sketch presented here and is not difficult.
Related large deviation principles have been obtained by numerous authors.

We now motivate the form of $I^u$.  Suppose that $\tilde{H}(\mu) = u$.
Then for all sufficiently large $n$ depending on $r$ the set
of $\zeta$ for which both $Y_n(\zeta) \in B(\mu,\alpha)$ and $H_n(\zeta) 
\in [u-r,u+r]$ is approximately equal to the set of $\zeta$ for
which both $Y_n(\zeta) \in B(\mu,\alpha)$ 
and $\tilde{H}(Y_n(\zeta)) \in [u-r,u+r]$.
Since $\tilde{H}$ is continuous and 
$\tilde{H}(\mu) = u$, for all sufficiently small $\alpha$ compared to $r$
this set reduces to $\{\zeta : Y_n(\zeta) \in B(\mu,\alpha)\}$.  
Hence for all sufficiently small $r$, all sufficiently large $n$ 
depending on $r$, and all sufficiently small $\alpha$ compared to $r$,
(\ref{eqn:sanov}) and the definition (\ref{eqn:microsz}) of $s$ yield
\begin{eqnarray*}
P_n^{u,r}\{Y_n \in B(\mu,\alpha)\} & \doteq &
\frac{P_n\{\{Y_n \in B(\mu,\alpha)\} \cap \{H_n \in [u-r,u+r]\}\}}
{P_n\{H_n \in [u-r,u+r]\}} \\
& \approx & \frac{P_n\{Y_n \in B(\mu,\alpha)\}}
{P_n\{H_n \in [u-r,u+r]\}} \\
& \approx & \exp[- a_n (I(\mu) + s(u))].
\end{eqnarray*}
On the other hand, if $\tilde{H}(\mu) \not = u$, then a similar calculation
shows that for all sufficiently small $r$, all sufficiently small $\alpha$,
and all sufficiently large $n$
$\, P_n^{u,r}\{Y_n \in B(\mu,\alpha)\} =0$.  Comparing these approximate
calculations with the desired asymptotic form (\ref{eqn:microasymp}) motivates
the formula
\be
\label{eqn:motivateiz}
I^{u}(\mu) \doteq \left\{ \begin{array}{ll} I(\mu) +
s(u) &  \trm{if }\; \tilde{H}(\mu) = u, \\ 
                                \infty & \trm{if }\; \tilde{H}(\mu) \not = u.
                               \end{array}  \right.  
\ee
In Section 3 we offer two proofs of the 
large deviation principle for the 
$P_n^{u,r}$-distributions of $Y_n$.  Both are straightforward; the first
follows fairly closely the heuristic sketch just given.  Forms of this
large deviation principle are given, for example, 
in \cite{DeuStrZes,LewPfiSul2,LewPfiSul}.

The asymptotic formulas (\ref{eqn:microasymp}) and (\ref{eqn:canasymp}) give rise to several interpretations of the rate functions. 
Through the distributions $P_n^{u,r}\{Y_n \in \cdot\}$ and $P_{n,a_n \beta}\{Y_n \in \cdot\}$, $I^u$ and $I_\beta$ measure the multiplicity
of microstates $\zeta \in \Y^{a_n}$ consistent with a given macrostate $\mu$.
Because of these asymptotic formulas, 
it also makes sense to say that for $i=I^u$
or $i=I_\beta$ a macrostate $\mu_1 \in \X$ 
is more predictable than a macrostate $\mu_2 \in \X$
if $i(\mu_1) < i(\mu_2)$.  Since $i$ is nonnegative, 
the most predictable or most
probable macrostates $\mu$
solve $i(\mu) = 0$.  It is natural to call such $\mu$ equilibrium macrostates.  
Specifically, $\mu \in \X$ satisfying $I^u(\mu) = 0$ is
called a microcanonical equilibrium macrostate; 
$\eu$ denotes the set of all such macrostates.  Analogously, a measure
$\mu \in \X$ satisfying $I_\beta(\mu) = 0$ is called
a canonical equilibrium macrostate; 
$\ebeta$ denotes the set of all such macrostates.   
In terms of equilibrium macrostates $\mu$, one can analyze the formation
of coherent structures by defining the mean vorticity as an appropriate
average of $\mu$ and comparing it, say by simulation, with the long-time
behavior of the vorticity $\omega(x,t) \doteq v_{2,x_1}(x,t) - v_{1,x_2}(x,t)$
as given by the Euler equations \cite{MilWeiCro,RobSom,TurWhi,WhiTur}.

%

\subsection{Equivalence and Nonequivalence of Ensembles}
\beginsec
\setcounter{defn}{0}

The microcanonical ensemble is physically fundamental, and the
canonical ensemble can be heuristically derived from it by considering
a small subsystem of a large reservoir \cite{Bal}.  Aside from
physical considerations concerning which ensemble is more appropriate
in the construction of a statistical model, the more mathematically
tractable canonical ensemble is often introduced as an approximation
to the microcanonical ensemble, which is somewhat difficult
to analyze. However, in order to justify this
use of the canonical ensemble, one must address a basic issue.  At the
level of equilibrium macrostates, do the two ensembles give equivalent
results?  This involves answering the following two questions.
\begin{enumerate}
\item For every $\beta$ and every $\mu$
in the set $\ebeta$ of canonical equilibrium macrostates, does
there exist a value of $u$ such that $\mu$ lies in the set $\eu$ 
of microcanonical equilibrium macrostates?  
\item Conversely, for every
$u$ and 
every $\mu \in \eu$ does there exist a value of $\beta$ such that
$\mu \in \ebeta$?   
\end{enumerate}
Whether or not the answers are yes, a more refined issue is to
determine the precise relationships between 
$\eu$ and $\ebeta$.  For example, if the answers are both yes,
then given $\beta$ in question 1 (resp., $u$ in question 2), how 
does one determine the corresponding value of $u$ (resp., $\beta$)?
It is with these issues, appropriately formulated
in terms of a general class of models
having multiple conserved quantities, that 
Sections 4 and 5 of the present paper is occupied.  
In those sections 
definitive and sharp results on the equivalence and nonequivalence
of ensembles are derived.

As we will see, in general question 1 in the preceding
paragraph has the answer yes; namely,
every $\mu \in \ebeta$ lies in $\eu$ for some value of $u$.  
As we illustrate by a number of examples given in Section 1.4,
question 2 can have the answer no; namely, it can
be the case that the
set of microcanonical equilibrium macrostates is 
richer than the set of canonical equilibrium macrostates. 
As we show in Theorem \ref{thm:z}, this behavior
has a striking formulation in terms of the microcanonical entropy
$s$, which is defined in
(\ref{eqn:microsz}).  If $s$ is not concave 
at a given value of $u$, then the ensembles are nonequivalent
in the sense that $\eu$ is disjoint from the sets
$\ebeta$ for all values of $\beta$.  

This general result has been anticipated in a number of works, including
those discussed in Section 4.2 of \cite{Thi} and in \cite{Kie2,KieNeu}.
These works exhibit nonconcave entropy curves for a number of physical
models that include a gravitating system of fermions and a system of circular
vortex filaments in an ideal fluid confined to a three-dimensional torus;
see Fig.\ 34 in \cite{Thi}, Fig.\ 3 in \cite{Kie2},
and Fig.\ 2 in \cite{KieNeu}.  They also point out that certain equilibrium 
macrostates corresponding to nonconcave portions of the entropy are
only realizable in the continuum limit of the microcanonical ensemble but not of the canonical ensemble.  Other examples of nonconcave entropies are given in 
Section 1.4 of the present paper.

The question as to whether the microcanonical and
canonical ensembles give equivalent results at the level of equilibrium
macrostates is formulated as a problem in global optimization.  
Let $u$ and $\beta$ be given.  By definition, 
a macrostate $\bar{\mu}$ belongs to 
$\E^u$ if and only if $I^u(\bar{\mu}) = 0$.  This is the case if and only if
$\bar{\mu}$ solves the following constrained minimization problem:
\be
\label{eqn:constr}
\mbox{ minimize } I(\mu) \mbox{ over } \mu \in \X 
\mbox{ subject to the constraint } \tilde{H}(\mu) = u;
\ee
it is worth noting that since the relative entropy $I(\mu)$ equals 
negative the physical entropy, this display defines a maximum
entropy principle with the energy constraint $\tilde{H}(\mu) = u$.
By definition, a macrostate $\bar{\mu}$ belongs to 
$\E_\beta$ if and only if $I_\beta(\bar{\mu}) = 0$.
This is the case if and only if
$\bar{\mu}$ solves the following unconstrained minimization problem:
\be
\label{eqn:unconstr}
\mbox{ minimize } (I(\mu) + \beta \tilde{H}(\mu)) 
 \mbox{ over } \mu \in \X.
\ee
In the unconstrained problem $\beta$ is a Lagrange multiplier dual to the
constraint $\tilde{H}(\mu) = u$ in (\ref{eqn:constr}).  Under general
conditions, solutions of the constrained minimization problem (\ref{eqn:constr})
are extremal points of $(I + \beta \tilde{H})$ on $\X$ 
\cite{IofTih,Zei}.  The question as to whether
the microcanonical and canonical ensembles give equivalent results
is equivalent to answering the following refined question related to this property.  What are the relationships between the sets of solutions
of the constrained and unconstrained minimization problems (\ref{eqn:constr})
and (\ref{eqn:unconstr})?

We now describe our results on the equivalence and nonequivalence
of ensembles by relating them to the behavior of the two basic thermodynamic
functions, $s$ and $\varphi$.  The following discussion
applies to the Miller-Robert model as well as to a class of
other models that have the Hamiltonian
as a single conserved quantity.
The discussion generalizes
to a wide class of other models having multiple conserved quantities.  
We first motivate a formula relating $s$ and $\varphi$.  To do this, we
use the definition of $s$, which we summarize by the formula
\[
P_n\{H_n \in du\} \approx \exp[a_n s(u)] \, du.
\]
We now calculate
\begin{eqnarray*}
\varphi(\beta) & = & - \lim_{n \goto \infty} \frac{1}{a_n} \log Z(n,a_n \beta) \\
& = &  -\lim_{n \goto \infty} \frac{1}{a_n} 
\log \int_{\Y^{a_n}} \exp[-a_n \beta H_n] \, dP_n \\
& = & -\lim_{n \goto \infty} \frac{1}{a_n} 
\log \int_{\R} \exp[-a_n \beta u] \, P_n\{H_n \in du\} \\
& \approx & -\lim_{n \goto \infty} \frac{1}{a_n} 
\log \int_{\R} \exp[-a_n (\beta u - s(u))] \, du.
\end{eqnarray*}
According to the heuristic reasoning that underlies Laplace's method, the
main contribution to the integral comes from the largest term.
This motivates
the relationship 
\be
\label{eqn:legfen}
\varphi(\beta) = \inf_{u \in \R} \{\beta u - s(u)\},
\ee
which expresses $\vphi$ as the Legendre-Fenchel transform $s^*$ of $s$.

For the Miller-Robert model and other models of turbulence considered in this
paper, $s$ is nonpositive and upper semicontinuous on $\R$ 
[Prop.\ \ref{prop:ldpj}(a)].  
If it is the case that $s$ is concave
on $\R$, then 
(\ref{eqn:legfen}) can be inverted to give 
$s$ in terms of $\varphi$; namely, for all $u \in \R$
\be
\label{eqn:invertlegfen}
s(u) = \inf_{\beta \in \R} \{\beta u - \varphi(\beta)\}.
\ee
Hence, when 
$s$ is concave on $\R$, each basic thermodynamic
function can be obtained from the other by a similar formula.
It is
natural to say that in this case the microcanonical ensemble and
the canonical ensemble are thermodynamically equivalent \cite{KieLeb,LewPfiSul}.
As we will see in Theorems \ref{thm:z} and \ref{thm:s-concave}, 
thermodynamic equivalence of ensembles is mirrored
by equivalence-of-ensemble relationships 
at the level of equilibrium macrostates.

By virtue of its definition (\ref{eqn:canonvphi}) or formula (\ref{eqn:legfen}),
$\varphi$ is a finite, concave, continuous function on $\R$.  In the case of
classical systems such as considered by Lanford \cite{Lan}, 
a superadditivity argument based on the fact that the underlying Hamiltonian
has finite range shows that the analogue of $s$ is an upper 
semicontinuous, concave function on $\R$.  
In general, however, because of the local mean-field, long-range nature
of the Hamiltonians in the Miller-Robert model and other models of turbulence
considered in this paper, 
the associated microcanonical entropies are typically not concave on subsets
of $\R$ corresponding to a range of negative temperatures.

In order to see how concavity properties of $s$ determine
relationships between the sets of equilibrium macrostates, we define
for $u \in \R$ the concave function \be
\label{eqn:sstarstar}
s^{**}(u) \doteq \inf_{\beta \in \R} \{\beta u - s^*(\beta)\} =
\inf_{\beta \in \R} \{\beta u - \vphi(\beta)\}.  \ee Because of
(\ref{eqn:invertlegfen}), it is obvious that $s$ is concave on $\R$ if
and only if $s$ and $s^{**}$ coincide.  Whenever $s(u) > -\infty$ and
$s(u) = s^{**}(u)$, we shall say that $s$ is concave at $u$.

Now assume that $s$ is not concave on $\R$; i.e., there exists $u \in
\R$ for which $-\infty < s(u) \not= s^{**}(u)$.  In this case, one
easily shows that ${s}^{**}$ equals the smallest upper semicontinuous,
concave function majorizing $s$.  In particular, when $s$ is not
concave on $\R$, it cannot be
recovered from $\varphi$ via a Legendre-Fenchel transform.

As we now explain, concavity and nonconcavity properties
of the microcanonical entropy $s$
have crucial implications for the
equivalence and nonequivalence of ensembles at the level
of equilibrium macrostates.  
In terms of such properties 
of $s$, we now give preliminary and incomplete
statements of the relationships between the sets $\eu$ and $\ebeta$ of
equilibrium macrostates for the two ensembles.  The reader is referred
to Theorems \ref{thm:z}, \ref{thm:beta}, and \ref{thm:full}
for precise statements.  For easy reference
they are summarized in Figure \ref{fig:equiv} in Section 4.

For a given value of $u$, there are three possible relationships that can
occur between $\E^u$ and $\E_\beta$.  If there exists a value 
of $\beta$ such that $\E^u=\E_\beta$, then the ensembles are
said to be fully equivalent.  If instead of equality $\E^u$ is a proper
subset of $\E_\beta$ for some $\beta$, 
then the ensembles are said to be partially
equivalent.  It may also happen that $\E^u \cap \E_\beta = \emptyset$
for all values of $\beta$.  If this occurs, then the
microcanonical ensemble is said to be nonequivalent to any canonical
ensemble or that nonequivalence of ensembles holds.  
It is convenient to group the first two cases together.  If for a given $u$
there exists $\beta$ such that either $\E^u$ equals $\E_\beta$
or $\E^u$ is a proper subset of $\E_\beta$, then the ensembles are said
to be equivalent.

The relationships between $\eu$ and $\ebeta$ depend on concavity 
and nonconcavity properties of $s$, expressed through the equality 
or nonequality of $s(u)$ and $s^{**}(u)$.   
These relationships are given next in items 1-3 together with references 
to where the results are stated precisely.
Criteria for equivalence of ensembles related to item 2 have been obtained
in various settings by a number of authors, including 
\cite{DeuStrZes,EyiSpo,LewPfiSul2,LewPfiSul}.  However, the results
underlying items 1 and 3 are new.
\begin{enumerate}
\item
{\bf Canonical is always microcanonical.}  For every $\beta$ 
and every $\mu \in \ebeta$, there exists $u$ such that $\mu \in \eu$
[Theorem \ref{thm:beta}].
\item
{\bf Equivalence.}  If $-\infty < s(u) = s^{**}(u)$---i.e., if $s$ is concave
at $u$---then there exists $\beta$ such that the ensembles are
equivalent [Remark \ref{rmk:dell} and Theorem \ref{thm:z}(a)].
\item
{\bf Nonequivalence.}  If $-\infty < s(u) \not = s^{**}(u)$---i.e., if $s$ 
is not concave at $u$---then the 
corresponding microcanonical
ensemble is nonequivalent 
to any canonical ensemble [Remark \ref{rmk:dell} and Theorem \ref{thm:z}(b)].
\end{enumerate}

Let $u$ be a point in $\R$ such that $s(u) > -\infty$.  
According to items 2 and 3, 
the ensembles are equivalent if and only if $s$ is concave at $u$.  
Under another natural hypothesis on $u$, one shows that $s$ is concave at $u$ if and only if there exists a supporting line to the graph of $s$ at $(u,s(u))$ [Lem.\
\ref{lem:CtoC}(a)]; i.e., there exists $\beta \in \R$ such that
\[
s(w) \leq s(u) + \beta(w - u) \ \mbox{ for all } w \in \R.
\]
In Theorem \ref{thm:full} we refine this necessary and sufficient condition
for equivalence of ensembles by showing that the ensembles are fully equivalent
if and only if there exists a supporting line to the graph of $s$ that touches
the graph of $s$ only at $(u,s(u))$; i.e., there exists $\beta \in \R$ such that
\[
s(w) < s(u) + \beta(w - u) \ \mbox{ for all } w \not= u.
\]
A sufficient condition that guarantees this property of $s$ is that $s(u)
= s^{**}(u)$ and $s^{**}$ is strictly concave in a neighborhood of $u$.

The relationships given in items 1-3
refine the relationships between the thermodynamic
functions $\varphi$ and $s$.  
In fact, the thermodynamic equivalence of ensembles that holds when 
$s = s^{**}$ on $\R$ is reflected
in the equivalence of ensembles for a given value of $u$
when $-\infty < s(u) = s^{**}(u)$ [item 2].   
On the other hand, when $-\infty < s(u) \not = s^{**}(u)$ 
for some value of $u$, the lack
of symmetry between $\varphi$ and $s$ as expressed by (\ref{eqn:legfen}) and (\ref{eqn:sstarstar}) is mirrored
by a lack of symmetry between the microcanonical and canonical
ensembles at the level of equilibrium macrostates.  
For each $\beta$, every canonical equilibrium macrostate
in $\ebeta$ lies in $\eu$ for some $u$ [item 1].  However,
for any $u$ for which $-\infty < s(u) \not = s^{**}(u)$ 
the corresponding microcanonical
ensemble is nonequivalent to any
canonical ensemble [item 3].  

We also prove a number of interesting results that follow easily from the main
theorems.  For example, in Corollary \ref{cor:beta} we show that if $\ebeta$ consists of a unique macrostate $\mu$, then $\eu$ consists of the unique macrostate $\mu$
for a corresponding value of $u$ ($u = \tilde{H}(\mu)$).  The uniqueness
of an equilibrium macrostate corresponds to the absence of a phase transition.

\subsection{Examples of Nonconcave Microcanonical Entropies}
\beginsec
\setcounter{defn}{0}

The most striking of our results on equivalence and nonequivalence of ensembles
is given in item 3 near the end of the preceding subsection.  If, for a given 
value of $u$, $-\infty < s(u) \not = s^{**}(u)$, then 
$\eu$ is disjoint from the sets $\ebeta$ for all values of $\beta$.
We next point out a number of statistical mechanical models having a nonconcave
microcanonical entropy and thus exhibiting, for a range of values of $u$,
the nonequivalence of ensembles that is formulated in item 3.  
\begin{enumerate}
\item {\bf Point vortex system.} This is the first
  statistical mechanical model proposed in the literature for studying the
  two-dimensional Euler equations.  It is defined in terms of a singular
interaction function, which is a Green's function. The model 
was introduced by Onsager
  \cite{Ons}; was further developed in the 1970's, notably by Joyce
  and Montgomery \cite{JoyMon}; and continues to be the subject of
  important studies, including
  \cite{BodGui,CagLioMarPul1,CagLioMarPul2,Kie,KieLeb}.  Proposition
  6.2 in \cite{CagLioMarPul2} isolates a class of flow domains for
  which the microcanonical entropy in the point vortex model
is not a concave function of its
  argument.  As pointed out in \cite[\S6]{KieLeb}, the Monte Carlo
  study of a point vortex system in a disk carried out in
  \cite{SmiOne} also displays a nonconcave microcanonical entropy.
  Strictly speaking, the results on nonequivalence of ensembles
given in the present paper apply only to a point vortex model in which the
  singular interaction function in the classical model has been
  regularized; see part (a) of Example \ref{exa:models}.
Nevertheless, special arguments can be invoked to
  extend them to the classical model with singular point vortices.
\item {\bf Two-dimensional turbulence.} A natural generalization,
  and also regularization, of the point vortex model is the
  Miller-Robert model.  In an unpublished numerical study, Turkington
  and Liang consider the Miller-Robert model in a disk with
  constraints on the energy, the total circulation, and the angular
  momentum (or impulse) and with a prior distribution on the vorticity
  that corresponds to vortex patch dynamics; this problem is the
  simplest Miller-Robert analogue of the problem studied in
  \cite{SmiOne} in the point-vortex formulation.  For fixed values of
  the total circulation and the angular momentum, 
Turkington and Liang compute microcanonical entropies as
  a function of energy using the algorithm developed in \cite{TurWhi}.
  They find that the microcanonical entropy-energy curve is concave
  on a certain interval and nonconcave on a complementary interval.
  These computations produce equilibrium macrostates that are vortices
  embedded in circular shear flows. 
\item {\bf Quasi-geostrophic turbulence on a $\beta$-plane.} The
  statistical equilibrium models proposed in \cite{Tur} are
  implemented in \cite{DibMajTur} for barotropic, quasi-geostrophic
  flow in a channel on the $\beta$-plane.  Various prior distributions
  on the potential vorticity are considered; these include a saturated
  model, in which the maximum and minimum of the potential vorticity
  constrain the microstate, and a dilute model, in which only the mean
  potential-vorticity magnitude is imposed.  Even in the absense of
  geophysical effects ($\beta=0$), the dilute model exhibits a
  nonconcave entropy-energy curve, as displayed in Figure 4 of
  \cite{DibMajTur}.  The equilibrium macrostates 
corresponding to values of the energy for which the entropy is nonconcave 
are shears that transition to monopolar vortices and then to
  dipolar vortices as the energy increases.  When the dilute model is
  replaced by the corresponding saturated model, in which an upper
  bound on the microscopic potential vorticity is enforced, the 
equilibrium macrostates are modified, particularly at high energies.  As is
  shown in Figure 16 of \cite{DibMajTur}, the nonconcavity of the
  entropy-energy curve persists at low energies; at high enough
  energies, however, it becomes concave, unlike in the dilute case.  At
  these high energies the equilibrium macrostates are not dipolar vortices, but
rather shear flows.
\item {\bf Quasi-geostrophic turbulence over topography.} A more
  complete study of the concavity of the microcanonical entropy is
  carried out in \cite{EllHavTur} for equivalent-barotropic,
  quasi-geostrophic flow over bottom topography on an $f$-plane.  As
  in \cite{DibMajTur} a channel geometry is imposed, but for simplicity
  only shear flows are considered.  Within this symmetry class, the
  topography is chosen to be sinusoidal, the energy and circulation
  are used as global invariants, and the prior distribution is taken
  to be a Gamma distribution with mean 0, variance 1, and nonzero skewness.  
As a function of the energy and the 
  circulation, the entropy is nonconcave in more than half of its domain.  These
  two-constraint results are described in detail in Section 6 of
  \cite{EllHavTur}.
\item {\bf Two-layer quasi-geostrophic turbulence.}  The one-layer 
model studied in
\cite{DibMajTur} is extended to a two-layer system in \cite{DibMaj}, 
where it is used to describe the physically important phenomenon 
of open-ocean convection.  In Figures 2 and 12 in that paper, the 
entropy-energy curve is seen to be nonconcave; the microcanonical
equilibrium macrostates corresponding to values of the energy
in the nonconcave region are asymmetric baroclinic monopoles.
%
%
\end{enumerate}

\subsection{Contents of This Paper} 
\beginsec
\setcounter{defn}{0}

In Section 2 we introduce 
the class of statistical mechanical models that will be analyzed
in this paper.
These models generalize the Miller-Robert model by incorporating
a finite sequence of interaction functions $H_{n,i}$ rather than
just the Hamiltonian.  In order to carry out the large deviation
analysis, we assume that there exists a hidden process $Y_n$ that takes
values in a complete separable metric space $\X$ and has the following two properties:
(a) for each interaction function there exists a representation function
$\tilde{H}_i$ such that uniformly over microstates $|H_{n,i} -
\tilde{H}_i \circ Y_n| \goto 0$ 
as $n \goto \infty$; (b) with respect to the prior
measure $P_n$ in the model, 
$Y_n$ satisfies the large deviation principle on $\X$. 
In Section 2 we show that with respect to the canonical ensemble 
$Y_n$ satisfies the large deviation principle, 
and we derive several properties of the set of 
canonical equilibrium macrostates.  

In Section 3 we consider the 
microcanonical ensemble, proving a large deviation principle
and studying properties of the set of microcanonical
equilibrium macrostates.  
We also 
point out the constrained 
maximum entropy principles that characterize microcanonical
equilibrium macrostates in certain models including the Miller-Robert model.

Section 4 is devoted to the presentation of our complete results
on the equivalence and nonequivalence of the two ensembles.  The results
are proved
in Theorems \ref{thm:z}, \ref{thm:beta}, and \ref{thm:full}
and are summarized in Figure \ref{fig:equiv}.

In Section 5.1 we introduce mixed ensembles
obtained by treating a subset of the
dynamical invariants canonically and the complementary
subset of dynamical invariants 
microcanonically.  We then prove the large deviation principle for
these ensembles.  Section 5.2 presents complete
equivalence and nonequivalence
results for the pure canonical and mixed ensembles while
Section 5.3 does the same for the mixed and the pure microcanonical
ensembles.  The results in Sections 5.2 and 5.3 follow from those in Section 4
with minimal changes in proof.  
They are summarized in Figures \ref{fig:fixed-beta-1} and 
\ref{fig:fixed-z-2}.  

\skp
\noi
{\bf Acknowledgement.}  We thank Michael Kiessling for a number of useful 
conversations.

\renewcommand{\theequation}{\arabic{section}.\arabic{equation}}
\renewcommand{\thedefn}{\arabic{section}.\arabic{defn}}

\skp
\skp
\section{Canonical Ensemble: LDP and Equilibrium Macrostates}
\beginsec

In this section we present a large deviation principle
for the canonical ensemble in a wide range of statistical mechanical models 
 [Thm.\ \ref{thm:canon}(b)].
In terms of that principle,
the set of canonical equilibrium macrostates is defined and some
of its properties derived [Thms.\ \ref{thm:canon}(c)-\ref{thm:weak-limits}].
After defining the class of models
under consideration, we specify in Example \ref{exa:models} a number
of specific models to which the theory applies.   

The models that we consider are defined
 in terms of the following quantities.
\bhyp
\label{hyp:prob}
\begin{itemize}
\item
A sequence of probability spaces $(\Omega_n,{\cal F}_n,P_n)$ 
indexed by $\nin$;
$\Omega_n$ are the configuration spaces for 
the statistical mechanical models.
\item
A positive integer $\sigma$ and for each $\nin$ 
a sequence of interaction functions
$\{H_{n,i}, i=1,\ldots,\sigma\}$, which are bounded 
measurable functions mapping $\Omega_n$ into $\R$.  
We define $H_n \doteq (H_{n,1},\ldots,H_{n,\sigma})$, which 
maps $\Omega_n$ into $\R^\sigma$.
\item
A sequence of positive scaling constants $a_n \goto \infty$.
\end{itemize}
\ehyp

Let $\lan \cdot,\cdot \ran$ denote the Euclidean inner product on $\R^\sigma$.
We define for each $\nin$,
$\beta = (\beta_1,\ldots,\beta_\sigma) \in \R^\sigma$, and set $B \in {\cal F}_n$
the partition function
\[
Z_n(\beta) \doteq \int_{\Omega_n}
\exp\!\left[- \sum_{i=1}^\sigma
\beta_i H_{n,i} \right]  dP_n = \int_{\Omega_n} 
\exp[- \lan \beta, H_n \ran] \, dP_n,
\]
which is well defined and finite, and the probability measure
\be
\label{eqn:definecan}
P_{n,\beta}\{B\} \doteq \frac{1}{Z_n(\beta)}
\int_B \exp[- \lan \beta, H_n \ran] \, dP_n.
\ee
The measures $P_{n,\beta}$ are Gibbs states
that define the canonical ensemble for the given model.  
For $\beta \in \R^\sigma$, we also define 
\[
\vphi(\beta) \doteq -\lim_{n \goto \infty} 
\frac{1}{a_n} \log Z_n(a_n \beta)
\]
if the limit exists and is nontrivial.  
In this formula $\beta$ is scaled with $a_n$, as is usual in studying
the continuum limit of models of turbulence
\cite[\S 3]{BouEllTur}.  
We refer to $\vphi(\beta)$ as the canonical free energy.  
If $\sigma = 1$ and $H_{n,1}$ is the Hamiltonian of the system,
then $\beta = \beta_1$ is the inverse temperature.

The first application of the theory of large deviations in this paper
is to express $\vphi(\beta)$ as a variational formula.
Let $\X$ be a Polish space (a complete separable metric space), 
$Y_n$ random variables mapping
$\Omega_n$ into $\X$, $Q_n$ probability measures
on $(\Omega_n,{\cal F}_n)$,
and $I$ a rate function on $\X$.  Thus 
$I$ maps $\X$ into $[0,\infty]$
and for each $M \in [0,\infty)$ the set
$\{x \in \X : I(x) \leq M\}$ is compact (compact level sets). 
For $A$ a subset of $\X$, we
define $I(A) \doteq \inf_{x \in A}I(x)$.  
We say that with respect to $Q_n$ the sequence $Y_n$ satisfies the large
deviation principle, or LDP, on $\X$ with scaling constants
$a_n$
and rate function $I$ if for any closed subset 
$F$ of $\X$ the large deviation
upper bound 
\be
\label{eqn:ldupper}
\limsup_{n \goto \infty} \frac{1}{a_n} \log Q_n\{Y_n \in F\} \leq - I(F)
\ee
is valid and for any open subset $G$ of $\X$ 
the large deviation lower bound
\be
\label{eqn:ldlower}
\liminf_{n \goto \infty} \frac{1}{a_n} \log Q_n\{Y_n \in F\} \geq - I(G)
\ee
is valid.  
We say that with respect to $Q_n$ the sequence $Y_n$ satisfies the 
Laplace principle on $\X$ with scaling constants $a_n$ 
and rate function $I$ 
if for all bounded continuous functions
$f$ mapping $\X$ into $\R$
\begin{eqnarray*}
\lefteqn{
\lim_{n \goto \infty} \frac{1}{a_n} \log \int_{\Omega_n} \exp[a_n f(Y_n)] \, dQ_n }
\\ && = \lim_{n \goto \infty} \frac{1}{a_n} \log \int_{\X} \exp[a_n f(x)]
\, Q_n\{Y_n \in dx\} = \sup_{x \in \X}\{f(x) - I(x)\}.
\end{eqnarray*}
As pointed out in Theorems 1.2.1 and 1.2.3 in \cite{DupEll}, 
$Y_n$ satisfies
the LDP with scaling constants $a_n$ and rate function $I$
if and only if $Y_n$ satisfies the Laplace principle 
with scaling constants $a_n$
and rate function $I$.  Evaluating the large deviation 
upper bound (\ref{eqn:ldupper}) for $F=\X$ and the 
large deviation lower bound (\ref{eqn:ldlower})
for $G = \X$ yields $I(\X) =0$, and since $I$ is nonnegative
and has compact level sets, the set of $x \in \X$ for which $I(x) = 0$
is nonempty and compact.  In the sequel we shall usually omit the phrase ``with scaling constants $a_n$'' in the statements of LDP's and Laplace principles.

A large deviation analysis of the 
general model is possible provided
we can find, as specified in Hypotheses \ref{hyp:general}, 
a hidden space, a hidden process, and 
a sequence of interaction representation functions, and provided the
hidden process satisfies the LDP on the hidden
space.   

\bhyp
\label{hyp:general}
\begin{itemize}
\item
{\bf Hidden space.} This is a Polish space $\X$.
\item
{\bf Hidden process.} This is a sequence $Y_n$,
where each $Y_n$ is a random variable
mapping $\Omega_n$ into $\X$.
\item
{\bf Interaction representation functions.}  This is a 
sequence $\{\tilde{H}_i, i=1,\ldots,\sigma\}$ of bounded
continuous functions mapping $\X$ into $\R$
such that as $n \goto \infty$ 
\be
\label{eqn:repr-fn}
H_{n,i}(\omega) = \tilde{H}_i(Y_n(\omega)) + \mbox{o}(1)
\: \mbox{ uniformly for } \omega \in \Omega_n;
\ee
i.e.,  $\lim_{n \goto \infty} \sup_{\omega \in \Omega_n} 
\left| H_{n,i}(\omega) - \tilde{H}_i(Y_n(\omega))\right| = 0$.
We define $\tilde{H} \doteq (\tilde{H}_1,\ldots,\tilde{H}_\sigma)$, which
maps $\X$ into $\R^\sigma$. 
\item
{\bf LDP for the hidden process.}
There exists a rate function $I$ mapping $\X$ into $[0,\infty]$
such that with respect to $P_n$ the sequence $Y_n$
satisfies the LDP on $\X$, or equivalently
the Laplace principle on $\X$, with rate function $I$.
\end{itemize}
In this context we use the term ``hidden'' because in many cases
the choices of the space $\X$ and the process $Y_n$ are far from obvious. \ink
\ehyp

We next present several models of turbulence
to which the results of this paper
can be applied.  

\bexa
\label{exa:models}
(a) {\bf Regularized Point Vortex Model.}  This model, analyzed in 
\cite{EyiSpo}, is an approximation
to the point vortex model, which we first define.  Let
$\Lambda$ be a smooth, bounded, connected, open subset of $\R^2$; 
$g(x,x')$ the Green's function for $- \! \bigtriangleup$ on $\Lambda$
with Dirichlet
boundary conditions; $h$ the continuous function mapping $\Lambda$ into $\R$
defined by $h(x) \doteq \frac{1}{2} \tilde{g}(x,x)$, 
where $\tilde{g}(x,x')$ is the regular part of the Green's
function $g(x,x')$; and $\theta$ normalized Lebesgue
measure on $\Lambda$ satisfying $\theta(\Lambda) =1$.  
For $\nin$ the point vortex model is defined on the configuration spaces
$\Omega_n \doteq \Lambda^n$ with the Borel $\sigma$-field.
$P_n$ equals the product measure on $\Omega_n$ with
identical one-dimensional marginals $\theta$, and $a_n \doteq n$.  
Configurations $\zeta \in \Lambda^n$ give the 
locations $\zeta_1,\ldots,\zeta_n$ of the
$n$ vortices.  
The interaction function for the point vortex model is the Hamiltonian 
\be
\label{eqn:hnptv}
H_n(\zeta) \doteq \frac{1}{2n^2} \sum_{1 \leq i < j \leq n} g(\zeta_i,\zeta_j) +
\frac{1}{n^2}\sum_{1 \leq i \leq n} h(\zeta_i).
\ee

Because $g(x,x')$ and $h(x)$ are not bounded continuous functions 
of $x$ and $x'$ in $\Lambda$, the point
vortex model cannot be studied by the methods of this paper,
but must be analyzed by other techniques \cite{BodGui,CagLioMarPul1,CagLioMarPul2,Kie,KieLeb}.  The regularized point 
vortex model is defined like the point vortex model except
that in the formula for $H_n$ $\, g(x,x')$
is replaced by a suitable
bounded continuous function $V(x,x')$ on $\Lambda^2$
and $h$ is replaced by a suitable bounded continuous $k$ on $\Lambda$.

For the regularized point vortex model
the hidden space is the space
$\X$ of probability measures on $\Lambda$ while the hidden
process is the sequence of empirical measures
\[
Y_n(\zeta) = Y_n(\zeta,dx) \doteq \frac{1}{n} 
\sum_{i = 1}^n \delta_{\zeta_i}(dx).
\]
By Sanov's Theorem, this sequence satisfies the large deviation
principle on $\X$ with rate function the relative
entropy $R(\mu|\theta)$ of $\mu$ with respect
to $\theta$ \cite{DemZei,DeuStr,DupEll}.  
For $\mu \in \X$ the interaction representation function is defined by
\[
\tilde{H}(\mu) \doteq \frac{1}{2} \int_{\Lambda \times \Lambda}
V(x,x') \, \mu(dx) \, \mu(dx').
\]
The approximation property (\ref{eqn:repr-fn}) is easily verified.

\skp
(b) {\bf Miller-Robert Model.}  This model of the two-dimensional
Euler equations is analyzed in
\cite{BouEllTur}, which explains in detail the physical
background.  For simplicity, let the flow domain be $T^2$, the unit torus
$[0,1) \times [0,1)$ with periodic boundary conditions.
For each $\nin$ let $\LL_n$ be a uniform lattice of 
$a_n \doteq 2^{2n}$ sites $t$ in 
$T^2$.  The intersite
spacing in each coordinate direction is $2^{-n}$. 
Each such lattice of $a_n$ sites induces a 
dyadic partition of $T^2$ into $a_n$ squares called microcells,
each having area $1/a_n$.  For each $s\in\LL_n$ we denote 
by $M(s)$ the unique
microcell containing the site $s$ in its lower left corner.
The configuration spaces for the Miller-Robert model are
$\Omega_n \doteq \Y^{a_n}$, where $\Y$ is a given compact subset of
$\R$.   Microstates are denoted by $\zeta = \{\zeta(s), s \in \LL_n\}$.  
Let $\rho$ be a probability measure
on $\R$ with support $\Y$.  
$P_n$ equals the product measure on $\Omega_n$ with
identical one-dimensional marginals $\rho$.  

There are two classes of interaction functions, the Hamiltonian
and the generalized enstrophies.  
For $\zeta \in \Omega_n$ the Hamiltonian is defined by
\[
\label{eqn:mrham} 
H_{n,1}(\zeta) \doteq \frac{1}{2 n^2} \sum_{s,s' \in \LL} 
g_n(s - s') \zeta(s) \zeta(s') \, , 
\] 
where $g_n(s-s')$ is a certain bounded continuous
approximation to the lattice Green's function 
\[ 
\label{eqn:lattice-green-fn}
g(s-s') \doteq \sum_{0 \not= \xi \in \ZZ^2}
   | 2 \pi \xi |^{-2} \exp[2\pi i \lan \xi,s-s'\ran].
\]
Fix $\alpha \in \N$.  For $i =2,\ldots,\alpha+1$
the generalized enstrophies are defined by 
\[ 
\label{eqn:genenstrophy} 
H_{n,i}(\zeta) \doteq \frac{1}{n} \sum_{s \in \LL_n} a_i(\zeta(s)), 
\]
where the $a_i$ are continuous functions mapping $\Y$ into $\R$.  

Hypotheses \ref{hyp:general} are verified in \cite{BouEllTur},
to which the reader is referred for details. 
Let $\theta$ denote Lebesgue measure on $T^2$.  
The hidden space is the space $\X$ of probability
measures $\mu(dx \times dy)$ on $T^2 \times \Y$ with first marginal
$\theta$.  The hidden process is the sequence of measures
\[
\label{eqn:mrhidden}
Y_n(dx \times dy) = Y_n(\zeta,dx \times dy) \doteq
\theta(dx) \otimes \sum_{s \in \LL_n} 1_{M(s)}(x) \, 
\delta_{\zeta(s)}(dy).
\]
For $\mu \in \X$ the Hamiltonian interaction
function is given by
\[
\label{eqn:mrhamrep}
\tilde{H}_1(\mu) \doteq \frac{1}{2} \int_{(T^2 \times \Y)^2} g(x-x') y y'
\, \mu(dx \times dy) \, \mu(dx' \times dy')
\]
while for $i =2,\ldots,\alpha+1$ the interaction functions for
the generalized enstrophies are given by
\[
\tilde{H}_i(\mu) \doteq \int_{T^2 \times \Y} a_i(y) \, \mu(dx \times dy).
\]
For $i=1$ one verifies (\ref{eqn:repr-fn}) by a detailed Fourier
analysis. 
For $i=2,\ldots,\alpha+1$ (\ref{eqn:repr-fn}) 
is easily verified to hold with no error term.

Given $\nin$ and an even integer $q < 2n$, 
we consider a dyadic partition of the lattice $\LL_n$ into
$2^q$ blocks, each block containing $a_n/2^q$ lattice sites. 
In correspondence with this partition we have a dyadic partition
$\{D_{q,k}, k =1,\ldots,2^q\}$ of $T^2$ into macrocells.  Each macrocell
is the union of $a_n/2^q$ microcells $M(s)$.  
The large deviation principle for $Y_n$ with respect to $P_n$ is verified
by comparing $Y_n$ with the two-component process
\[
\label{eqn:mrtwocomp}
W_{n,q}(dx \times dy) = W_{n,q}(\zeta,dx \times dy) \doteq
\theta(dx) \otimes \sum_{k =1}^{2^q} 1_{D_{q,k}}(x) \, L_{n,q,k}(\zeta,dy), 
\]
where $L_{n,q,k}$ denotes the empirical measure $\frac{1}{a_n/2^q} \sum_{s \in D_{q,k}} \delta_{\zeta(s)}(dy)$.  Through these empirical measures, $W_{n,q}$ introduces an averaging over the intermediate scale of the macrocells and 
thus corresponds to a coarse graining of the vorticity field.  
Using Sanov's Theorem, one verifies
that as $n \goto \infty$, $q \goto \infty$, $W_{n,q}$ satisfies the two-parameter LDP on $\X$ with rate function the relative
entropy $R(\mu|\theta \times \rho)$ of $\mu(dx \times dy)$ with respect to the 
product measure $\theta(dx) \times \rho(dy)$ \cite[\S 5]{BouEllTur}.  An 
approximation result relating $Y_n$ and $W_{n,q}$ then allows 
one to prove that $Y_n$
satisfies the LDP on $\X$ with the same rate function.  

\skp
(c)  {\bf Quasi-geostrophic potential vorticity model.}  This model of the 
quasi-geostrophic potential vorticity equation, described in detail 
in \cite{DibMajTur} and \cite{EllHavTur}, incorporates
the geophysical terms associated with the Coriolis effect,
the deformation of an upper free surface, and bottom topography.
The large deviation analysis of the model is carried out in \cite{EllHavTur}.

\skp
(d)  {\bf Dispersive wave model for the nonlinear Schr\"{o}dinger
equation.}  This model is defined in \cite{JorJos,JorTurZir}, to which the 
reader is referred for details.  The hidden process is a Gaussian 
process taking values in $L^2[0,1]$ and satisfying the LDP with respect to
the prior distribution that is proved in 
\cite{EllRos}.  The large deviation analysis of this model is the
subject of \cite{EllJorTur}.
\ink
\eexa 

We now return to the general model.  Its large deviation analysis 
with respect to the canonical ensemble is summarized
in the next theorem.  Part (a) states a variational formula
for $\vphi(\beta)$, and part (b)
gives the LDP for the hidden process $Y_n$
with respect to the sequence of Gibbs measures $P_{n,\beta}$.  
Part (c) describes
the set ${\cal E}_{\beta}$ consisting of
points at which the rate function in part (b)
attains its minimum of 0.  Part (d) gives a concentration
property of $\E_\beta$.  As we point out after the statement of the theorem,
${\cal E}_{\beta}$ can be identified with the
set of equilibrium macrostates of the statistical mechanical model.
The mathematical tractability of the canonical ensemble is
reflected in the simplicity of the proof of Theorem \ref{thm:canon}.

\begin{thm} \per
\label{thm:canon}
We assume Hypotheses {\em \ref{hyp:prob}} and {\em \ref{hyp:general}}.
For $\beta \in \R^\sigma$
the following conclusions hold.

{\em (a)} $\vphi(\beta) \doteq - \lim_{n \goto \infty} 
\frac{1}{a_n} \log Z_n(a_n \beta)$ exists and is given by
\be
\label{eqn:vphivar}
\vphi(\beta) = \inf_{x \in \X}
\{\lan \beta, \tilde{H}(x) \ran + I(x)\};
\ee
$\vphi(\beta)$ is a finite, concave, continuous function on $\R^\sigma$.

{\em (b)}  With respect to 
$P_{n,a_n \beta}$, $Y_n$
satisfies the LDP on $\X$ with rate function
\[
I_{\beta}(x) \doteq I(x) + \lan \beta , \tilde{H}(x) \ran
- \inf_{y \in \X}\{I(y) + \lan \beta,\tilde{H}(y) \ran \}
= I(x) + \lan \beta , \tilde{H}(x) \ran - \varphi(\beta).
\] 

{\em (c)} The set ${\cal E}_{\beta } \doteq
\{x \in \X: I_{\beta}(x) = 0\}$ is a nonempty,
compact subset of $\X$.  A point $\bar{x}$ lies in $\E_\beta$ if and only if
\[
I(\bar{x}) + \lan \beta, \tilde{H}(\bar{x}) \ran = \inf_{y \in \X}
\{I(y) + \lan \beta,\tilde{H}(y) \ran \} = \vphi(\beta);
\]
equivalently, if and only if $\bar{x}$ solves the following unconstrained
minimization problem:
\[
\mbox{ minimize } (I(x) + \lan \beta, \tilde{H}(x) \ran) \mbox{ over }
x \in \X.
\]

{\em (d)} If $A$ is any Borel subset
of $\X$ whose closure $\bar{A}$ satisfies $\bar{A} \cap {\cal E}_{\beta} = \emptyset$, then $I_{\beta}(\bar{A}) > 0$ and
for some $C < \infty$ 
\[
P_{n,a_n \beta }\{Y_n \in A\} \leq C \exp[-a_n I_{\beta}(\bar{A})/2] 
\goto 0 \: \mbox{ as } n \goto \infty.
\]
\end{thm}

\noi
{\bf Proof.}  (a) Since $Y_n$ satisfies the LDP with respect to $P_n$, 
$Y_n$ satisfies the Laplace principle with respect to $P_n$ with the
same rate function $I$.  Hence by the approximation 
property (\ref{eqn:repr-fn}) and the boundedness and continuity of 
the function mapping $x \mapsto \lan \beta, \tilde{H}(x) \ran$, 
\begin{eqnarray*}
\vphi(\beta) & = & -\lim_{n \goto \infty} \frac{1}{a_n} \log Z_n(a_n \beta) \\
& = & -\lim_{n \goto \infty} \frac{1}{a_n} \log
\int_{\Omega_n} \exp[-a_n \lan \beta,H_n \ran] \, dP_n \\
& = & - \lim_{n \goto \infty} \frac{1}{a_n} \log \int_{\Omega_n}
\exp[- a_n \lan \beta, \tilde{H}(Y_n) \ran ] \, dP_n \\
& = & \inf_{x \in \X} \{\lan \beta, \tilde{H}(x) \ran + I(x)\}.
\end{eqnarray*}
This formula exhibits $\vphi$ as a finite, concave function on $\R^\sigma$,
which is therefore continuous on $\R^\sigma$.

(b) $I_\beta$ is a rate function since $I$ is a rate function
and the function mapping $x \mapsto \lan \beta, \tilde{H}(x) \ran$ is
bounded and continuous.  In order to prove that 
with respect to $P_{n,a_n \beta}$ 
$\, Y_n$ satisfies the LDP with rate function $I_\beta$, it suffices
to prove that with respect to $P_{n,a_n \beta}$
$\, Y_n$ satisfies the Laplace principle with rate function $I_\beta$.
This is an immediate consequence of (\ref{eqn:repr-fn}) and part (a);
for details, see the proof of part (b) of Theorem 3.1 in \cite{BouEllTur}.

(c) ${\cal E}_\beta$ is a nonempty, compact subset of $\X$ because
$I_\beta$ is a rate function.  The equivalent characterizations of $\bar{x} \in \E_\beta$ follow from the definition of $I_\beta$.

(d) If $\bar{A} \cap {\cal E}_{\beta} = \emptyset$, 
then for each $x \in A$ we have 
$I_{\beta}(x)>0$.  Since $I_{\beta}$ is a rate function,
it follows that $I_{\beta}(\bar{A}) > 0$.
The large deviation upper bound in part (b) yields the display
in part (d) for some $C < \infty$. The proof of the theorem
is complete.  \ink

\skp

Part (d) of Theorem \ref{thm:canon}
can be regarded as a concentration property
of the $P_{n,a_n \beta}$-distributions of $Y_n$.  This
property justifies calling ${\cal E}_{\beta}$ the set
of equilibrium macrostates with respect to $P_{n,a_n \beta}\{Y_n \in dx\}$
or, for short, as the set of canonical equilibrium macrostates.

The next theorem further justifies the designation of ${\cal E}_{\beta}$ 
as the set of canonical equilibrium macrostates by relating weak
limits of subsequences of $P_{n,a_n \beta}\{Y_n \in \cdot\}$
to $\ebeta$.  For example, if one knows that ${\cal E}_{\beta}$ 
consists of a unique point $\tilde{x}$, then it follows
that the entire sequence $P_{n,a_n \beta}\{Y_n \in \cdot\}$ 
converges weakly to $\delta_{\tilde{x}}$.  This situation
corresponds to the absence of a phase transition.  For specific models,
 more detailed
information about weak limits of subsequences of $P_{n,a_n \beta}$ 
have been obtained by a number of authors including
\cite{CagLioMarPul2,EyiSpo,Kie,MesSpo}.

\begin{thm} \per
\label{thm:weak-limits}
We assume Hypotheses {\em \ref{hyp:prob}} and {\em \ref{hyp:general}}.
For $\beta \in \R^\sigma$, any subsequence
of $P_{n,a_n \beta}\{Y_n \in \cdot\}$ has 
a subsubsequence converging weakly to a probability measure
$\Pi_{\beta}$ on $\X$ that is concentrated
on ${\cal E}_{\beta} \doteq \{x \in \X : I_\beta(x) = 0\}$; i.e., $\Pi_{\beta}\{({\cal E}_{\beta})^c\} =0$.
If ${\cal E}_{\beta}$
consists of a unique point $\tilde{x}$, then the entire 
sequence $P_{n,a_n \beta}\{Y_n \in \cdot\}$ 
converges weakly to $\delta_{\tilde{x}}$. 
\end{thm}

\noi
{\bf Proof.} 
Define $a^* \doteq \min_{\nin} a_n > 0$.  
As shown in the proof of Lemma 2.6 in \cite{LynSet},
the large deviation upper bound
given in part (b) of Theorem \ref{thm:canon} implies that
for each $M \in (0,\infty)$ there
exists a compact subset $K$ of $\X$ such that for all $\nin$
\[
P_{n,a_n \beta}\{Y_n \in K^c\} \leq
\frac{e^{-a_n M}}{1 - e^{-M}} \leq \frac{e^{-a^* M}}{1-e^{-M}} \; .
\]
It follows that the sequence $P_{n,a_n \beta}\{Y_n \in \cdot\}$ 
is tight and therefore that any subsequence has a subsubsequence
$P_{n',a_{n'} \beta}\{Y_{n'} \in \cdot\}$ converging
weakly as $n' \goto \infty$ to a probability measure $\Pi_{\beta}$ 
on $\X$ [Prohorov's Theorem].  In order to show that 
$\Pi_{\beta}$ is concentrated on ${\cal E}_{\beta}$, 
we write the open set $({\cal E}_{\beta})^c$ 
as a union of countably many open balls $V_j$
such that the closure $\bar{V}_j$ of each $V_j$ has empty intersection
with ${\cal E}_{\beta}$.  By part (c) of Theorem
\ref{thm:canon} $P_{n',a_{n'} \beta}\{Y_{n'} \in V_j\} \goto 0$
as $n' \goto \infty$, and so
\[
0 = \liminf_{n' \goto \infty} P_{n',a_{n'}\beta}\{Y_{n'} \in V_j\}
\geq \Pi_{\beta}\{V_j\}.
\]
It follows that $\Pi_{\beta}\{V_j\} = 0$ and thus
that $\Pi_{\beta}\{({\cal E}_{\beta})^c\} = 0$, as claimed. 

Now assume that ${\cal E}_{\beta} = \{\tilde{x}\}$.
Then the only probability measure on $\X$ that is concentrated on
${\cal E}_{\beta}$ is $\delta_{\tilde{x}}$.  
Since by the first part of the proof any subsequence of 
$P_{n,a_n \beta}\{Y_n \in \cdot\}$ has a subsubsequence
converging weakly to $\delta_{\tilde{x}}$, it follows that the entire
sequence $P_{n,a_n \beta}\{Y_n \in \cdot\}$
converges weakly to $\delta_{\tilde{x}}$.  This completes the proof.  \ink

\skp
In the next section we consider the LDP for $Y_n$
when conditioning is present.

\skp
\skp
\section{Microcanonical Ensemble: LDP and Equilibrium Macrostates}
\beginsec

As in the preceding section, we consider 
models defined in
terms of a sequence of 
interaction functions $\{H_{n,i}, i=1\ldots,\sigma\}$, which are
bounded measurable functions mapping $\Omega_n$ into $\R$. 
In general, the interaction 
functions represent conserved quantities with respect to some dynamics
that underlies the model.  For suitable
values of $(u_1,\ldots,u_\sigma) \in \R^\sigma$ the ideal way to define
the microcanonical ensemble
is to condition the probability measure $P_n$ on the set
$\{  H_{n,1} = u_1,\ldots,  H_{n,\sigma} = u_\sigma\}$.  
However, in order to avoid problems concerning 
the existence of regular conditional probability distributions,
we shall condition $P_n$ on 
$\{H_{n,1} \in [u_1-r,u_1+r],\ldots, 
H_{n,\sigma} \in [u_\sigma - r,u_\sigma+r]\}$,
where $r \in (0,1)$.
These conditioned measures, given in (\ref{eqn:definemicro}), 
define the microcanonical ensemble.
Theorem \ref{thm:microcan}
proves the LDP for the 
distributions of $Y_n$ with respect to the microcanonical
ensemble in the double limit 
obtained by sending first $n \goto \infty$ and then $r \goto 0$.
We then define, in terms of the rate function in this LDP,
the set of microcanonical equilibrium macrostates and derive some of its properties.

For $u = (u_1,\ldots,u_\sigma) \in \R^\sigma$ a key role in the large deviation
analysis of the microcanonical ensemble is played by 
\be 
\label{eqn:j}
J(u) \doteq \inf\{I(x) : x \in \X, \tilde{H}(x) = u\} .
\ee
In terms of $J$ the canonical free energy $\vphi(\beta)$, 
given in part (a) of Theorem \ref{thm:canon} by
\[
\vphi(\beta) = \inf_{x \in \X} \{
\lan \beta , \tilde{H}(x) \ran + I(x) \},
\]
can be rewritten as
\begin{eqnarray*}
\vphi(\beta) & = &
\inf_{u \in \R^\sigma} \left\{ \inf\{
\lan \beta , \tilde{H}(x) \ran + I(x) : x \in \X,
\tilde{H}(x) = u \}\right\}  \\
& = &   \inf_{u \in \R^\sigma} \{\lan \beta , u \ran
+ J(u)\}. \nonumber
\end{eqnarray*}
Introducing the microcanonical entropy
\be
\label{eqn:s(z)}
s(u) \doteq -J(u) = - \inf\{I(x) : x \in \X, \tilde{H}(x) = u\},
\ee
we have
\be
\label{eqn:fenchel}
\vphi(\beta) = \inf_{u \in \R^\sigma} \{\lan \beta,u \ran - s(u)\}.
\ee
This formula expresses $\vphi$
as the Legendre-Fenchel transform of $s$.
The microcanonical entropy 
will play a central role in the results on equivalence 
and nonequivalence of the canonical and microcanonical
ensembles to be presented in Section 4. 

The function $J$ plays other roles in the theory.  Since each
$\tilde{H}_i$ is a bounded continuous function mapping $\X$
into $\R$ and since with respect to $P_n$ $\, Y_n$ satisfies the LDP 
on $\X$ with 
rate function $I$, it follows from the contraction principle that
with respect to $P_n$ 
$\, \tilde{H}(Y_n) = 
(\tilde{H}_1(Y_n),\ldots,\tilde{H}_\sigma(Y_n))$ satisfies the
LDP on $\R^\sigma$ with rate function $J$ \cite[Thm.\ 4.2.1]{DemZei}.
When expressed in terms of the equivalent Laplace principle,
this means that for any bounded continuous function $g$ mapping
$\R^\sigma$ into $\R$
\[
\lim_{n \goto \infty} \frac{1}{a_n} \log \int_{\Omega_n} 
\exp[a_n \, g(\tilde{H}(Y_n))] \, dP_n 
= \sup_{u \in \R^\sigma} \{g(u) - J(u)\}.
\]
Because of the approximation property (\ref{eqn:repr-fn}),
this readily extends to the Laplace principle on $\R^\sigma$, 
and thus the LDP on $\R^\sigma$, for 
$H_n \doteq (H_{n,1},\ldots, H_{n,\sigma})$.  

In part (a) of the next proposition we record the LDP's just discussed and
two properties of the microcanonical entropy.
When applied to the regularized point vortex model, the LDP for
the $P_n$-distributions of $H_n$ generalizes
the large deviation estimates obtained in \cite[Thm.\ 2.1]{EyiSpo}.
In parts (b) and (c) of the proposition some related
facts needed later in this section are given.  We define $\, \mbox{dom} \, J \,$ to be the set of $u \in \R^\sigma$ for which $J(u) < \infty$.  For $r \in (0,1)$ and $u \in \mbox{dom} \, J$, we also define 
\[
\{u\}^{(r)} \doteq [u_1-r,u_1+r] \times \cdots \times [u_\sigma - r,u_\sigma +r].
\]
Part (b) is a consequence of the LDP for $H_n$ given in 
part (a) and of the bound
$J(\mbox{int}(\{u\}^{(r)})) \leq J(u)$.  
Part (c) follows from the lower semicontinuity
of $J$ and from part (b).

\begin{prop} \per
\label{prop:ldpj}  We assume Hypotheses {\em \ref{hyp:prob}} and
{\em \ref{hyp:general}}.  The following conclusions hold.

{\em (a)}  With respect to $P_n$, the sequences $\tilde{H}(Y_n)$
and $ H_n$ satisfy the LDP on $\R^\sigma$
with rate function $J$.  Hence $s \doteq -J$ is nonpositive
and upper semicontinuous.

{\em (b)}  For $u \in \mbox{{\em dom}} \, J$
and any $r \in (0,1)$ 
\begin{eqnarray*}
-J(u) & \leq & \liminf_{n \goto \infty} 
\frac{1}{a_n} \log P_n\{H_n \in \{u\}^{(r)}\} \\
& \leq & \limsup_{n \goto \infty} \frac{1}{a_n} \log P_n\{H_n \in \{u\}^{(r)}\} 
\; \leq \; -J(\{u\}^{(r)}).
\end{eqnarray*}

{\em (c)}  As $r \goto 0$, $J(\{u\}^{(r)}) \nearrow J(u)$.  Hence
\[
\lim_{r \goto 0} \, \lim_{n \goto \infty} \frac{1}{a_n} \log P_n\{H_n \in \{u\}^{(r)}\} = -J(u).
\]
\end{prop}

\skp
The main theorem of this section is the LDP
for $Y_n$ with respect to the microcanonical ensemble,
given in Theorem \ref{thm:microcan}.
For $A \in {\cal F}_n$ this ensemble is defined by the conditioned measures
\be
\label{eqn:definemicro}
P_{n}^{u,r}\{A\} \doteq P_n\{A \, |   H_n \in \{u\}^{(r)}\},
\ee
where $u \in \mbox{dom} \, J$ and $r \in (0,1)$.
For all sufficiently large $n$ it follows from part (b) of Proposition \ref{prop:ldpj} that $P_n\{H_n \in \{u\}^{(r)}\}
> 0$ and hence that $P_n^{u,r}$ is well defined.

\begin{thm} \per
\label{thm:microcan}
Take $u \in \mbox{{\em dom}} \, J$
and assume Hypotheses {\em \ref{hyp:prob}} and {\em \ref{hyp:general}}.
With respect to the conditioned measures 
$P_{n}^{u ,r}$, $Y_n$ satisfies the LDP
on $\X$, in the double limit $n \goto \infty$ and
$r \goto 0$, with rate function
\[
I^{u}(x) \doteq 
\left\{ \begin{array}{ll} I(x) -
J(u) &  \trm{if }\; \tilde{H}(x) = u, \\ 
                                \infty & \trm{otherwise.}   
                               \end{array}  \right.
\]  
That is, for any closed subset $F$ of $\X$
\be 
\label{eqn:microldupper}
\lim_{r \goto 0} \, \limsup_{n \goto \infty}
\frac{1}{a_n} \log P_n^{u ,r}\{Y_n \in F\} 
\leq - I^{u }(F)  
\ee
and for any open subset $G$ of $\X$
\be 
\label{eqn:microldlower}
\lim_{r \goto 0} \, \liminf_{n \goto \infty}
\frac{1}{a_n} \log P_n^{u ,r}\{Y_n \in G \}
\geq - I^{u }(G).   
\ee 
\end{thm}

\skp
We first prove that $I^u$ defines a rate function.
Clearly $I^u$ is nonnegative.
For $u \in \mbox{dom} \, J$
and $M < \infty$
\[
\{x \in \X : I^u(x) \leq M\} = \{x \in \X : I(x) \leq M + J(u)\}
\cap \tilde{H}^{-1}(\{u\}).
\]
Since $J(u) < \infty$, $I$ has compact level sets, 
and $\tilde{H}^{-1}(\{u\})$ is closed,
it follows that $I^u$ has compact level sets.

Concerning the large deviation bounds in Theorem \ref{thm:microcan}, 
we offer two proofs.
The first is preferred because it is close to the heuristic sketch 
of the LDP given in the introduction.  
Throughout the two proofs we fix $u \in \mbox{dom} \, J$.

The first proof of the large deviation upper bound 
actually derives a stronger inequality.   
Namely, for all sufficiently small 
$r \in (0,1)$ and any closed subset $F$ of $\X$
\be
\label{eqn:thisisupper}
\limsup_{n \goto \infty} \frac{1}{a_n} \log P_n^{u,r}\{Y_n \in F\} 
\leq -I^u(F).
\ee
For any $x \in \X$ and $\alpha > 0$ we denote by $\bar{B}(x,\alpha)$ and
$B(x,\alpha)$ the closed ball and the open ball in $\X$ with center $x$
and radius $\alpha$.  Let $\delta > 0$ be given.
Since $I$ is lower semicontinuous, for any $x \in \X$ and all sufficiently
small $\alpha > 0$ we have $I(\bar{B}(x,\alpha)) \geq I(x) - \delta$.
Now take any $x \in \X$ such that $\tilde{H}(x) = u$.
For any 
$r \in (0,1)$ and all sufficiently small $\alpha$ 
the large deviation upper bound for $Y_n$ with
respect to $P_n$ and part (b) of Proposition \ref{prop:ldpj} yield
\begin{eqnarray}
\label{eqn:hotsoup}
\lefteqn{
\limsup_{n \goto \infty} \frac{1}{a_n} \log P_n^{u,r}\{Y_n \in \bar{B}(x,\alpha)\}}
 \\ && \leq \limsup_{n \goto \infty} \frac{1}{a_n} 
\log P_n\{\{Y_n \in \bar{B}(x,\alpha)\} \cap \{H_n \in \{u\}^{(r)}\}\} \nonumber \\ && \hspace{.35in} - \liminf_{n \goto \infty} \frac{1}{a_n} \log P_n\{H_n \in \{u\}^{(r)}\} \nonumber \\
&& \leq \limsup_{n \goto \infty} \frac{1}{a_n} 
\log P_n\{Y_n \in \bar{B}(x,\alpha) \} \nonumber \\
&& \hspace{.35in} - \liminf_{n \goto \infty} \frac{1}{a_n} \log P_n\{H_n \in \{u\}^{(r)}\} \nonumber \\
&& \leq -I(\bar{B}(x,\alpha)) + J(u) \nonumber \\
&& \leq -I(x) + J(u) + \delta \nonumber \\
&& = -I^u(x)+ \delta.  \nonumber
\end{eqnarray}

Now take any $x \in \X$ such that $\tilde{H}(x) \not = u$.
Thus $I^u(x) = \infty$, and there exists $t \in (0,1)$ such that
$\tilde{H}(x) \not \in \{u\}^{(t)}$.
By the approximation property (\ref{eqn:repr-fn}) and the 
continuity of $\tilde{H}$, 
for any $r \in (0,t)$, 
all sufficiently small $\alpha > 0$,
and all sufficiently large $n$ we have 
\[
\{Y_n \in \bar{B}(x,\alpha)\} \cap \{H_n \in \{u\}^{(r)}\} 
 \subset \{Y_n \in \bar{B}(x,\alpha)\} \cap 
\{\tilde{H}(Y_n) \in \{u\}^{(t)} \} = \emptyset.
\]
Hence for such $r$ and $\alpha$
\begin{eqnarray*}
\lefteqn{
\limsup_{n \goto \infty} \frac{1}{a_n} \log P_n^{u,r}\{Y_n \in \bar{B}(x,\alpha)\}}
 \\ && \leq \limsup_{n \goto \infty} \frac{1}{a_n} 
\log P_n\{\{Y_n \in \bar{B}(x,\alpha)\} \cap \{H_n \in \{u\}^{(r)}\}\} \\
&& \hspace{.35in} - \liminf_{n \goto \infty} \frac{1}{a_n} \log P_n\{H_n \in \{u\}^{(r)}\} \\  && = -\infty = -I^u(x).
\end{eqnarray*}

We have proved that for any $x \in \X$,
all sufficiently small $r \in (0,1)$, and all sufficiently small $\alpha > 0$
\[
\limsup_{n \goto \infty} \frac{1}{a_n} \log P_n^{u,r}\{Y_n \in \bar{B}(x,\alpha)\} \leq -I^u(x) + \delta.
\]
Let $F$ be a compact subset of $\X$.  We can cover $F$ with finitely
many closed balls $\bar{B}(x_i,\alpha_i)$ with $x_i \in F$ and $\alpha_i > 0$
so small that the last display is valid for $x = x_i$,
all sufficiently small $r \in (0,1)$, and $\alpha = \alpha_i$.
It follows that for all sufficiently small $r \in (0,1)$
\[
\limsup_{n \goto \infty} \frac{1}{a_n} \log P_n^{u,r}\{Y_n \in F\}
\leq - \min_i I^u(x_i) + \delta \leq -I(F) + \delta.
\]
Sending $\delta \goto 0$ yields the upper bound (\ref{eqn:thisisupper}).
Finally, for any closed set $F$ the upper bound (\ref{eqn:thisisupper}) 
is a consequence of the following uniform
exponential tightness estimate.

\begin{lemma} \per
\label{lem:exptight}
Fix $u \in \mbox{{\em dom}} \, J$.
Then for all sufficiently large
$M \in (0,\infty)$ there exists a compact subset $D$ of $\X$
such that for every $r \in (0,1)$
\[
\limsup_{n \goto \infty} \frac{1}{a_n} \log P_n^{u,r}\{Y_n \in D^c\} 
\leq -M.
\]
\end{lemma}

\noi
{\bf Proof.}  Given $u \in \mbox{dom} \, J$, we take $M > J(u)$.
As shown in the proof of Lemma 2.6 in \cite{LynSet},
the large deviation upper bound satisfied by $Y_n$ with respect to $P_n$
implies that there exists a compact
subset $D$ of $\X$ such that
\[
\limsup_{n \goto \infty} \frac{1}{a_n} \log P_n\{Y_n \in D^c\} \leq -
2M.
\]
Since for every $r \in (0,1)$ 
\[
P_n^{u,r}\{Y_n \in D^c\} \leq 
\frac{P_n\{Y_n \in D^c\}}
{P_n\{H_n \in \{u\}^{(r)}\}},
\]
it follows from part (b) of Proposition \ref{prop:ldpj} that
\begin{eqnarray*}
\lefteqn{
\limsup_{n \goto \infty} \frac{1}{a_n} \log P_n^{u,r}\{Y_n \in D^c\} } 
\\
&& \leq  \limsup_{n \goto \infty} \frac{1}{a_n} \log P_n\{Y_n \in 
D^c\} - \liminf_{n \goto \infty} \frac{1}{a_n} \log 
P_n\{H_n \in 
 \{u\}^{(r)}\} \\ && \leq -2M + J(u) \leq -M.
\end{eqnarray*}
This completes the proof.  \ink
 
\skp
We next prove the large deviation lower bound in Theorem \ref{thm:microcan}
by showing that for any fixed $r \in (0,1)$ and any open subset $G$ of $\X$
\be
\label{eqn:thisislower}
\liminf_{n \goto \infty} \frac{1}{a_n} \log P_n^{u,r}\{Y_n \in G\} 
\geq -I^u(G) + J(\{u\}^{(r)}) - J(u).
\ee
Sending $r \goto 0$ and using part (c) of Proposition \ref{prop:ldpj}
yields the large deviation lower bound in Theorem \ref{thm:microcan}.

Let $x$ be any point in $G$ such that $\tilde{H}(x) = u$. 
By the approximation property (\ref{eqn:repr-fn}) and the 
continuity of $\tilde{H}$, for any number $r^-$
satisfying $0 < r^- < r$ and all sufficiently large $n$, we can choose 
$\alpha > 0$ to be so small that $B(x,\alpha) \subset G$ and 
\begin{eqnarray*}
\{Y_n \in {B}(x,\alpha)\} \cap \{H_n \in \{u\}^{(r)} \}
& \supset & \{Y_n \in {B}(x,\alpha)\} \cap 
\{\tilde{H}(Y_n) \in \{u\}^{(r^-)} \} \\ & = & \{Y_n \in {B}(x,\alpha)\}.
\end{eqnarray*}
Hence for such $\alpha$, the large deviation lower bound for $Y_n$ with
respect to $P_n$ and part (b) of Proposition \ref{prop:ldpj} yield 
\begin{eqnarray*}
\lefteqn{
\liminf_{n \goto \infty} \frac{1}{a_n} 
\log P_n^{u,r}\{Y_n \in G\} }
\\ && \geq \liminf_{n \goto \infty} \frac{1}{a_n} 
\log P_n^{u,r}\{Y_n \in {B}(x,\alpha)\}
 \\ && \geq \liminf_{n \goto \infty} \frac{1}{a_n} 
\log P_n\{\{Y_n \in {B}(x,\alpha)\} \cap \{H_n \in \{u\}^{(r)}\}\} \\
&& \hspace{.35in} - \limsup_{n \goto \infty} \frac{1}{a_n} \log P_n\{H_n \in \{u\}^{(r)}\} \\
&& \geq \liminf_{n \goto \infty} \frac{1}{a_n} 
\log P_n\{Y_n \in {B}(x,\alpha)\} \\
&& \hspace{.35in} - \limsup_{n \goto \infty} \frac{1}{a_n} \log P_n\{H_n \in \{u\}^{(r)}\} \\
&& \geq -I({B}(x,\alpha)) + J(\{u\}^{(r)}) \\ 
&& \geq -I(x) + J(\{u\}^{(r)}) \\
&& = -I^u(x) + J(\{u\}^{(r)}) - J(u).
\end{eqnarray*}
Now take any $x \in \X$ such that $\tilde{H}(x) \not = u$.   Since
$I^u(x) = \infty$, it follows that
\[
\liminf_{n \goto \infty} \frac{1}{a_n} 
\log P_n^{u,r}\{Y_n \in G\} \geq -\infty = 
-I^u(x) + J(\{u\}^{(r)}) - J(u).
\]
We have thus obtained the same lower bound for all $x \in G$.
We conclude that 
\begin{eqnarray*}
\liminf_{n \goto \infty} \frac{1}{a_n} 
\log P_n^{u,r}\{Y_n \in G\} & \geq & 
\sup_{x \in G}\{-I^u(x)\} + J(\{u\}^{(r)}) - J(u) \\
& = & -I^u(G) + J(\{u\}^{(r)}) - J(u).
\end{eqnarray*}
This completes the proof of the large deviation lower bound
(\ref{eqn:thisislower}).  The proof of Theorem \ref{thm:microcan}
is done.

The second proof of the large deviation bounds in Theorem \ref{thm:microcan}
uses the following alternate representation for the rate function:
\[
I^u(x) = I(\{x\} \cap \tilde{H}^{-1}(\{u\})).
\]
Let $F$ be any closed subset of $\X$.  We choose
$\psi$ to be any function mapping 
$(0,1)$ onto $(0,1)$ with the properties that $\psi(r) > r$  for all $r \in (0,1)$ and $\lim_{r \goto 0} \psi(r) =0$.
Clearly, as $r \downarrow 0$, $\{u\}^{(\psi(r))} \downarrow \{u\}$.
We need the limit 
\[
\lim_{r \rightarrow 0} I(F \cap \tilde{H}^{-1}(\{u\}^{(\psi(r))})) \\
= I(F \cap \tilde{H}^{-1}(u)),
\]
which follows from routine calculations using the continuity of $\tilde{H}$
and the fact that $I^u$ is a rate function.  The proof of this limit is omitted.  The rest of the proof of the large deviation upper bound is straightforward.   
By the approximation property (\ref{eqn:repr-fn}) and the 
continuity of $\tilde{H}$, 
for any $r \in (0,1)$ and all sufficiently large $n$
\[
P_n\{\{Y_n \in F \} \cap \{H_n \in \{u\}^{(r)}\} \}
 \leq P_n\{\{Y_n \in F\} \cap 
\{\tilde{H}(Y_n) \in \{u\}^{(\psi(r))} \}\}.
\]
Then the large deviation upper
bound for $Y_n$ with respect to $P_n$ and part (c) of Proposition \ref{prop:ldpj} yield
\begin{eqnarray*}
 \lefteqn{ \lim_{r \rightarrow 0} \limsup_{n \rightarrow \infty} \frac{1}{a_n}
  \log P_n^{u,r}\{ Y_n \in F\} } \\
 && \leq  \lim_{r \rightarrow 0} \limsup_{n \rightarrow \infty} \frac{1}{a_n} 
  \log P_n\{ Y_n \in [F \cap
  \tilde{H}^{-1}(\{u\}^{(\psi(r))})]\}
   \\  && \hspace{.3in} - \ 
\lim_{r \goto 0} \, \liminf_{n \goto \infty} \frac{1}{a_n} \log P_n\{H_n \in \{u\}^{(r)}\} \\
 && \leq  -\lim_{r \rightarrow 0} 
I(F \cap \tilde{H}^{-1}(\{u\}^{(\psi(r))})) + J(u) \\
 && = -I(F \cap \tilde{H}^{-1}(u)) + J(u) \\
 && = -I^u(F).
\end{eqnarray*}
This is the large deviation upper bound (\ref{eqn:microldupper}).

Now let $G$ be any open
subset of $\X$.  Again by the approximation property 
(\ref{eqn:repr-fn}) and the 
continuity of $\tilde{H}$, 
for any number $r^-$ satisfying $0 < r^- < r$ and all sufficiently large $n$
\begin{eqnarray*}
\lefteqn{
P_n\{\{Y_n \in G \} \cap \{H_n \in \{u\}^{(r)}\} \}} \\
&&  \geq P_n\{\{Y_n \in G\} \cap 
\{\tilde{H}(Y_n) \in \{u\}^{(r^-)} \}\} \\
&& \geq P_n\{Y_n \in G \cap \tilde{H}^{-1}(\mbox{int}\{u\}^{(r^-)})\}.
\end{eqnarray*}
The large deviation lower bound for $Y_n$ with
respect to $P_n$ and part (c) of Proposition \ref{prop:ldpj} yield
\begin{eqnarray*}
 \lefteqn{ \lim_{r \rightarrow 0} \liminf_{n \rightarrow \infty} \frac{1}{a_n}
  \log P_n^{u,r}\{ Y_n \in G\}}\\
 && \geq  \lim_{r \rightarrow 0} \liminf_{n \rightarrow \infty} \frac{1}{a_n} 
  \log P_n\{ Y_n \in [G \cap
  \tilde{H}^{-1}(\mbox{int}\{u\}^{(r^-)})]\}
  \\ && \hspace{.3in} - \ \lim_{r \goto 0} \, \lim_{n \goto \infty} \frac{1}{a_n} \log P_n\{H_n \in \{u\}^{(r)}\} \\
 && \geq  -\lim_{r \rightarrow 0} I(G \cap \tilde{H}^{-1}
 (\mbox{int}\{u\}^{(r^-)})) + J(u)\\
 && \geq -I(G \cap \tilde{H}^{-1}(u)) + J(u)\\
 && = -I^u(G).
\end{eqnarray*}
This is the large deviation lower bound (\ref{eqn:microldlower}),
completing the second proof of the large deviation bounds in Theorem 
\ref{thm:microcan}.  The proof of Theorem \ref{thm:microcan} is done.

\skp
In Section 2 the large deviation analysis of the canonical
ensemble led us to define, in terms of the rate function in the
corresponding LDP, the set of canonical
equilibrium macrostates.  
Analogously, for $u \in \mbox{dom} \, J$ we define, in terms of the rate function $I^u$ in Theorem \ref{thm:microcan}, the set of microcanonical
equilibrium macrostates
\[
{\cal E}^u \doteq \{x \in \X : I^u(x) = 0\}.
\]
Thus $\bar{x} \in \eu$ if and only if $I(\bar{x}) = J(u)$ and 
$\tilde{H}(\bar{x}) = u$.
We next point out that in certain models elements of $\eu$ have an 
equivalent characterization in terms of constrained maximum entropy principles.

\brmk 
\label{rmk:equiv}
{\bf Equivalent characterization in terms of constrained maximum entropy principles.}
Since $J(u)$ equals the infimum of $I$ over all elements $x$ satisfying
the constraint $\tilde{H}(x) = u$, we see that $\bar{x} \in \eu$ if and only if
$\bar{x}$ solves the following constrained minimization problem:
\[
\mbox{ minimize } I(x) \mbox{ over } x \in \X \mbox{ subject to the constraint } \tilde{H}(x) = u.
\]
Both for 
the regulariued point vortex model and the Miller-Robert model
the rate function $I$ equals a relative entropy, which in turn
equals minus the physical entropy.  Hence for these models the last
display gives an equivalent characteriuation of microcanonical
equilibrium macrostates in terms of a constrained maximum entropy principle.  \ink
\ermk

Parts (c) and (d) 
of Theorem \ref{thm:canon} state several properties of the set 
${\cal E}_\beta$ of canonical equilibrium macrostates.  The next
theorem gives analogous properties of $\eu$.  The second of these
properties is slightly more complicated than in the canonical case
because the microcanonical measures $P_n^{u,r}$ 
depend on the two parameters $n \in \N$
and $r \in (0,1)$.   

\begin{thm} \per
\label{thm:iz}
We assume Hypotheses {\em \ref{hyp:prob}} and {\em \ref{hyp:general}}.
For $u \in \mbox{{\em dom}} \, J$ the following conclusions hold.

{\em (a)} ${\cal E}^u \doteq
\{x \in \X: I^u(x) = 0\}$ is a nonempty,
compact subset of $\X$.  A point $\bar{x} \in \X$ lies in $\E^u$
if and only if $I(\bar{x}) = J(u)$ and $\tilde{H}(\bar{x}) = u$;
equivalently, if and only if $\bar{x}$ solves the following constrained
minimization problem:
\[
\mbox{ minimize } I(x) \mbox{ over } x \in \X
\mbox{ subject to the constraint } \tilde{H}(x) = u.
\]

{\em (b)} Let $A$ be any Borel subset
of $\X$ whose closure $\bar{A}$ satisfies $\bar{A} \cap {\cal E}^u 
= \emptyset$. Then $I^u(\bar{A}) > 0$.  In addition,
there exists $r_0 \in (0,1)$ and for all $r \in (0,r_0]$
there exists $C_r < \infty$ such that
\[
P_n^{u,r}\{Y_n \in A\} \leq C_r \exp[-a_n I^u(\bar{A})/2] 
\goto 0 \: \mbox{ as } n \goto \infty.
\]
\end{thm}

\noi
{\bf Proof.} (a) $\E^u$ is a nonempty, compact subset of $\X$ because $I^u$
is a rate function.  The equivalent characterizations of $\bar{x} \in \E^u$
follow from the formula for $I^u$.  

(b) If $\bar{A} \cap {\cal E}^u = \emptyset$, 
then for each $x \in A$ we have 
$I^u(x)>0$.  Since $I^u$ is a rate function,
it follows that $I^u(\bar{A}) > 0$.
The large deviation upper bound for the
$P_n^{u,r}$-distributions of 
$Y_n$ given in (\ref{eqn:microldupper}) 
completes the proof. \ink

\skp
Part (b) of Theorem \ref{thm:iz} can be regarded as a concentration property
of the $P_{n}^{u,r}$-distributions of $Y_n$.  This
property justifies calling ${\cal E}^u$ 
the set of microcanonical equilibrium macrostates.

Theorem \ref{thm:weak-limits} studies compactness properties
of the sequence of $P_{n,a_n \beta}$-distributions of $Y_n$ 
and shows that any weak limit
of a convergent subsequence of this sequence 
is concentrated on ${\cal E}_\beta$.
In the next theorem we formulate an analogue
for the microcanonical ensemble, studying compactness and weak limit 
properties
of the $P_n^{u,r}$-distributions of $Y_n$. 
In the case of the classical lattice gas, a related result
is given, for example, in \cite[Lem.\ 4.1]{DeuStrZes}.
\begin{thm} \per
\label{thm:microlimits}
We assume Hypotheses {\em \ref{hyp:prob}} and {\em \ref{hyp:general}}.
For $u \in \mbox{{\em dom}} \, J$ the following conclusions hold.

{\em (a)}  For $r \in (0,1)$, any subsequence of 
$P_n^{u,r}\{Y_n \in \cdot\}$ has a subsubsequence 
$P_{n'}^{u,r}\{Y_{n'} \in \cdot\}$ converging
weakly to a probability measure $\Pi^{u,r}$ on $\X$ as $n' \goto \infty$.

{\em (b)} There exists $r_0 \in (0,1)$ such that for all $r \in (0,r_0]$
 $\, \Pi^{u,r}$ is concentrated on ${\cal E}^u$; i.e., 
$\Pi^{u,r}\{({\cal E}^u)^c\} = 0$.  Thus if $\E^u$ consists of a
unique point $\tilde{x}$, then for all $r \in (0,r_0]$ the entire
sequence $P_n^{u,r}\{Y_n \in \cdot\}$ converges weakly to
$\delta_{\tilde{x}}$ as $n \goto \infty$.

{\em (c)} For any sequence $r_k \subset (0,1)$ converging to $0$, any
subsequence of $\Pi^{u,r_k}$ has a subsubsequence converging weakly to a
probability measure $\Pi^u$ on $\X$ that is concentrated on $\E^u$.
\end{thm}

\noi
{\bf Proof.}  
(a)  Define $a^* \doteq \min_{n \in \N}a_n > 0$. 
The exponential tightness estimate in Lemma \ref{lem:exptight}
implies that for all sufficiently large $M \in (0,\infty)$
there exists a compact subset $D$ of $\X$ such that
for all $r \in (0,1)$ and all sufficiently large $n$
\be
\label{eqn:compact}
P_n^{u,r}\{Y_n \in D^c\} \leq \exp[-a_n M/2] \leq
\exp[-a^* M/2].
\ee
Since $M$ can be taken to be arbitrarily large, this yields the
tightness of the sequence $P_n^{u,r}\{Y_n \in \cdot\}$.  The tightness
implies that any subsequence of $P_n^{u,r}\{Y_n \in \cdot\}$ has a
subsubsequence $P_{n'}^{u,r}\{Y_{n'} \in \cdot\}$ converging weakly to
a probability measure $\Pi^{u,r}$ on $\X$ as $n' \goto \infty$
[Prohorov's Theorem].  This completes the proof of part (a).

(b) We use the value of $r_0$ from part (b) of Theorem \ref{thm:iz}.
As in the proof of Theorem \ref{thm:weak-limits}, in order to prove
the concentration property of $\Pi^{u,r}$, we 
write the open set $({\cal E}^u)^c$ as a
union of countably many open balls $V_j$ such that the closure
$\bar{V}_j$ of each $V_j$ has empty intersection with ${\cal E}^u$.
Let $P_{n'}^{u,r}\{Y_{n'} \in \cdot\} \Rightarrow \Pi^{u,r}$ be the
subsubsequence arising in the proof of part (a) of the present theorem.
For $r \in (0,r_0]$, part (b) of Theorem \ref{thm:iz} implies that
$P_{n'}^{u,r} \{Y_{n'} \in V_j\} \rightarrow 0$ as $n' \rightarrow
\infty$, and so
 \[
0 = \liminf_{n' \goto \infty} P_{n'}^{u,r}\{Y_{n'} \in V_j\}
\geq \Pi^{u,r}\{V_j\}.
\]
It follows that $\Pi^{u,r}\{V_j\} = 0$ and thus
that $\Pi^{u,r}\{({\cal E}^u)^c\} = 0$, as claimed.
If $\E^u$ consists of a unique point $\tilde{x}$, then
as in the proof of Theorem \ref{thm:weak-limits},
one shows that as $n \rightarrow \infty$ $\, P_{n}^{u,r}\{Y_{n} \in
\cdot \} \Rightarrow 
\delta_{\tilde{x}}$. 
This completes the proof of part (b).

(c) This follows from part (b), Prohorov's Theorem, and the compactness of
$\E^u$.  The proof of Theorem 
\ref{thm:microlimits} is complete.
\ink

\skp

\skp
\skp
\section{Equivalence and Nonequivalence of Ensembles }
\beginsec
\label{section:equiv}

In the preceding section we presented, for the microcanonical
ensemble, analogues of results proved for the canonical ensemble in
Section 2.  These include large deviation theorems and properties of
the set of equilibrium macrostates.  Such analogues of results for the
two ensembles point to a much deeper relationship between them. 
As we will soon see, the two ensembles are intimately related
both at the level of thermodynamic functions and
at the level of equilibrium macrostates, and the results at
these two levels mirror each other.  

\begin{figure}[t]
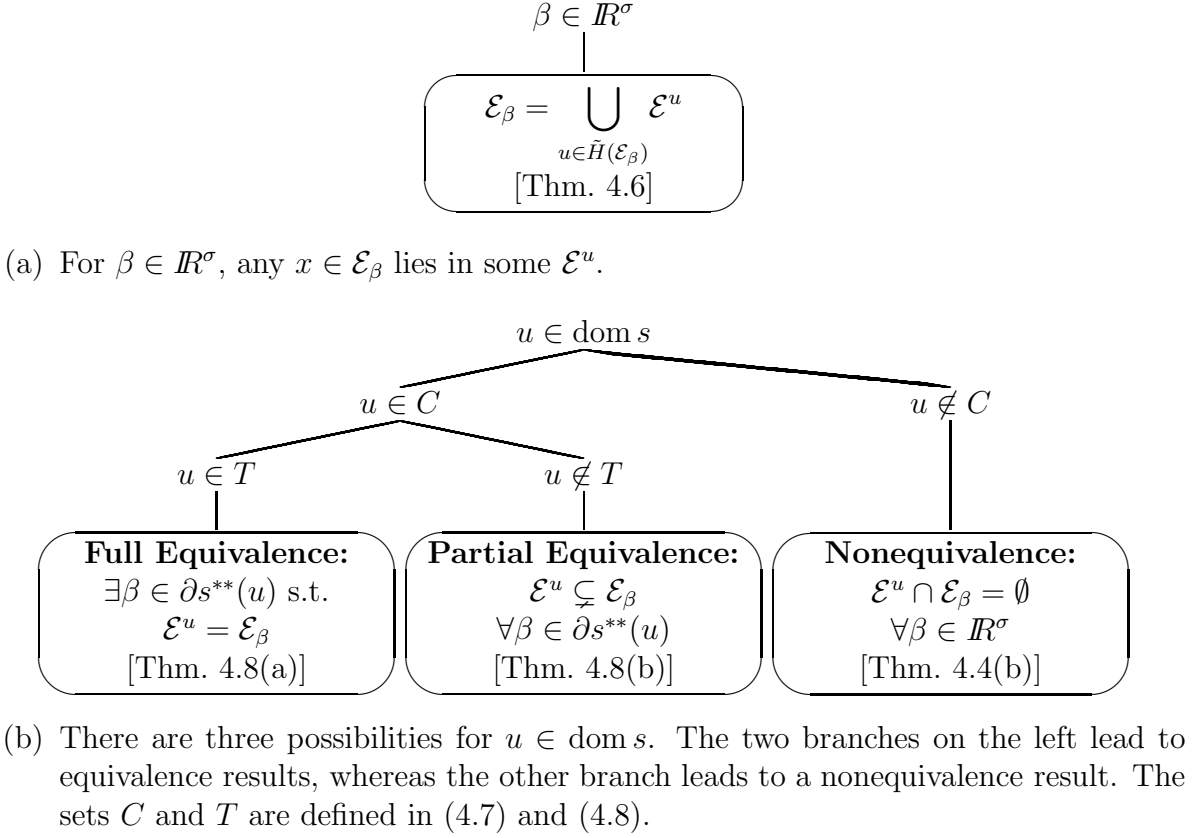

\begin{center}
\begin{Bundle}{$\beta \in \R^\sigma$}
  \chunk{\ovalbox{\begin{Bcenter} $\E_\beta = \displaystyle{\bigcup_{u
            \in \tilde{H}(\E_\beta)}} \E^u$\\ 
        \makebox[4cm][c]{$[$Thm.\ \ref{thm:beta}$]$}
      \end{Bcenter}}}
\end{Bundle}\\
\begin{itemize}
\item[(a)] For $\beta \in \R^\sigma$, any $x \in \E_\beta$ lies in
  some $\E^u$.  
\end{itemize}
\setlength{\GapDepth}{5mm} \drawwith{\drawline[10000]}
\begin{Bundle}{$u \in \mbox{dom} \, s$}
\chunk{\begin{Bundle}{$u \in C$}    
    \chunk{\begin{Bundle}{$u \in T$}
        \chunk{\ovalbox{\begin{Bcenter}
              \textbf{Full Equivalence:}\\
              \makebox[4.5cm][c]{$\exists \beta \in \partial s^{**}(u) \mbox{ s.t.}$}\\
              $\E^u = \E_\beta$\\
              $[$Thm.\ \ref{thm:full}(a)$]$
            \end{Bcenter}}}
      \end{Bundle}}
    \chunk{\begin{Bundle}{$u \not\in T$}
        \chunk{\ovalbox{\begin{Bcenter}
              \textbf{Partial Equivalence:}\\
              \makebox[4.5cm][c]{$\E^u \subsetneq \E_\beta$}\\
              $\forall \beta \in \partial s^{**}(u)$\\ 
              $[$Thm.\ \ref{thm:full}(b)$]$
            \end{Bcenter}}}
      \end{Bundle}}
  \end{Bundle}}
\chunk{%
  \setlength{\GapDepth}{14.5mm}
  \begin{Bundle}{$u \not\in C$}
    \chunk{\ovalbox{\begin{Bcenter}
          \textbf{Nonequivalence:}\\
          \makebox[4.5cm][c]{$\E^u \cap \E_\beta = \emptyset$}\\
          $\forall \beta \in \R^\sigma$\\
          $[$Thm.\ \ref{thm:z}(b)$]$
        \end{Bcenter}}}
  \end{Bundle}}
\end{Bundle}
\begin{itemize}
\item[(b)] There are three possibilities for $u \in \mbox{dom} \, s$.  The
  two branches on the left lead to equivalence results, whereas the
  other branch leads to a nonequivalence result.  The sets $C$ and $T$
  are defined in (\ref{eqn:C}) and (\ref{eqn:T}).
\end{itemize}
\caption{{}Equivalence and nonequivalence of ensembles.}
\label{fig:equiv}
\end{center}
\end{figure}

Our main results on equivalence and nonequivalence of ensembles at the
level of equilibrium macrostates are presented in Theorems
\ref{thm:z}, \ref{thm:beta}, and \ref{thm:full} and are
summarized in Figure \ref{fig:equiv}.  Definitive and complete, they
express, in terms of concavity and other properties of the
microcanonical entropy, relationships between the sets of
canonical and microcanonical
equilibrium macrostates.  The proofs of these relationships are based
on straightforward concave analysis.  Other results in this section
explore related issues.  For example, Corollary \ref{cor:beta} is a
uniqueness result for equilibrium macrostates, Theorem
\ref{thm:refinement} relates the equivalence of ensembles to the
differentiability of the canonical free energy, and Theorem
\ref{thm:beta-concave-s} shows that a certain equivalence-of-ensemble
relationship implies a concavity property of the microcanonical
entropy.  

We start our presentation by recalling an elementary result at the
level of thermodynamic functions.
The microcanonical entropy is the nonpositive function defined
for $u \in \R^\sigma$ by 
\[
s(u) \doteq -J(u) \doteq -\inf \{I(x) : x \in \X, \tilde{H}(x) = u\}.
\] 
We define $\mbox{dom} \, s$ as the set of $u \in \R^\sigma$ for which
$s(u) > -\infty$.  As shown in (\ref{eqn:fenchel}), 
the canonical free energy $\varphi(\beta)$ can be obtained
from $s$ by the formula
\be
\label{eqn:vphi-convex}
\vphi(\beta) = \inf_{u \in \R^\sigma} \{\lan \beta, u \ran - s(u)\},
\ee
which expresses $\vphi$ as the Legendre-Fenchel transform $s^*$ of $s$.
In general, $\varphi = s^*$ is finite, concave, and continuous on
$\R^\sigma$ [Thm.\ \ref{thm:canon}(a)], 
and $s$ is upper semicontinuous 
[Prop.\ \ref{prop:ldpj}(a)].  If it is the case that $s$ is concave
on $\R^\sigma$, then concave function theory implies that 
$s$ equals the Legendre-Fenchel transform of $\vphi$ \cite[p.\ 104]{Roc}; 
viz., for $u \in \R^\sigma$
\be
\label{eqn:s-concave}
s(u)= \vphi^*(u)= \inf_{\beta \in \R^\sigma} \{ \langle \beta, u \rangle -
\varphi(\beta) \}. 
\ee

If $s$ is concave on $\R^\sigma$, then
following standard terminology in the statistical mechanical
literature, we say that the canonical ensemble and the microcanonical
ensemble are thermodynamically equivalent \cite{KieLeb,LewPfiSul}.
As we will see, when properly interpreted, 
the nonconcavity of $s$ at points $u \in \R^\sigma$
will imply that the ensembles are nonequivalent at the level of
equilibrium macrostates for those values of $u$ [Thm.\ 
\ref{thm:z}(b)].  Further connections between thermodynamic
equivalence of ensembles and equivalence of ensembles at the level
of equilibrium macrostates are made explicit in Theorem \ref{thm:s-concave}.
In particular, under a hypothesis on the domains of 
various functions that is not necessarily satisfied in all models
of interest, thermodynamic equivalence of ensembles is a necessary and
sufficient condition for equivalence of ensembles to hold at the
level of equilibrium macrostates [Thm.\ \ref{thm:s-concave}(c)].

The concavity of $s$ on $\R^\sigma$ depends on the nature of $I$ and $\tilde{H}$.  For example, if $I$ is concave on $\X$ and $\tilde{H}$
is affine, then $s$ is concave on $\R^\sigma$.  However, in general the concavity of $s$ is not valid.  In fact, because of the local mean-field,
long-range nature of the Hamiltonians arising in many models of turbulence,
including the Miller-Robert model [Example \ref{exa:models}(b)], the associated
microcanonical entropies are typically not concave on subsets of $\R^\sigma$ corresponding to a range of negative temperatures.  

In order to see how concavity properties of $s$ determine
relationships between the sets of
equilibrium macrostates, we define for $u \in \R^\sigma$
the concave function
\[
s^{**}(u) \doteq \inf_{\beta \in \R^\sigma} \{\lan \beta, u \ran - 
s^*(\beta)\} = \inf_{\beta \in \R^\sigma} \{\lan \beta, u \ran - 
\vphi(\beta)\}.
\]
Because of (\ref{eqn:s-concave}), it is obvious that $s$ is concave on
$\R^\sigma$ if and only if $s$ and $s^{**}$ coincide.  
Whenever $s(u) > -\infty$ and $s(u) = s^{**}(u)$, we shall
say that $s$ is concave at $u$.  

Now assume that $s$ is not concave on $\R^\sigma$.
Since for any
$u \in \mbox{dom} \, s$ and all $\beta \in \R^\sigma $ we have
$s(u) \leq \lan \beta,u \ran - s^*(\beta)$, it follows that
for all $u \in \R^\sigma$ 
\be
\label{eqn:majorize}
s(u) \leq \inf_{\beta \in \R^\sigma}
\{\lan \beta,u \ran - s^*(\beta)\} = s^{**}(u).
\ee
In addition, if $f$ is any upper semicontinuous, concave function
satisfying $s(u) \leq f(u)$ for all $u \in \R^\sigma$, then for all
$\beta \in \R^\sigma$ $\, s^*(\beta) \geq f^*(\beta)$ and thus
$s^{**}(u) \leq f^{**}(u) = f(u)$ for all $u \in \R^\sigma$.  It
follows that if $s$ is not concave on $\R^\sigma$,
then $s^{**}$ is the upper semicontinuous, concave hull of
$s$; i.e., the smallest upper semicontinuous, concave function on
$\R^\sigma$ that majorizes $s$.  In particular, if $s(u) > -\infty$,
then $s^{**}(u) > -\infty$; thus $\mbox{dom} \, s
\subset \mbox{dom} \, s^{**}$.

Since $s^{**}$ is an upper semicontinuous, concave function, we can
introduce a basic concept in concave function theory that
will play a key role in our results on equivalence and nonequivalence
of ensembles.
For $u \in \mbox{dom}\, s^{**}$ the superdifferential of
$s^{**}$ at $u$ is defined as the set $\partial s^{**}(u)$ consisting of $\beta
\in \R^\sigma$ such that
\be
\label{eqn:chopin}
s^{**}(w) \leq s^{**}(u) + \langle \beta, w-u \rangle \mbox{ for all } w
  \in \R^\sigma;
\ee
any such $\beta$ is called a supergradient of $s^{**}$ at $u$.  The
effective domain of the superdifferential of $s^{**}$ is defined to be
the set $\mbox{dom}\, \partial s^{**}$ consisting of $u \in \R^\sigma$
for which $\partial s^{**}(u)$ is nonempty.  It can be shown that
\cite[p.\ 217]{Roc}
\be \label{eqn:domain-inclusions}
\mbox{ri(dom}\, s^{**}) \subset \mbox{dom}\, \partial s^{**} \subset
\mbox{dom}\, s^{**},
\ee
where for $A$ a subset of $\R^\sigma$ $\mbox{ri(dom} \, A)$ denotes
the relative interior of $A$. 
These relationships imply that $\partial s^{**}(u)$ is nonempty for $u \in \mbox{dom}\, s^{**}$ except possibly for $u$ in the relative boundary of
$\mbox{dom}\, s^{**}$.

The purpose of this section is to investigate,
in terms of concavity properties of $s$ and $s^{**}$, relationships between the
set $\E_\beta$ of canonical equilibrium macrostates and the set 
$\E^u$ of microcanonical equilibrium macrostates.  
We recall that for $\beta \in \R^\sigma$ and
$u \in \mbox{dom} \, s$ these
sets are defined by 
\begin{eqnarray*}
{\cal E}_{\beta} & = & \{x \in \X : I_{\beta}(x)=0\} \\ 
& = & \left\{x \in \X : I(x) + \lan 
\beta , \tilde{H}(x) \ran = \inf_{y \in \X} \{I(y) + \lan \beta,
\tilde{H}(y) \ran\} = \varphi(\beta)\right\}
\end{eqnarray*}
and 
\[
{\cal E}^u \doteq \{x \in \X : I^u(x)=0 \}=\{x \in \X : \tilde{H}(x)=u,
 I(x)= -s(u) \}.
\]
$I_\beta$ is the rate function in the
LDP for the canonical ensemble [Thm.\ 
\ref{thm:canon}], and $I^u$ is the rate
function in the LDP for the microcanonical
ensemble [Thm.\ \ref{thm:microcan}].  As the sets of points at which
the corresponding rate functions attain their minimum of 0, 
both $\ebeta$ for $\beta \in
\R^\sigma$ and $\eu$ for $u \in \mbox{dom} \, s$ are nonempty and
compact.  It is convenient to extend the definition 
of $\eu$ to all $u \in \R^\sigma$ by defining
$\eu = \emptyset$ for $u \in \R^\sigma \setminus \mbox{dom} \, s$.  

First-order differentiability conditions show that
relationships between $\E_\beta$ and $\E^u$ are plausible.  In fact,
the first-order condition for $x^*\in \X$ to be in $\E_\beta$ is \be
\label{eqn:canonicalcondition}
I'(x^*) + \langle \beta, \tilde{H}'(x^*) \rangle = 0, \ee where $'$
denotes the Frechet derivative and we assume that $I$ and $\tilde{H}$
are Frechet-differentiable.  The first-order condition for $x^* \in
\X$ to be in $\E^u$ is also (\ref{eqn:canonicalcondition}), where
$\beta$ is a Lagrange multiplier dual to the constraint
$\tilde{H}(x^*)=u$.  In order to see the precise relationships between
$\E^u$ and $\E_\beta$, we need a more detailed analysis.  

As we will see, there are three possible relationships that can
occur between $\E^u$ and $\E_\beta$.  If for a given $u \in \mbox{dom}\,
s$ there exists $\beta \in \R^\sigma$ such that 
$\E^u=\E_\beta$, then the ensembles are
said to be fully equivalent or that full equivalence of ensembles holds.  
If instead of equality $\E^u$ is a proper
subset of $\E_\beta$ for some $\beta \in \R^\sigma$, 
then the ensembles are said to be partially
equivalent or that partial equivalence of ensembles holds.  
It may also happen that $\E^u \cap \E_\beta = \emptyset$
for all $\beta \in \R^\sigma$.  If this occurs, then the
microcanonical ensemble is said to be nonequivalent to any canonical
ensemble or that nonequivalence of ensembles holds.   It is convenient
to group the first two cases together.  If for a given $u$
there exists $\beta$ such that either $\E^u$ equals $\E_\beta$ or
$\E^u$ is a proper subset of $\E_\beta$, then the ensembles are said
to be equivalent or that equivalence of ensembles holds. 

The probabilistic role played by 
$\E^u$ and $\E_\beta$ should be kept in mind when
interpreting these relationships.  According to part (c) of Theorem
\ref{thm:canon}, for any Borel subset $A$ whose closure is disjoint
from $\E_\beta$, $P_{n, a_n\beta}\{Y_n \in A \} \rightarrow 0$.
Theorem \ref{thm:weak-limits} refines this by showing that convergent
subsequences of $P_{n, a_n\beta}\{Y_n \in \cdot \}$ have weak limits
with support in $\E_\beta$.  Theorems \ref{thm:iz} and
\ref{thm:microlimits} do the same for the microcanonical ensemble.
Only when $\ebeta = \eu = \{x\}$ can we be sure that the two ensembles
give the same prediction in the sense of weak convergence.  
A condition implying these equalities is
given in Corollary \ref{cor:beta}.

A key insight revealed by our results is that the set $\eu$ of
microcanonical equilibrium macrostates can be richer than the set
${\cal E}_\beta$ of canonical equilibrium macrostates.  Specifically,
every $x \in \E_\beta$ is also in $\E^u$ for some $u$, but if the
microcanonical entropy $s$ is not concave at some $u$, then any $x \in
\E^u$ does not lie in $\E_\beta$ for any $\beta$ (nonequivalence
of ensembles).  This verbal
description is made precise in Theorems \ref{thm:z} and \ref{thm:beta}, 
while Theorems \ref{thm:z} and \ref{thm:full} give necessary and sufficient
conditions for equivalence of ensembles to hold.  The content of Theorem
\ref{thm:beta} is summarized in Figure \ref{fig:equiv}(a).  The contents
of Theorems \ref{thm:z} and \ref{thm:full} are summarized in
Figure \ref{fig:equiv}(b). 

Theorem \ref{thm:z} gives a geometric condition that is necessary and sufficient
for equivalence of ensembles to hold.  We define $C$ to be the set
of $u \in \R^\sigma$ for which there exists a supporting hyperplane to
the graph of $s$ at $(u,s(u))$.  In symbols,
\be 
\label{eqn:C}
C \doteq \{u \in \R^\sigma : \exists \beta \in \R^\sigma \ni s(w) \leq s(u) + \langle \beta, w-u \rangle \mbox{ for
  all } w \in \R^\sigma \}.
\ee
If $u \in C$, then the $\beta$ appearing in this display 
is a normal vector to the supporting hyperplane. 
According to part (a) of Theorem \ref{thm:z}, for a particular $u \in
\mbox{dom} \, s$ equivalence of ensembles holds if and 
only if $u \in C$.  According to part (b) of the theorem, for a particular
$u \in \mbox{dom} \, s$ 
nonequivalence of ensembles holds if and only if $u \not \in C$.

Theorem \ref{thm:full} refines part (a) of Theorem \ref{thm:z} by giving a geometric condition that is necessary and sufficient
for full equivalence of ensembles to hold.
We define $T$ to be the set of $u \in \R^\sigma$ for which there exists a
supporting hyperplane to the graph of $s$ that touches the graph of $s$
only at $(u,s(u))$.  In symbols, 
\be 
\label{eqn:T}
T \doteq \{u \in \R^\sigma : \exists \beta \in \R^\sigma \ni
  s(w) < s(u) + \langle \beta, w-u \rangle \mbox{ for
  all } w \not=u \}.
\ee
Clearly, $T$ is a subset of $C$, which is the set of $u$ for which
equivalence of ensembles holds [Thm.\ \ref{thm:z}(a)].
According to Theorem \ref{thm:full},
for a particular $u \in \mbox{dom} \, s$ full equivalence of ensembles
holds if and only if $u \in T$.

Before proving any results on the equivalence and nonequivalence
of ensembles, we point out an alternate representation of $C$
that will elucidate the connection between these results and concavity
properties of $s$ and $s^{**}$.  In general $s$ is not concave on $\R^\sigma$.
According to part (b) of Lemma \ref{lem:CtoC}, $C$ equals the set of $u \in
\mbox{dom} \, \partial s^{**}$ at which $s$ is concave; i.e., the set of 
$u \in \mbox{dom} \, \partial s^{**}$ such that $s(u)$ equals the value at $u$
of the concave function $s^{**}$.  
It follows 
from part (b) of Lemma \ref{lem:CtoC} that if $s$ is not concave at some 
$u \in \mbox{dom} \, s$, then $u \not \in C$ and so nonequivalence of
ensembles holds [Thm.\ \ref{thm:z} (b)].

It is easy to find a sufficient condition on $s^{**}$ for
full equivalence of ensembles to hold.
Suppose that for some $u \in \R^\sigma$
$s(u) = s^{**}(u)$ and that there exists $\beta \in 
\R^\sigma$ such that 
\be
\label{eqn:bach} 
s^{**}(w) < s^{**}(u) + \langle \beta, w-u \rangle \mbox{ for all } w
  \not= u;
\ee
i.e., the inequality (\ref{eqn:chopin}) defining $\beta \in \partial s^{**}(u)$ 
holds with strict inequality for all $w \not= u$.
Since $s(w) \leq s^{**}(w)$, it follows that 
\be
\label{eqn:ineqt}
s(w) < s(u) + \langle \beta, w-u \rangle \mbox{ for all } w
  \not= u.
\ee
That is, $u$ lies in $T$, which according to Theorem \ref{thm:full}
is the subset of $\R^\sigma$ for which full equivalence of ensembles holds.
If, for example, $s^{**}$ is strictly concave in a neighborhood of $u$,
then (\ref{eqn:bach}) holds for any $\beta \in \partial s^{**}(u)$ and
thus we have full equivalence of ensembles.  

In order to find a sufficient condition on $s^{**}$ for
partial equivalence of ensembles to hold, let $u$ be a point in $\R^\sigma$
such that $s^{**}$ is affine in a neighborhood of $u$. 
Then except in pathological cases, 
for any $\beta \in \R^\sigma$ the strict inequality
(\ref{eqn:ineqt}) cannot be valid for all $w \not = u$, and so partial
equivalence of ensembles holds.  

Part (b) of the next lemma gives the alternate representation
of $C$ to which we referred three paragraphs earlier.  This representation
involves the set 
\[
\Gamma \doteq \{u \in \R^\sigma : s(u) = s^{**}(u)\}.
\]

\begin{lemma} \per 
\label{lem:CtoC}
{\em (a)} For $u$ and $\beta$ in $\R^\sigma$, $s(w) \leq s(u) + \langle
\beta, w-u \rangle$ for all $w \in \R^\sigma$ if and only if both $s(u)=
s^{**}(u)$ and $\beta \in \partial s^{**}(u)$.

{\em (b)} $C = \Gamma \cap \mbox{\em dom}\, \partial s^{**}$, and
$C \subset \Gamma \cap \mbox{\em dom}\, s$.
\end{lemma}

\brmk
\label{rmk:dell}
It is not difficult to refine the second assertion in part (b) of
this lemma by showing that
\[
\Gamma \cap \mbox{ri(dom} \, s) \subset C = 
\Gamma \cap \mbox{dom}\, \partial s^{**} \subset \Gamma \cap 
\mbox{dom} \, s.
\]
This relationship implies that, except possibly for relative boundary
points of $\mbox{dom} \, s$, $C$ consists of $u \in \mbox{dom} \, s$ for 
which $s(u) = s^{**}(u)$.  According to Theorem \ref{thm:z}, equivalence
of ensembles holds for a particular $u \in \mbox{dom} \, s$ if and only
if $u \in C$.  Combining this with the observation in the preceding 
sentence, we see that, except possibly for relative boundary
points of $\mbox{dom} \, s$, equivalence of ensembles holds
for $u \in \mbox{dom} \, s$ if and only if $s(u) = s^{**}(u)$.
\ermk

\noi \textbf{Proof of Lemma \ref{lem:CtoC}.}  (a) We start the proof 
by first assuming that $s(w) \leq s(u) + \langle \beta,
w-u \rangle$ for all $w \in \R^\sigma$.  It follows that $u \in \mbox{dom} \, s$ and that $\langle \beta, u \rangle - s(u) \leq \langle
\beta, w \rangle - s(w)$ for all $w \in \R^\sigma$.  Therefore
\[
\langle \beta, u \rangle - s(u) = \inf_{w \in \R^\sigma} \{ \langle
\beta, w \rangle - s(w) \} = \varphi(\beta).
\]
Since $s^{**}(w) = \inf_{\gamma \in \R^\sigma}\{\langle \gamma, w
\rangle - \varphi(\gamma)\} \leq \langle \beta, w \rangle -
\varphi(\beta)$, the last display and the
inequality $s(u) \leq s^{**}(u)$ imply that for all $w \in \R^\sigma$
\begin{eqnarray*}
s^{**}(w) &\leq& \langle \beta, w \rangle -\varphi(\beta)
= s(u) +\langle \beta, w \rangle -\langle \beta, u \rangle\\
&
\leq& s^{**}(u) + \langle \beta, w-u \rangle.
\end{eqnarray*}
Thus $\beta \in \partial s^{**}(u)$.  Setting $w=u$ 
yields $s(u)=s^{**}(u)$.

Now assume that $s(u) = s^{**}(u)$ and that $\beta \in \partial s^{**}(u)$; 
thus for all $w \in \R^\sigma$
\[
s^{**}(w) \leq s^{**}(u) + \langle \beta, w -u \rangle=s(u) + \langle
\beta, w-u \rangle.
\]
Since $s(w) \leq s^{**}(w)$ for all $w \in
\R^\sigma$, it follows that for all $ w \in \R^\sigma$
\[
s(w) \leq s(u) + \langle \beta , w-u \rangle.
\]
This completes the proof of part (a).

(b) The first assertion is an immediate consequence of part (a).
As mentioned in the proof of part (a), if $u \in C$, then $u \in
\mbox{dom} \, s$.  We conclude that $C \subset \Gamma \cap 
\mbox{dom} \, s$, as claimed.  \ink

\skp
The next lemma will facilitate the proofs of a number of our results on the equivalence and nonequivalence of ensembles.  Part (b) refines 
one of the conditions 
in part (a), substituting a weaker hypothesis that leads to
the same conclusion.

\begin{lemma} \per \label{lem:z}
For $u$ and $\beta \in \R^\sigma$ the following conclusions hold.

{\em (a)} The inequality $s(w) \leq s(u) + \langle \beta , w-u \rangle$ is valid
for all $w
\in \R^\sigma$ if and only if $\E^u \not= \emptyset$ and $\E^u \subset
\E_\beta$.

{\em (b)}  If $\E^u \cap \E_\beta \not= \emptyset$, then 
$s(w) \leq s(u) + \langle \beta , w-u \rangle$ for all $w
\in \R^\sigma$.  
\end{lemma}

\noi \textbf{Proof.}  We first prove that if
$s(w) \leq s(u) + \langle \beta,
w-u \rangle$ for all $w \in \R^\sigma$, then $\eu \not = \emptyset$
and $\eu \subset \ebeta$.  The hypothesis implies that $u \in 
\mbox{dom} \, s$ and that $\langle \beta, u \rangle -s(u) \leq \langle
\beta ,w \rangle - s(w)$ for all $ w \in \R^\sigma$.  Therefore
\[
\langle \beta, u \rangle -s(u) = \inf_{w \in \R^\sigma} \{ \langle
\beta ,w \rangle - s(w)\} = \varphi(\beta)= \inf_{y \in \X} \{ \langle
\beta , \tilde{H}(y) \rangle + I(y) \}.
\]
The fact that $u$ is an element of $\mbox{dom}\, s$ implies that $\E^u \not=
\emptyset$.  Let $x$ be an arbitrary element in $\E^u$.  Since
$\tilde{H}(x)=u$ and $I(x)= -s(u)$, the display implies 
that \[
\langle \beta, \tilde{H}(x) \rangle +I(x) = \inf_{y \in \X} \{
\langle \beta, \tilde{H}(y) \rangle + I(y) \}
\]
 and thus that $x \in
\E_\beta$.  Since $x$ is an arbitrary element in $\E^u$, it follows that $\E^u \subset \E_\beta$.  

In order to complete the proof of part (a), it suffices
to prove part (b).
Thus suppose that $\E^u \cap \E_\beta \not= \emptyset$ and 
let $x$ be an arbitrary element in 
$\E^u \cap \E_\beta$.  Since $\E^u \not = \emptyset$, 
we have $u \in \mbox{dom} \, s$. In addition, 
since $\tilde{H}(x) = u$, $I(x) =-s(u)$, and 
\[
\langle
\beta, \tilde{H}(x) \rangle + I(x) = \inf_{y \in \X} \{ \langle \beta,
\tilde{H}(y) \rangle +I(y) \} = \varphi(\beta),
\]
it follows that for all $w \in \R^\sigma$
\[
\langle \beta, u \rangle - s(u)=
\varphi(\beta)
=\inf_{w' \in \R^\sigma} \{ \langle \beta, w' \rangle -s(w') \} \leq
\langle \beta, w \rangle - s(w).
\]
Therefore $s(w) \leq s(u) + \langle \beta, w-u \rangle$ for all $w \in
\R^\sigma$, as claimed.  \ink

\skp

The next theorem is our first main result.  Part (a) states that for a
particular $u\in \mbox{dom}\, s$ equivalence of
ensembles holds if and only if $u \in C$.  In Theorem \ref{thm:s-concave}
we make explicit the connection between part (a) and the relationship
between thermodynamic equivalence of ensembles and equivalence
of ensembles at the level of equilibrium macrostates.  
Part (b) of the next theorem states that for a
particular $u\in \mbox{dom}\, s$ nonequivalence of ensembles holds if
and only if $u \not\in C$.  In particular, if $s$ is not concave at some
$u \in \mbox{dom} \, \partial s^{**}$, then the ensembles are nonequivalent at 
the level of equilibrium macrostates.  Theorem \ref{thm:z} was
inspired by, and greatly improves upon, the presentation on pages 857-859
of \cite{EyiSpo}, which treats the regularized point vortex model.
While part (b) of Theorem \ref{thm:z} is related to part (b) of Lemma
5.1 in \cite{LewPfiSul2}, our Theorem \ref{thm:z} makes the nonequivalence
of ensembles more explicit.

\begin{thm} \per \label{thm:z}
  We assume Hypotheses {\em \ref{hyp:prob}} and {\em
    \ref{hyp:general}}. For $u \in \mbox{\em dom}\, s$ the following
conclusions hold.  

{\em (a)}  $u \in C$ if and only if $\E^u \subset \E_\beta$ for some
$\beta \in \R^\sigma$. 

{\em (b)}  $u \not\in C$ if and only if  $\E^u \cap \E_\beta = \emptyset$ for
all $\beta \in \R^\sigma$.
\end{thm}

\noi
{\bf Proof.}  (a) This is an immediate
 consequence of part (a) of Lemma \ref{lem:z}.

(b) If $u \not \in C$, then for any $\beta \in \R^\sigma$
the inequality $s(w) \leq s(u) + \lan \beta, w-u \ran$ does
not hold for all $w \in \R^\sigma$.  Part (b) of Lemma \ref{lem:z} implies
that $\eu \cap \ebeta = \emptyset$ for all $\beta \in \R^\sigma$.
To show the converse, assume that 
$\E^u \cap \E_\beta = \emptyset$ for
all $\beta \in \R^\sigma$ and that $u \in C$.
But if $u \in C$, then part (a) of Lemma \ref{lem:z} implies that 
$\E^u \subset \E_\beta$ for some $\beta \in \R^\sigma$.  This contradiction
shows that $u \not \in C$, completing the proof.  \ink 

\skp

In the next proposition we refine part (a) of Theorem \ref{thm:z}
by specifying the set of $\beta$ for which $\eu \subset \ebeta$.

\begin{prop} \per \label{prop:z}
  We assume Hypotheses {\em \ref{hyp:prob}} and {\em
    \ref{hyp:general}}. Then for $u \in C$, 
$ \E^u \subset \E_\beta$ for all $\beta \in \partial
s^{**}(u)$ and $\E^u \cap \E_\beta = \emptyset$ for all $\beta \not\in
\partial s^{**}(u)$.
\end{prop}

\noi \textbf{Proof.} For $u \in C$, part (b) of Lemma
\ref{lem:CtoC} implies that $s(u) = s^{**}(u)$ and $\partial s^{**}(u)
\not= \emptyset$. If $\beta \in \partial s^{**}(u)$, then part (a) of
the same lemma implies that $s(w) \leq s(u) + \langle \beta, w-u
\rangle$ for all $w \in \R^\sigma$.   Part (a) of Lemma
\ref{lem:z} then implies that $\E^u \subset \E_\beta$.
This proves the first half of the proposition.
On the other hand, if $\beta \not\in \partial s^{**}(u)$, then 
it is not true that $s(w) \leq s(u) +
\langle \beta , w-u \rangle$ for all $w \in \R^\sigma$ 
[Lem.\ \ref{lem:CtoC}(a)].  It follows from part (b) of Lemma \ref{lem:z}
that $\E^u \cap \E_\beta = \emptyset$.
\ink

\skp
Theorem \ref{thm:z} considers $u \in \mbox{dom} \, s$, 
proving that partial or full equivalence of ensembles holds if and only
if $u \in C$.  The next theorem is our second main result.
It shifts focus from $u \in \mbox{dom} \, s$
to $\beta \in \R^\sigma$, proving that every set 
$\E_\beta$ of canonical equilibrium
macrostates is a disjoint union of
$\E^u$ for $u$ in a particular index set that depends on $\beta$.  

\begin{thm} \per
\label{thm:beta}  
We assume Hypotheses {\em \ref{hyp:prob}} and 
{\em \ref{hyp:general}}.  Then for all $\beta \in \R^\sigma$,
$\tilde{H}(\E_\beta) \subset \mbox{{\em dom}}\, s$ and
\[ 
  \E_\beta= \bigcup_{u \in \tilde{H}(\E_\beta)} \E^u. 
\] 
The sets $\E^u$, $u \in \tilde{H}(\E_\beta)$, are nonempty and disjoint.
\end{thm}

\noindent
\textbf{Proof.} Let $x$ be an arbitrary element in $\E_\beta$ and
define $\tilde{u} \doteq \tilde{H}(x)$.  Since $I_\beta(x) = 0$, we have
\[
I(x) + \langle \beta,\tilde{H}(x) \rangle = \inf_{y \in \X}\{I(y) + 
\langle \beta,\tilde{H}(y) \rangle\} < \infty,
\]
and so $s(\tilde{u}) \geq -I(x) > -\infty$. Thus $\tilde{u} \in
\mbox{dom} \, s$.  Because $x$ is an arbitrary element in $\E_\beta$, 
this proves that $\tilde{H}(\E_\beta) \subset
\mbox{dom}\, s$.  Since $\tilde{u} \in \mbox{dom}\, s$,
$\E^{\tilde{u}}$ can be characterized as the set of $x \in \X$
satisfying $\tilde{H}(x) = \tilde{u}$ and $I(x) = -s(\tilde{u})$.

We now prove that $x \in \E^{\tilde{u}}$.
Since $x \in \E_{\beta}$, it follows that for any $y \in \X$
\[
I(x) + \lan \beta,\tilde{u} \ran = 
I(x) + \langle \beta, \tilde{H}(x) \rangle \leq I(y)
+ \langle \beta,\tilde{H}(y) \rangle,
\] 
and thus for any $y \in \X$ satisfying $\tilde{H}(y) = \tilde{u}$, we have
$I(x) \leq I(y)$.  This implies that 
\[
I(x) \leq \inf\{I(y) : y \in \X, \tilde{H}(y) = \tilde{u}\} 
= -s(\tilde{u}) \leq I(x),
\]
and so $I(x) = -s(\tilde{u})$.
It follows that $x \in \E^{\tilde{u}}$.  Since $x$ is an arbitrary
element in $\E_\beta$, we have shown that
\[
\E_\beta \subset \bigcup_{u \in \tilde{H}(\E_\beta)} \E^u.
\]

In order to prove the reverse inclusion, we show that for any $u \in \tilde{H}(\E_\beta)$ we have $\E^u \subset \E_\beta$.  Any such $u$
has the form $u=\tilde{H}(y)$ for some $y \in \E_\beta$.  From our work in the
preceding two paragraphs we know that $u \in \mbox{dom}\, s$ and $y
\in \E^u$.  Thus $y \in \E^u \cap \E_\beta$.  Since 
$\E^u \cap \E_\beta \not= \emptyset$, it follows from Theorem \ref{thm:z}
that $\E^u \subset \E_\beta$.  This completes the proof of the
display in the theorem.

The sets $\E^u, u \in \tilde{H}(\E_\beta)$, are nonempty since any such
$u$ lies in $\mbox{dom}\, s$. The sets are also disjoint since for $u
\not= u', x \in \E^u \cap \E^{u'}$ implies that $\tilde{H}(x)$ 
equals both $u$ and $u'$.  The proof of the theorem is complete.
\ink

\skp

The following useful corollary states that when $\E_\beta$ consists of
a unique point $x$, then with $\tilde{u} \doteq \tilde{H}(x)$,
$\E^{\tilde{u}}$ consists of the unique point $x$.  This follows from
Theorem \ref{thm:beta} since $\tilde{H}(\E_\beta)= \{\tilde{H}(x)\}$.
The corollary sharpens the result on page 861 of \cite{EyiSpo}, which
needs the additional hypotheses that $s$ is strictly concave and essentially
smooth in order to reach the same conclusion.

\begin{cor} \per 
\label{cor:beta}
Suppose that $\E_{\beta}=\{x\}$ for some $\beta \in \R^{\sigma}$.  Then
$\E^{\tilde{u}}=\{x\}$, where $\tilde{u} \doteq 
\tilde{H}(x)$.
\end{cor}


We now turn our attention to a criterion for full equivalence of
ensembles, which is stated in terms of the set $T$ defined in 
(\ref{eqn:T}).  Part (a) of Theorem \ref{thm:z} states that for a
particular $u \in \mbox{dom}\, s$ equivalence of
ensembles holds if and only if $u \in C$.  The next
theorem refines this by showing 
that full equivalence of ensembles holds if and only
if $u \in T$.  Part (a) gives the sufficiency and part (b) the necessity.

\begin{thm} \per \label{thm:full}
  We assume Hypotheses {\em \ref{hyp:prob}} and {\em
    \ref{hyp:general}}.  The following conclusions hold.

{\em (a)} If $u \in T$, then there exists $\beta \in \partial s^{**}(u)$ such
that $\E^u=\E_\beta$.

{\em (b)} If $u \in C \setminus T$, then $\E^u \subsetneq \E_\beta$ for
all $\beta \in \partial s^{**}(u)$ and $\E^u \cap \E_\beta =
\emptyset$ for all $\beta \not\in \partial s^{**}(u)$. 
\end{thm}

\noi \textbf{Proof.} (a)  If  $u \in T$, then there exists
$\beta \in \R^\sigma$ such that $s(w) < s(u) + \langle \beta, w-u
\rangle$ for all $w \not=u$.  Part (a) of Lemma \ref{lem:z} implies
that $\E^u \subset \E_\beta$.  Suppose that $\E^u$ is a proper subset
of $\E_\beta$.  Then Theorem \ref{thm:beta} implies the existence of 
$u' \not=u$ such that $\E^{u'} \not= \emptyset$ and $\E^{u'} \subset
\E_\beta$, and part (a) of Lemma \ref{lem:z} yields
\[
s(w) \leq s(u') + \langle \beta, w-u' \rangle \mbox{ for all } w \in
\R^\sigma. 
\]
Setting $w=u$ 
and using the fact that $s(u') <
s(u) + \langle \beta, u'-u \rangle$, we see that
\[
s(u) \leq s(u') + \langle \beta , u-u' \rangle
< s(u) + \langle \beta, u'-u \rangle + \langle \beta, u-u' \rangle
=s(u).
\]
This contradiction shows that the assumption that
$\E^u$ is a proper subset of $\E_\beta$ is false.  The proof of part (a)
is complete.

(b)  For $u \in C \setminus T$, Proposition \ref{prop:z} implies that
$\E^u \subset \E_\beta$ for all $\beta \in \partial s^{**}(u)$ and
$\E^u \cap \E_\beta = \emptyset$ for all $\beta \not\in \partial
s^{**}(u)$.  We now show that for any $\beta \in \partial
s^{**}(u)$, $\E^u$ is a proper subset of $\E_\beta$.  
Since $\E^u \subset \E_\beta$, part (a) of Lemma \ref{lem:z}
implies that $s(w) \leq s(u) + \langle \beta, w-u \rangle$ for all $w
\in \R^\sigma$.  Since $u \not\in T$, there exists $u' \not= u$ 
such that $s(u')=s(u) + \langle \beta, u'-u \rangle$.  Then for
all $w \in \R^\sigma$
\[
s(w) \leq s(u) + \langle \beta, w-u \rangle = s(u') + \langle \beta,
w- u' \rangle.
\]
It now follows from part (a) of Lemma \ref{lem:z} that $\E^{u'} \not= \emptyset$
and $\E^{u'} \subset \E_\beta$.  Thus $\E^u$ is a proper subset of
$\E_\beta$, as claimed.  \ink

\skp

We recall that thermodynamic equivalence of ensembles
is said to hold when $s$ is concave on $\R^\sigma$.
The next theorem addresses the issue of how thermodynamic equivalence
of ensembles 
mirrors equivalence of ensembles at the level of 
equilibrium macrostates.  Part (a) shows that thermodynamic equivalence
is a sufficient condition for macroscopic equivalence to hold for all
$u \in \mbox{dom} \, \partial s$.  Since when $s$ is concave on $\R^\sigma$
we have $\mbox{ri(dom} \, s) \subset 
\mbox{dom} \, \partial s \subset \mbox{dom} \, s$, it follows that thermodynamic
equivalence is a sufficient condition for macroscopic equivalence to hold 
for all $u \in \mbox{dom} \, s$ except possibly 
for relative boundary points.  Part (b) proves a partial converse to (a).
In part (c) we point out that thermodynamic equivalence is equivalent
to macroscopic equivalence under an extra hypothesis on the domains of
$s$, $s^{**}$, and $\partial s^{**}$.  The proof of the theorem follows
readily from our previous results.  The theorem is related to
Lemma 6.2 and Theorem 6.1 in \cite{LewPfiSul2}.

\begin{thm} \per
\label{thm:s-concave}
{\em (a)} Assume that $s$ is concave on $\R^\sigma$.  Then for all
$u \in \mbox{{\em dom}} \, \partial s$, $\eu \subset \ebeta$ for
some $\beta \in \partial s(u)$.  Thus, thermodynamic equivalence
of ensembles implies equivalence of ensembles at the level of equilibrium
macrostates for all $u \in \mbox{{\em dom}} \, \partial s$.

{\em (b)}  Assume that $\mbox{{\em dom}} \, s = \mbox{{\em dom}} \, s^{**}$
and that for all $u \in \mbox{{\em dom}} \, s$ there exists $\beta \in \R^\sigma$
such that $\eu \subset \ebeta$.  Then $s$ is concave on $\R^\sigma$.
Thus, under the hypothesis that 
$\mbox{{\em dom}} \, s = \mbox{{\em dom}} \, s^{**}$, equivalence of ensembles
at the level of equilibrium macrostates for all $u \in \mbox{{\em dom}} \, s$
implies thermodynamic equivalence of ensembles.

{\em (c)}  Assume that $\mbox{{\em dom}} \, s = \mbox{{\em dom}} \, s^{**}
= \mbox{{\em dom}} \, \partial s^{**}$.  Then thermodynamic equivalence of
ensembles holds if and only if the ensembles are equivalent at the level of
equilibrium macrostates.
\end{thm}

\noi
{\bf Proof.}  (a) If $s$ is concave on $\R^\sigma$, then $s = s^{**}$ on $\R^\sigma$ and $C = \mbox{dom} \, \partial s^{**} = \mbox{dom} \, \partial s$
[Lem.\ \ref{lem:CtoC}(b)].  Part (a) of
Theorem \ref{lem:z} completes the proof of part (a). 

(b) The hypotheses imply that any element of $\mbox{dom} \, s$ is an 
element of $C$, which in turn is a subset of $\Gamma \doteq
\{u \in \R^\sigma : s(u) = s^{**}(u)\}$.  It follows that $s$ and $s^{**}$
agree on $\mbox{dom} \, s = \mbox{dom} \, s^{**}$ and thus that $s$ is concave
on $\R^\sigma$.

(c)  This follows from parts (a) and (b).  \ink

\skp

With Theorem \ref{thm:s-concave} the presentation of the main results
in this section is complete.  We end this section by giving two additional
theorems in which we explore further relationships involving $\ebeta$,
$\eu$, and the thermodynamic functions $\vphi$ and $s$.

In part (a) of the next theorem we refine Theorem \ref{thm:beta} by proving 
that $\E_\beta =\bigcup_{u \in \partial \varphi(\beta) \cap \Gamma}
\E^u$, where $\partial \varphi(\beta)$ denotes the superdifferential
at $\beta$ of the concave function $\varphi$ and, as introduced
in Lemma \ref{lem:CtoC}, $\Gamma \doteq \{u \in \R^\sigma: s(u) = s^{**}(u)\}$.  This in turn allows us
to give, in part (b), a necessary and sufficient condition for the
differentiability of $\varphi$ at a point $\beta$.  Part (c) is a special
case of part (b).

\begin{thm} \per \label{thm:refinement}
We assume Hypotheses {\em \ref{hyp:prob}} and 
{\em \ref{hyp:general}}.  The following conclusions hold.

{\em (a)}  For all $\beta \in \R^\sigma$
\[ 
  \E_\beta= \bigcup_{u \in \tilde{H}(\E_\beta)} \E^u =\bigcup_{u \in
    \partial \varphi(\beta) \cap \Gamma} \E^u.
\] 

{\em (b)}  $\varphi$ is differentiable at $\beta$ if and only if both
  $\E_\beta= \E^u$ for some $u$ and $\partial \varphi(\beta) \subset
  \Gamma$.

{\em (c)} If $s$ is concave on $\R^\sigma$, then $\varphi$ is
  differentiable at $\beta$ if and only if $\E_\beta= \E^u$ for some $u$.
\end{thm}

\noi \textbf{Proof.} (a) It follows from
part (a) of Lemma \ref{lem:z} and part (a) of Lemma \ref{lem:CtoC}
that
\[
\E^u \not= \emptyset \mbox{ and } \E^u \subset \E_\beta \mbox{ if and
  only if } s(u)=s^{**}(u) \mbox{ and } \beta \in \partial s^{**}(u).
\]
Since $\beta \in \partial s^{**}(u)$ if and only if $u
\in \partial s^*(\beta)= \partial \varphi(\beta)$ \cite[p.\ 218]{Roc},
it follows that
\be
\label{eqn:kyle}
\E^u \not= \emptyset \mbox{ and } \E^u \subset \E_\beta
\mbox{ if and only if } u \in \partial \varphi(\beta) \cap \Gamma.
\ee
Thus 
\[
\bigcup_{u \in
    \partial \varphi(\beta) \cap \Gamma} \E^u \subset \E_\beta.
\]

We complete the proof of part (a) by showing that we have equality in this
display.  By Theorem \ref{thm:beta} $\E_\beta$ is a disjoint union of $\E^u$ for
$u \in \tilde{H}(\E_\beta) \subset \mbox{dom} \, s$.   Hence for each 
$u \in \tilde{H}(\E_\beta)$, $\E^u \not= \emptyset$ and $\E^u \subset \E_\beta$.
Thus (\ref{eqn:kyle}) implies that $\tilde{H}(\E_\beta) \subset
\partial \varphi(\beta) \cap \Gamma$.  We conclude that
\[
\bigcup_{u \in
    \partial \varphi(\beta) \cap \Gamma} \E^u \subset \E_\beta
= \bigcup_{u \in \tilde{H}(\E_\beta)} \E^u \subset
\bigcup_{u \in
    \partial \varphi(\beta) \cap \Gamma} \E^u,
\]
and therefore $\bigcup_{u \in
    \partial \varphi(\beta) \cap \Gamma} \E^u = \E_\beta$.

(b) We first assume that $\varphi$ is differentiable at $\beta$.
Since by part (a) $\partial
\varphi(\beta) \cap \Gamma \not= \emptyset$ for any $\beta$, the differentiability of $\vphi$ at $\beta$ implies that
$\partial \varphi(\beta)= \{\nabla \varphi(\beta)\} \subset \Gamma$ and
that $\E_\beta=\E^{\nabla \varphi(\beta)}$.
We now assume that $\E_\beta=\E^u$ for some $u$ and $\partial
\varphi(\beta) \subset \Gamma$.  Since part (a) implies that $\partial
\varphi(\beta) \cap \Gamma =\{ u \}$,
we conclude that $\partial \varphi(\beta)= 
\partial \varphi(\beta) \cap \Gamma = \{u\}$
and therefore that $\varphi$ is differentiable at $\beta$.

(c) This follows from part (b) since the concavity of $s$ on
$\R^\sigma$ implies that $\Gamma=\R^\sigma$, and so $\partial
\varphi(\beta) \subset \Gamma$ is always true.
\ink

\skp

The next theorem is the final result in this section.
Under the hypothesis that $s$ is concave on $\R^\sigma$,
part (a) gives a simpler form of the 
representation in part (a) of Theorem
\ref{thm:refinement}.  Part (b) is a partial converse of part (a).
\begin{thm} \per
\label{thm:beta-concave-s}
We assume Hypotheses {\em \ref{hyp:prob}} and 
{\em \ref{hyp:general}}.   The following conclusions hold.

{\em (a)} Assume that $s$ is concave on $\R^\sigma$.  
Then for all $\beta \in \R^\sigma$
\[
\E_\beta = \bigcup_{u \in \tilde{H}(\E_\beta)} \E^u=
\bigcup_{u \in \partial \varphi(\beta)} \E^u.
\]

{\em (b)}  Now assume that for all $\beta \in \R^\sigma$
\[
\E_\beta = \bigcup_{u \in \partial \varphi(\beta)} \E^u.
\]
Then $s$ is a finite concave function on any convex subset of
$\mbox{\em ri(dom}\, s)$.
\end{thm}

\noi \textbf{Proof.} (a) Since $s$ is concave on $\R^\sigma$, $\Gamma$ equals
$\R^\sigma$ and thus $\partial \varphi(\beta) \cap \Gamma = \partial
\varphi(\beta)$ for all $\beta \in \R^\sigma$.  Hence part (a) follows
from part (a) of Theorem \ref{thm:refinement}.

(b)  Since by definition
$\E^u = \emptyset$ for all $u \not\in \mbox{dom}\, s$,
it follows from the hypothesis in part (b)
and from part (a) of Theorem \ref{thm:refinement} that for
all $\beta \in \R^\sigma$
\[
\E_\beta = \bigcup_{u \in \partial \varphi(\beta) \cap
  \mbox{\scriptsize dom}\, s} \E^u=\bigcup_{u \in \partial
  \varphi(\beta) \cap \Gamma} \E^u.
\]
Thus $\partial \varphi(\beta) \cap \mbox{dom}\, s=\partial
\varphi(\beta) \cap \Gamma$.  Taking the union over all $\beta \in
\R^\sigma$ yields
\[
\bigcup_{\beta \in \R^\sigma} \partial \varphi(\beta) \cap
\mbox{dom}\, s = \bigcup_{\beta \in \R^\sigma} \partial \varphi(\beta)
\cap \Gamma \subset \Gamma.
\]
By standard duality theory for upper semicontinuous, concave
functions on $\R^\sigma$ \cite[p.\ 218]{Roc}, $\bigcup_{\beta \in \R^\sigma} \partial \varphi(\beta) =
\mbox{dom}\, \partial s^{**}$.  Thus
\[
(\mbox{dom}\, \partial s^{**}) \cap (\mbox{dom} \, s) \subset \Gamma.
\]
Since
$\mbox{ri(dom}\, s) \subset \mbox{ri(dom}\, s^{**}) \subset
\mbox{dom}\, \partial s^{**}$, we conclude that
$\mbox{ri(dom}\, s) \subset \Gamma$ and therefore that $s$ is concave on
any convex subset of $\mbox{ri(dom}\, s)$.  The proof of the theorem
is complete. \ink

\skp

In the next section we extend the large deviation  
theorems in Sections 2 and 3 and the duality theorems in the
present section to the study of mixed ensembles.

\renewcommand{\theequation}{\arabic{section}.\arabic{subsection}.\arabic{equation}}
\renewcommand{\thedefn}{\arabic{section}.\arabic{subsection}.\arabic{defn}}

\newpage
\section{Mixed Ensembles}
\beginsec
\label{section:mixed}

In broad terms the canonical ensemble differs from the microcanonical
ensemble by the manner in which the dynamical invariants are incorporated
in the respective probability measures: exponentiation in the former ensemble 
and conditioning in the latter ensemble. 
In Section 5.1 we define two classes of mixed ensembles, 
a mixed canonical-microcanonical ensemble and a mixed microcanonical-canonical
ensemble, which differ only in the order in which
the exponentiation and the conditioning are performed.  
In part (b) of Theorem \ref{thm:ldpmixed} we show that
with respect to both of these ensembles the hidden process
$Y_n$ satisfies the large deviation principle with the same rate function.
Hence the sets of equilibrium macrostates for both of these
ensembles are the same.  In Section 5.2 we present 
complete equivalence and nonequivalence results relating
the sets of equilibrium macrostates for the mixed 
and the pure canonical ensembles.  In Section 5.3, we do the same
for the sets of equilibrium
macrostates for the mixed and the pure microcanonical ensembles.
These results will
be applied in future work to a number of problems, including
soliton turbulence for the nonlinear Schr\"{o}dinger equation
\cite{EllJorTur}. 

\subsection{Properties of the Mixed Ensembles}
\beginsec

The definitions of the mixed ensembles involve quantities introduced
in Hypotheses \ref{hyp:prob} and \ref{hyp:general}.
We shall use
the notation $\can(H_n;P_n)_\beta$ to denote the canonical ensemble
$P_{n,\beta}$, which is defined in (\ref{eqn:definecan}), and the notation 
$\micro(H_n;P_n)^{u,r}$
to denote the microcanonical ensemble $P_n^{u,r}$, which is defined in 
(\ref{eqn:definemicro}). 
The LDP's for $Y_n$ with respect to the canonical
ensemble and with respect to the microcanonical ensemble are given
in Theorems \ref{thm:canon} and \ref{thm:microcan}, respectively.
The respective rate functions are
\[
I_{\beta}(x) \doteq I(x) + \lan \beta , \tilde{H}(x) \ran
- \inf_{y \in \X}\{I(y) + \lan \beta,\tilde{H}(y) \ran \}
\]
and for $u \in \mbox{dom} \, J$
\[
\label{eqn:microrate}
I^{u}(x) \doteq \left\{ \begin{array}{ll} I(x) -
J(u) &  \trm{if }\; \tilde{H}(x) = u, \\ 
                                \infty & \trm{otherwise.}   
                               \end{array}  \right.  
\]
In the sequel we shall use the following alternate formula for $I^u$:
\[
I^u(x) = I(\{x\} \cap \tilde{H}^{-1}(\{u\})) - J(u).
\] 
Analogous formulas will arise in the study of the mixed ensembles.

In order to introduce the mixed ensembles, we assume that
$\sigma \geq 2$.  Let $\tau$ be an 
integer satisfying $1 \leq \tau \leq \sigma$ and consider decompositions
of $H_n$ and of $\tilde{H}$ defined as follows:
\[
H_n = (H_n^1,H_n^2), \mbox{ where } H_n^1 \doteq (H_{n,1},\ldots,
H_{n,\tau}) \mbox{ and }  H_n^2 \doteq (H_{n,\tau + 1},\ldots,H_{n,\sigma}),
\]
\[
\tilde{H} = (\tilde{H}^1,\tilde{H}^2), \mbox{ where } 
\tilde{H}^1 \doteq (\tilde{H}_1,\ldots,
\tilde{H}_\tau) \mbox{ and }  \tilde{H}^2 \doteq 
(\tilde{H}_{\tau + 1},\ldots,\tilde{H}_\sigma).
\]
Writing $\beta = (\beta^1,\beta^2) \in \R^\tau
\times \R^{\sigma - \tau}$ and $u = (u^1,u^2) \in
\R^\tau \times \R^{\sigma - \tau}$, we define
\begin{eqnarray*}
\can(H_n^1,H_n^2; P_n)_{\beta^1,\beta^2}(d\omega) & \doteq &
\can(H_n;P_n)_\beta(d\omega) \\
& = & \frac{1}{Z_n(\beta^1,\beta^2)} \, \exp[-\lan \beta^1,H_n^1(\omega) \ran
-\lan \beta^2,H_n^2(\omega) \ran] \, P_n(d\omega),
\end{eqnarray*}
where $Z_n(\beta^1,\beta^2) \doteq Z_n(\beta)$, 
and we define
\begin{eqnarray*}
\micro(H_n^1,H_n^2; P_n)^{u^1,u^2,r}(d\omega) & \doteq  & \micro(H_n;P_n)^{u,r}(d\omega) \\ & = &
P_n(d\omega | H_n^1 \in \{u^1\}^{(r)}, H_n^2 \in \{u^2\}^{(r)}).
\end{eqnarray*}
The function 
$J(u) \doteq \inf\{I(x) : x \in \X, \tilde{H}(x) = u\}$ plays
a key role in the large deviation analysis of the microcanonical ensemble.
We rewrite this function as
\be
\label{eqn:jz1z2}
J(u^1,u^2) \doteq
\inf\{I(x) : x \in \X, \tilde{H}^1(x) = u^1,
\tilde{H}^2(x) = u^2\}.
\ee
 
The innovation of the present subsection is to consider the asymptotic
properties of two mixed ensembles, both at the level of thermodynamic
functions and at the level of equilibrium macrostates.
We define a mixed canonical-microcanonical
ensemble by replacing the measure $P_n$ in 
the canonical ensemble $\can(H_n^1;P_n)_{\beta^1}$ by the microcanonical
ensemble $\micro(H_n^2;P_n)^{u^2,r}$.  For $u^2 \in \R^{\sigma - \tau}$
and $\beta^1 \in \R^\tau$, the resulting measure is given
by
\begin{eqnarray*}
\lefteqn{
\can(H_n^1;\micro(H_n^2;P_n)^{u^2,r})_{\beta^1}(d\omega)} \\
&& \doteq \frac{1}{Z_n(\beta^1,\{u^2\}^{(r)})} \exp[-\lan \beta^1,H_n^1(\omega)
 \ran] \, P_n(d\omega | H_n^2 \in \{u^2\}^{(r)}),
\end{eqnarray*}
where
\[
Z_n(\beta^1,\{u^2\}^{(r)}) \doteq \int_{\Omega_n}
\exp[-\lan \beta^1,H_n^1(\omega) \ran] \, P_n(d\omega|H_n^2 \in 
\{u^2\}^{(r)}).
\]
By a similar verification as in the paragraph after Proposition \ref{prop:ldpj}, 
the microcanonical ensemble $\micro(H_n^2;P_n)^{u^2,r}$,
and thus this mixed ensemble, are well defined for all sufficiently large $n$
provided $u^2$ lies in the domain of 
\be
\label{eqn:jz2}
J^2(u^2) \doteq \inf\{I(x) : x \in \X, \tilde{H}^2(x) = u^2\}.
\ee

In an analogous way, 
we define a mixed microcanonical-canonical ensemble
by replacing the measure $P_n$ in the microcanonical ensemble
$\micro(H_n^2;P_n)^{\utwo,r}$ by the canonical ensemble
$\can(H_n;P_n)_{\bone}$.  For $\beta^1 \in \R^\tau$
and $u^2 \in \R^{\sigma -\tau}$, the resulting measure is given by
\[ 
\micro(H_n^2;\can(H_n^1;P_n)_{\beta^1})^{u^2,r}(d\omega) \doteq
Q_{n,\beta^1}(d\omega | H_n^2 \in \{u^2\}^{(r)}),
\]
where
\[
Q_{n,\beta^1}(d\omega) \doteq \frac{1}{Z_n(\beta^1)} \, 
\exp[-\lan \beta^1,H_n^1(\omega) \ran]
\, P_n(d\omega).
\]
This mixed ensemble is well defined for all sufficiently large
$n$ provided $u^2$ lies in the domain of the function $J_{\beta^1}$
that stands in the same relationship to the mixed ensemble
as the function $J$ in (\ref{eqn:jz1z2}) stands to the microcanonical
ensemble.  Since $J$ is defined in terms of $I$, which is the rate
function in the LDP for $Y_n$ with respect to
$P_n$, $J_{\beta^1}$ is defined in terms of 
the rate function for $Y_n$ with respect
to the canonical ensemble $\can(H_n^1;P_n)_{a_n \beta^1}$.  By Theorem 
\ref{thm:canon}, this rate function is given by
\[
I_{\beta^1}(x) \doteq I(x) + \lan \beta^1,\tilde{H}^1(x)\ran -
\inf_{y \in \X} \{I(y) + \lan \beta^1,\tilde{H}^1(y)\ran \}.
\]
It follows that
\bea
\nonumber 
J_{\beta^1}(u^2) & \doteq &   \inf\{I_{\beta^1}(x) : x \in \X, \tilde{H}^2(x) =
u^2\} \\ \label{eqn:jbz}
& =  & \inf\{I(x) + \lan \beta^1, \tilde{H}^1(x)\ran 
: x \in \X, \tilde{H}^2(x) = u^2\} \\
\nonumber & & \hspace{.15in}
- \inf_{y \in \X}\{I(y) + \lan \beta^1, \tilde{H}^1(y)\ran\}.
\eea
By the discussion earlier in this paragraph, the mixed ensemble 
$\micro(H_n^2;$ $\can(H_n^1;P_n)_{\beta^1})^{u^2,r}$ is well-defined
for all sufficiently large $n$ provided $\utwo$ lies in the domain
of $J_{\beta^1}$.  
Since $\tilde{H}^1(x)$ is finite for all $x \in \X$, 
$u^2 \in \mbox{dom} \, J_{\beta^1}$ if and only if 
$u^2 \in \mbox{dom} \, J^2$.
By the same proof as that of 
Proposition \ref{prop:ldpj}, with respect to $P_n$, 
the sequences $\tilde{H}^2(Y_n)$ and $H_n^2$ satisfy the LDP
on $\R^{\sigma - \tau}$ with rate function $J^2$.
As a consequence, $\mbox{dom} \, J^2$ is nonempty as is
$\mbox{dom} \, \jbone$.

We recall from Section 4 that 
\[
s(u) \doteq -J(u) = - \inf\{I(x) : x \in \X, \tilde{H}(x) = u\}
\] 
defines the
microcanonical entropy and that its Legendre-Fenchel transform
gives the canonical free energy.  Both functions appear in relationships
involving $\ebeta$ and $\eu$ that appear in that section.
In an analogous way, for $\bone \in \R^\tau$
and $\utwo \in \R^{\sigma - \tau}$, we define the entropy 
with respect to the mixed ensemble
$\micro(H_n^2;$ $\can(H_n^1;P_n)_{a_n \beta^1})^{u^2,r}$ to be
\be
\label{eqn:sbone}
s_{\bone}(\utwo) \doteq -J_{\bone}(\utwo).
\ee
This entropy and the associated free energy will appear in 
the results on equivalence and nonequivalence of ensembles 
to be given in Section 5.2.

In order to complete the definitions of the various ensembles,
we also consider the pure ensembles
\[
\can(H_n^1;\can(H_n^2;P_n)_{\beta^2})_{\beta^1} \ \mbox{ and } \
\micro(H_n^1;\micro(H_n^2;P_n)^{u^2,r})^{u^1,r},
\]
which are defined similarly as above.
We omit the simple calculation showing that for all $n$ and $r$
\be
\label{eqn:cancan}
\can(H_n^1;\can(H_n^2;P_n)_{\beta^2})_{\beta^1}(d\omega) =
\can(H_n^1,H_n^2;P_n)_{\beta^1,\beta^2}(d\omega)
\ee
and
\be
\label{eqn:micromicro}
\micro(H_n^1;\micro(H_n^2;P_n)^{u^2,r})^{u^1,r}(d\omega) =
\micro(H_n^1,H_n^2;P_n)^{u^1,u^2,r}(d\omega).
\ee
On the other hand, for all $n$ and $r$ the mixed canonical-microcanonical
ensemble and the mixed microcanonical-canonical ensemble are different.
In the next theorem we record the LDP's satisfied
by $Y_n$ with respect to the various ensembles
introduced in this subsection.  The pleasant surprise is that although
the two mixed ensembles are different for all $n$ and $r$, with respect
to each of them, with $\beta^1$ replaced
by $a_n \beta^1$, $Y_n$ satisfies the LDP with the
identical rate function.  

Before stating the theorem, we define the rate functions for 
each ensemble. 
For $\beta = (\beta^1,\beta^2) 
\in \R^\tau \times \R^{\sigma - \tau}$,
$u^2 \in \mbox{dom} \, J^2$, 
and $u = (u^1,u^2) \in \mbox{dom} \, J$, we define
the following functions mapping $\X$ into $[0,\infty]$:
\bea
\label{eqn:ibetabeta}
I_{\beta^1,\beta^2}(x) & \doteq & I(x) + \lan \beta^1, \tilde{H}^1(x) \ran
+ \lan \beta^2, \tilde{H}^2(x) \ran \\ \nonumber
& & \hspace{.15in} - \inf_{y \in \X} \{ I(y) + \lan \beta^1, \tilde{H}^1(y) \ran
+ \lan \beta^2, \tilde{H}^2(y) \ran\},
\eea
\begin{eqnarray}
\label{eqn:ibetaz}
I_{\beta^1}^{u^2}(x) & \doteq & I(\{x\} \cap (\thtwo)^{-1}(\{\utwo\})) + \lan \beta^1,\tilde{H}^1(x) \ran \\
&& \hspace{.15in} - \inf\{I(y) + \lan \beta^1,\tilde{H}^1(y) \ran : y \in \X,
\tilde{H}^2(y) = u^2\}, \nonumber
\end{eqnarray} 
and 
\be
I^{u^1,u^2}(x) \doteq I(\{x\} \cap 
(\thone)^{-1}(\{\uone\}) \cap (\thtwo)^{-1}(\{\utwo\})) - J(u^1,u^2).
\label{eqn:izz}
\ee

\begin{thm} \per
\label{thm:ldpmixed}
We assume Hypotheses {\em \ref{hyp:prob}} 
and {\em \ref{hyp:general}}.
For $(\bone,\btwo) \in \R^\tau \times \R^{\sigma - \tau}$ 
the following conclusions hold.

{\em (a)}   With respect to the canonical ensemble 
$\ccan({H}_n^1,{H}_n^2;P_n)_{a_n \beta^1,a_n \beta^2}$, 
$Y_n$ satisfies the LDP on $\X$
with rate function 
$I_{\beta^1,\beta^2}$ given in {\em (\ref{eqn:ibetabeta})}.  

{\em (b)} Take $u^2 \in \mbox{{\em dom}} \, J^2$ {\em [}see 
{\em (\ref{eqn:jz2})}{\em ]}.
Both with respect to the mixed canonical-microcanonical ensemble
$\ccan(H_n^1;\mmicro(H_n^2;P_n)^{u^2,r})_{a_n \beta^1}$
and with respect to the mixed microcanonical-canonical ensemble
$\mmicro(H_n^2;\ccan(H_n^1;P_n)_{a_n \beta^1})^{u^2,r}$, $Y_n$ satisfies
the LDP on $\X$, in the double limit
$\ngi$ and $r \goto 0$, with rate function $I_{\beta^1}^{u^2}$ given in 
{\em (\ref{eqn:ibetaz})}.  

{\em (c)} Take $u = (u^1,u^2) \in \mbox{{\em dom}} \, J$ {\em [}see {\em 
(\ref{eqn:jz1z2})}{\em ]}.  With respect to the microcanonical
ensemble $\mmicro(H_n^1,H^2_n;P_n)^{u^1,u^2,r}$, $Y_n$
satisfies the LDP on $\X$, in the double
limit $\ngi$ and $r \goto 0$, with 
rate function $I^{u^1,u^2}$ given in {\em (\ref{eqn:izz})}.
\end{thm}

\noi
{\bf Proof.}  Part (a) is proved in Theorem \ref{thm:canon},
and part (c) is proved 
in Theorem \ref{thm:microcan}.  In part (b) we first prove
the LDP for $Y_n$ with respect to
$\micro(H_n^2;$$\can(H_n^1;P_n)_{a_n \beta^1})^{u^2,r}$.
Theorem \ref{thm:canon} implies that with respect to
$\can(H_n^1;P_n)_{a_n \beta^1}$, $Y_n$ satisfies the LDP with rate function
\[
I_{\beta^1}(x) \doteq I(x) + \lan \beta^1,\tilde{H}^1(x) \ran 
- \inf_{y \in \X} \{I(y) + \lan \beta^1,\tilde{H}^1(y)\}.
\]
With $P_n$ replaced by $\can(H_n^1;P_n)_{a_n \beta^1}$
and $I$ replaced by $I_{\beta^1}$, Theorem \ref{thm:microcan}
guarantees that if $u^2 \in \mbox{dom} \, J_{\beta^1} =
\mbox{dom} \, J^2$, then with respect to 
$\micro(H_n^2;$$\can(H_n^1;P_n)_{a_n \beta^1})^{u^2,r}$,
$Y_n$ satisfies the LDP,
in the double limit $\ngi$ and $r \goto 0$, with rate function
\[
(I_{\beta^1})^{u^2}(x) \doteq  I_{\beta^1}(\{x\} \cap (\thtwo)^{-1}(\{\utwo\})) 
- \inf\{I_{\beta^1}(y) : y \in \X,
\tilde{H}^2(y) = u^2\}.
\]
Substituting the definition of $I_{\beta^1}$, we see that
\beas
(I_{\beta^1})^{u^2}(x) & = & I(\{x\} \cap (\thtwo)^{-1}(\{\utwo\})) + \lan \beta^1,\tilde{H}^1(x) \ran  \\ 
& & \hspace{.15in} - \inf\{I(y) + \lan \beta^1,\tilde{H}^1(y)\ran : y \in \X, \tilde{H}^2(y) = u^2\}.
\eeas
This is the function $I_{\beta^1}^{u^2}$ defined in (\ref{eqn:ibetaz}).
We have proved that with respect to 
$\micro(H_n^2;$ $\can(H_n^1;P_n)_{a_n \beta^1})^{u^2,r}$,
$Y_n$ satisfies the LDP, in the double
limit $\ngi$ and $r \goto 0$, with rate function $I_{\beta^1}^{u^2}$.

We next consider the LDP for $Y_n$ with
respect to 
$\can(H_n^1;$ $\micro(H_n^2;P_n)^{u^2,r})_{a_n \beta^1}$.  Since
$u^2 \in \mbox{dom} \, J^2$, Theorem \ref{thm:microcan}
implies that with respect to $\micro(H_n^2;P_n)^{u^2}$, $Y_n$
satisfies the LDP, in the double limit
$\ngi$ and $r \goto 0$, with rate function 
\[
I^{u^2}(x) \doteq I(\{x\} \cap (\thtwo)^{-1}(\{\utwo\})) - J^2(u^2).
\]
One can easily modify the proof of Theorem \ref{thm:canon} to handle
the situation in which $P_n$ is replaced by a doubly indexed
class of probability measures such as $\micro(H_n^2;P_n)^{u^2,r}$ with 
the property 
that with respect to these measures $Y_n$ satisfies the LDP.
With this modification, 
replacing $P_n$ by $\micro(H_n^2;P_n)^{u^2,r}$ and $I$
by $I^{u^2}$, we see that with respect to 
$\can(H_n^1;$ $\micro(H_n^2;P_n)^{u^2,r})_{a_n \beta^1}$,
$Y_n$ satisfies the LDP,
in the double limit $\ngi$ and $r \goto 0$, with rate function
\begin{eqnarray*}
(I^{u^2})_{\beta^1}(x) & \doteq &
I^{u^2}(x) + \lan \beta^1,\tilde{H}^1(x) \ran  - 
\inf_{y \in \X}\{I^{u^2}(y) + \lan \beta^1,\tilde{H}^1(y) \ran\} \\
& = & I(\{x\} \cap (\thtwo)^{-1}(\{\utwo\})) 
+ \lan \beta^1,\tilde{H}^1(x) \ran \\
&& \hspace{.15in} - \ 
\inf\{I(y) + \lan \beta^1,\tilde{H}^1(y) \ran
: y \in \X, \tilde{H}^2(y) = u^2\}.
\end{eqnarray*}
This is the function $I_{\beta^1}^{u^2}$ defined in (\ref{eqn:ibetaz}).
We have shown that with respect to 
$\can(H_n^1;$ $\micro(H_n^2;P_n)^{u^2,r})_{a_n \beta^1}$,
$Y_n$ satisfies the LDP, in the double
limit $\ngi$ and $r \goto 0$, with rate function $I_{\beta^1}^{u^2}$.
The proof of the theorem is complete.  \ink

\skp

In the next two subsections, we consider equivalence and nonequivalence
results for the ensembles
whose LDP's are derived 
in Theorem \ref{thm:ldpmixed}.   These results
are derived as immediate consequences of our work in Section 4,
where equivalence and nonequivalence results for the canonical
and microcanonical ensembles were derived.

\subsection{Equivalence and Nonequivalence of the Canonical and Mixed Ensembles}
\beginsec
\setcounter{defn}{0}

\begin{figure}[t]
\begin{center}
\begin{Bundle}{$(\beta^1,\beta^2) \in \R^\tau \times \R^{\sigma-\tau}$}
  \chunk{\ovalbox{\makebox[4.5cm][c]{$\E_{\beta^1,\beta^2} =
        \displaystyle{\bigcup_{u^2 \in
            \tilde{H}^2(\E_{\beta^1,\beta^2})}} \E_{\beta^1}^{u^2}$ }}}
\end{Bundle}\\
\begin{itemize}
\item[(a)] For $(\beta^1,\beta^2) \in \R^\tau \times \R^{\sigma-\tau}$, any $x
  \in \E_{\beta^1,\beta^2}$ lies in some $\E_{\beta^1}^{u^2}$.
\end{itemize}
\setlength{\GapDepth}{5mm} \drawwith{\drawline[10000]}
\begin{Bundle}{$(\beta^1,u^2) \in \R^\tau \times \mbox{dom} \, s_{\beta^1}$}
  \chunk{\begin{Bundle}{$u^2 \in C_{\beta^1}$}
      \chunk{\begin{Bundle}{$u^2 \in T_{\beta^1}$}
          \chunk{\ovalbox{\begin{Bcenter} \textbf{Full Equivalence:}\\ 
                \makebox[4.5cm][c]{$\exists \beta^2 \in \partial
                  s_{\beta^1}^{**}(u^2)$ s.t.}\\ 
                $\E_{\beta^1}^{u^2} = \E_{\beta^1,\beta^2}$
            \end{Bcenter}}}
      \end{Bundle}}
    \chunk{\begin{Bundle}{$u^2 \not\in T_{\beta^1}$}
        \chunk{\ovalbox{\begin{Bcenter}
              \textbf{Partial Equivalence:}\\
              \makebox[4.5cm][c]{$\E_{\beta^1}^{u^2} \subsetneq
                \E_{\beta^1, \beta^2}$}\\
              $\forall \beta^2 \in \partial s_{\beta^1}^{**}(u^2)$ 
            \end{Bcenter}}}
      \end{Bundle}}
  \end{Bundle}}
\chunk{%
  \setlength{\GapDepth}{14.5mm}
  \begin{Bundle}{$u^2 \not\in C_{\beta^1}$}
    \chunk{\ovalbox{\begin{Bcenter}
          \textbf{Nonequivalence:}\\
          \makebox[4.5cm][c]{$\E_{\beta^1}^{u^2} \cap \E_{\beta^1,
              \beta^2} = \emptyset$}\\
          $\forall \beta^2 \in \R^{\sigma-\tau}$
        \end{Bcenter}}}
  \end{Bundle}}
\end{Bundle}
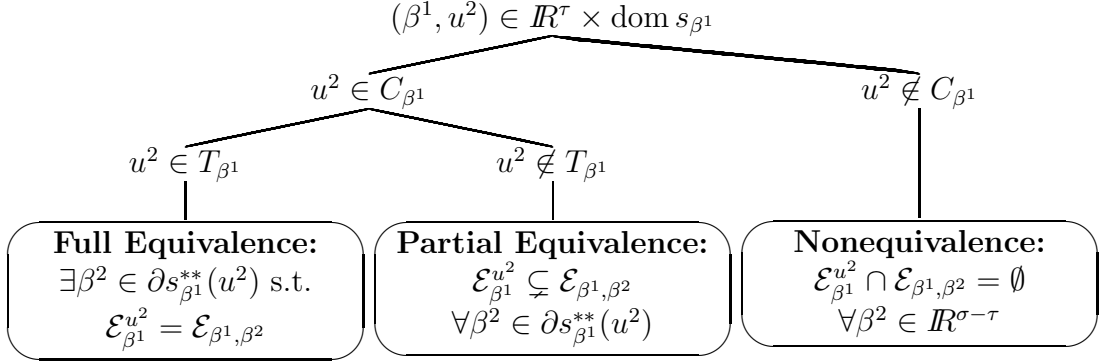
\begin{itemize}
\item[(b)] For $\beta^1 \in \R^\tau$, there are three possibilities for 
  $u^2 \in \mbox{dom} \, s_{\beta^1}$.  The two branches on the left
  lead to equivalence results, whereas the other branch leads to a
  nonequivalence result.  The sets $C_{\beta^1}$ and $T_{\beta^1}$ are
  defined in the last paragraph of Section 5.2.
\end{itemize}
\caption{{}Equivalence and nonequivalence of canonical and mixed ensembles.}
\label{fig:fixed-beta-1}
\end{center}
\end{figure}

In this subsection we study, at the level of equilibrium macrostates, the equivalence and nonequivalence of the canonical
ensemble $\can(\hnone,\hntwo;P_n)_{a_n \beta^1,a_n \beta^2}$ 
and the mixed ensemble
$\micro(\hntwo;\can(\hnone;P_n)_{a_n \beta^1})^{u^2,r}$.  The parameters
$\beta^1$, $\beta^2$, and $u^2$ satisfy
$\beta^1 \in \R^\tau$, $\beta^2 \in \R^{\sigma - \tau}$, 
and $u^2 \in \mbox{dom} \, J^2$, where
\[
J^2(u^2) \doteq \inf\{I(x) : \tilde{H}^2(x) = u^2\}.
\]
By a similar verification as in the paragraph after Proposition \ref{prop:ldpj},
this condition on $\utwo$
guarantees that the mixed ensemble is well defined
for all sufficiently large $n$.
The relationships between the sets of equilibrium macrostates
for the two ensembles follow immediately from Theorems \ref{thm:z},
\ref{thm:beta}, and \ref{thm:full} with minimal changes in proof.  
Hence we shall only
summarize them in Figure \ref{fig:fixed-beta-1}.

By Theorem \ref{thm:ldpmixed}, for $(\bone,\btwo) \in
\R^\tau \times \R^{\sigma - \tau}$, with respect to
$\can(\hnone,\hntwo;$ $P_n)_{a_n \beta^1,a_n \beta^2}$
$\, Y_n$ satisfies the LDP with rate function
\bea
\label{eqn:ibetabeta2}
I_{\beta^1,\beta^2}(x) & \doteq & I(x) + \lan \beta^1, \tilde{H}^1(x) \ran
+ \lan \beta^2, \tilde{H}^2(x) \ran \\ \nonumber
& & \hspace{.15in} - \inf_{y \in \X} \{ I(y) + \lan \beta^1, \tilde{H}^1(y) \ran
+ \lan \beta^2, \tilde{H}^2(y) \ran\}.
\eea
In addition, for $(\bone,\utwo) \in \R^\tau \times \mbox{dom} \, J^2$, 
with respect to 
$\micro(\hntwo;\can(\hnone;$ $P_n)_{a_n \beta^1})^{u^2,r}$ $\, Y_n$ satisfies
the LDP with rate function
\be
\label{eqn:ibetaz2}
I_{\beta^1}^{u^2}(x) \doteq I(\{x\} \cap (\thtwo)^{-1}(\{\utwo\})) + \lan \beta^1,\tilde{H}^1(x) \ran) - \psi_{\beta^1}^{u^2},
\ee
where 
\[
\psi_{\beta^1}(\utwo) \doteq  \inf\{I(y) + \lan \beta^1,\tilde{H}^1(y)\ran : y \in \X, \tilde{H}^2(y)= u^2\}.
\]
For $\bone \in \R^\tau$, $\btwo \in \R^{\sigma - \tau}$, and 
$\utwo \in \mbox{dom} \, J^2$,
we define the corresponding sets of equilibrium macrostates
\[
{\cal E}_{\beta^1,\beta^2} \doteq 
\{x \in \X : I_{\beta^1,\beta^2}(x) = 0\}
\]
and
\beas
{\cal E}_{\beta^1}^{u^2} & \doteq & \{x \in \X :
I_{\beta^1}^{u^2}(x) = 0\} \\
& = & \{x \in \X : \tilde{H}^2(x) = \utwo,
I(x) + \lan \bone,\thone(x) \ran = \psi_{\beta^1}^{\utwo} \}.
\eeas
As the sets of points at which the corresponding rate functions attain their
minimum of 0, 
both ${\cal E}_{\bone,\btwo}$ and ${\cal E}_{\beta^1}^{\utwo}$ 
are nonempty, compact subsets of $\X$
for $\bone \in \R^\tau$, $\btwo \in \R^{\sigma - \tau}$,
and $\utwo \in \mbox{dom} \, J^2$.
The main purpose of this subsection is to record the relationships
between these sets.

Before doing so, we point out a concentration property, relative
to the set ${\cal E}_{\bone}^{\utwo}$, of the distributions of $Y_n$
with respect to the mixed ensemble 
$\micro(H_n^2;\can(H_n^1;P_n)_{a_n \beta^1})^{u^2,r}$.  This concentration
property is an immediate consequence of the LDP
proved in part (b) of Theorem \ref{thm:ldpmixed}.  It 
justifies calling ${\cal E}_{\bone}^{\utwo}$ the set of
equilibrium macrostates with respect to the mixed ensemble.
This concentration property is analogous to those
for the canonical ensemble and for the microcanonical ensemble
given in part (c) of Theorem \ref{thm:canon} and in 
part (b) of Theorem \ref{thm:iz};
the proof is omitted.

\begin{thm} \per
\label{eqn:concmixed}
We assume Hypotheses {\em \ref{hyp:prob}} and {\em \ref{hyp:general}}.
For $\bone \in \R^\tau$, $\utwo \in 
\mbox{{\em dom}} \, J^2$, and $A$ any Borel subset
of $\X$ whose closure $\bar{A}$ satisfies 
$\bar{A} \cap {\cal E}_{\bone}^{\utwo}
= \emptyset$, we have $I_{\bone}^{\utwo}(\bar{A}) > 0$.  In addition,
there exists $r_0 \in (0,1)$ and for all $r \in (0,r_0]$
there exists $C_r < \infty$ such that
\[
\mmicro(H_n^2;\ccan(H_n^1;P_n)_{a_n \beta^1})^{u^2,r}\{Y_n \in A\} 
\leq C_r \exp[-a_n I_{\bone}^{\utwo}(\bar{A})/2] 
\goto 0 \: \mbox{ as } n \goto \infty.
\]
\end{thm}

As in Theorem \ref{thm:microlimits}, one can also study compactness
and weak limit properties of the distributions of $Y_n$ with respect
to $\micro(H_n^2;\can(H_n^1;P_n)_{a_n \beta^1})^{u^2,r}$.  We shall omit
this topic.

We return to the relationships between
${\cal E}_{\beta^1,\beta^2}$ and ${\cal E}_{\beta^1}^{u^2}$.
Since for each $n$ 
\[
\can(H_n^1,H_n^2;P_n)_{\beta^1,\beta^2} \ \mbox{ and } \
\can(H_n^2;\can(H_n^1;P_n)_{\beta^1})_{\beta^2}
\]
are equal, we can derive the relationships between these
sets of equilibrium macrostates
by applying the results of Section 4 to the canonical ensemble
and microcanonical ensemble
\[
\can(H_n^2;Q_n)_{a_n \beta^2} \ \mbox{ and } \
\micro(H_n^2;Q_n)^{u^2}, \mbox{ with } Q_n \doteq
\can(H_n^1;P_n)_{a_n \beta^1}.
\]
To this end, we introduce the relevant thermodynamic functions.
With respect to $\can(H_n^2;$ $\can(H_n^1;P_n)_{a_n \beta^1})_{a_n \beta^2}$
the free energy is given by
\bea
\label{eqn:vphibone}
\vphi_{\bone}(\btwo) & = & -\lim_{\ngi} \frac{1}{a_n} \log
\int_{\Omega_n} \exp[-a_n \lan \btwo,\hntwo \ran ] \, 
d\!\left(\can\rule{0mm}{4mm}(\hnone;P_n)_{a_n \bone}\right) 
\\ \nonumber & = & \inf_{x \in \X} \{I(x) + \lan \bone,\thone(x) \ran +
\lan \btwo,\thtwo(x) \ran \} - \vphi^1(\beta^1),
\eea
where 
\bea
\label{eqn:vphi1}
\vphi^1(\beta^1) & \doteq &
-\lim_{\ngi} \frac{1}{a_n} \log
\int_{\Omega^n} \exp[-a_n \lan \bone,\hnone \ran] \, dP_n  \\
\nonumber & = & \inf_{y \in \X}\{I(y) + \lan \bone,\thone(y)\ran \}.
\eea
The function $\vphi_{\bone}$ is finite, concave, and continuous
on $\R^{\sigma - \tau}$.
In (\ref{eqn:sbone}) we identified the entropy with respect to 
$\micro(\hntwo;\can(\hnone;P_n)_{a_n \bone})^{\utwo,r}$ to be
\bea
s_{\bone}(\utwo) & \doteq & -J_{\bone}(\utwo) \nonumber \\
& = & -\inf\{I_{\bone}(x) : x \in \X, \thtwo(x) = \utwo\}
\label{eqn:jbtwo}  \\ \nonumber & = & 
-\inf\{I(x) + \lan \bone,\thone(x) \ran : x \in \X,
\thtwo(x) = \utwo\} + \vphi^1(\bone);
\eea
$\utwo \in \mbox{dom} \, s_{\beta^1}$ if and only if $u^2 \in
\mbox{dom} \, J^2$. 

As in Section \ref{section:equiv}, whether or not
the entropy $s_{\bone}$ 
is concave on $\R^{\sigma-\tau}$, its Legendre-Fenchel
transform $s_{\bone}^*$ equals $\vphi_{\bone}$.
If in addition $s_{\bone}$ is concave on $\R^{\sigma-\tau}$, then this 
formula can be inverted to give $s_{\bone} = \vphi_{\bone}^*$.

For $\beta^1 \in \R^\tau$ the relationships between $\E_{\beta^1,\beta^2}$ 
and $\E_{\beta^1}^{u^2}$ are summarized in Figure \ref{fig:fixed-beta-1}.
These relationships depend on two sets that are the analogues
of the sets $C$ and $T$ defined in (\ref{eqn:C}) and (\ref{eqn:T}).
For $\beta^1 \in \R^\tau$ we define
$C_{\beta^1}$ to be the set of $u^2 \in \R^{\sigma-\tau}$ for which there
exists $\beta^2 \in
\R^{\sigma-\tau}$ such that 
\[
s_{\beta^1}(w) \leq s_{\beta^1}(u^2) +
\langle \beta^2, w-u^2 \rangle \mbox{ for all } w \in \R^{\sigma-\tau}.
\]  
We also define $T_{\beta^1}$
to be the set of $u^2 \in \R^{\sigma - \tau}$ for which there
exists $\beta^2 \in
\R^{\sigma - \tau}$ such that 
\[
s_{\beta^1}(w) < s_{\beta^1}(u^2) +
\langle \beta^2, w-u^2 \rangle \mbox{ for all } w \not= u^2.
\]  
As in Lemma \ref{lem:CtoC}, it can be shown that
$C_{\beta^1} = \Gamma_{\beta^1} \cap \mbox{dom}\, \partial
s_{\beta^1}^{**}$,
where
$\Gamma_{\beta^1} \doteq \{ u^2 \in \R^{\sigma - \tau}: s_{\beta^1}(u^2)=
s_{\beta^1}^{**}(u^2) \}$.

\subsection{Equivalence and Nonequivalence of the Mixed and Microcanonical Ensembles}
\beginsec
\setcounter{defn}{0}

\begin{figure}[t]
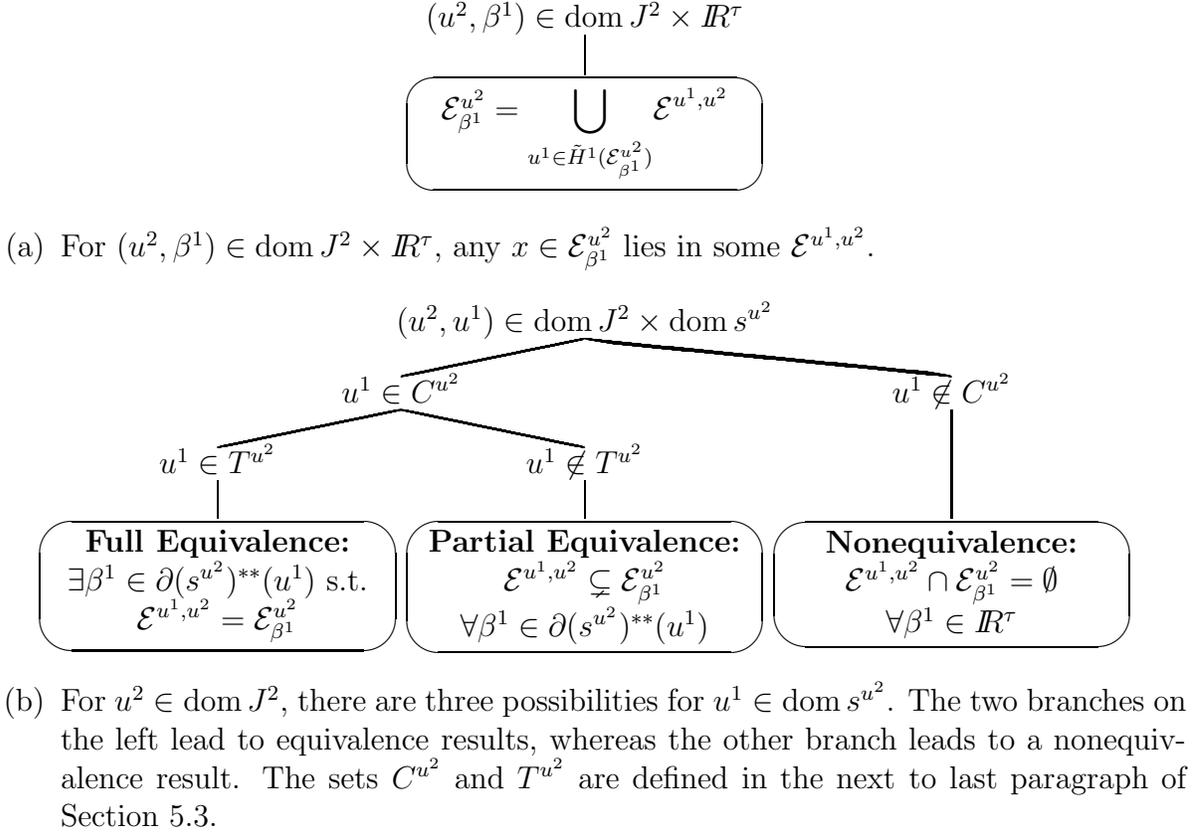

\begin{center}
\begin{Bundle}{($u^2,\bone) \in \mbox{dom} \, J^2 \times \R^\tau$}
  \chunk{\ovalbox{\makebox[4.5cm][c]{ 
         $\E_{\beta^1}^{u^2}= \displaystyle{\bigcup_{u^1 \in
            \tilde{H}^1(\E_{\beta^1}^{u^2})}} \E^{u^1,u^2}$}}}
\end{Bundle}\\
\begin{itemize}
\item[(a)] For $(u^2,\bone) \in\mbox{dom}\, J^2 \times \R^\tau$, any $x
  \in \E_{\beta^1}^{u^2}$ lies in some $\E^{u^1,u^2}$. 
\end{itemize}
\setlength{\GapDepth}{5mm} \drawwith{\drawline[10000]}
\begin{Bundle}{$(u^2,u^1) \in \mbox{dom}\, J^2 \times \mbox{dom} \, s^{u^2}$}
  \chunk{\begin{Bundle}{$u^1 \in C^{u^2}$} \chunk{\begin{Bundle}{$u^1
            \in T^{u^2}$} \chunk{\ovalbox{\begin{Bcenter} 
                \textbf{Full Equivalence:}\\ 
                \makebox[4.5cm][c]{$\exists \beta^1 \in \partial
                  (s^{u^2})^{**}(u^1)$ s.t.}\\ $\E^{u^1,u^2} =
                \E_{\beta^1}^{u^2}$
            \end{Bcenter}}}
      \end{Bundle}}
    \chunk{\begin{Bundle}{$u^1 \not\in T^{u^2}$}
        \chunk{\ovalbox{\begin{Bcenter}
              \textbf{Partial Equivalence:}\\
              \makebox[4.5cm][c]{$\E^{u^1,u^2} \subsetneq
                \E_{\beta^1}^{u^2}$}\\
              $\forall \beta^1 \in \partial (s^{u^2})^{**}(u^1)$
            \end{Bcenter}}}
      \end{Bundle}}
  \end{Bundle}}
\chunk{%
  \setlength{\GapDepth}{14.5mm}
  \begin{Bundle}{$u^1 \not\in C^{u^2}$}
    \chunk{\ovalbox{\begin{Bcenter} \textbf{Nonequivalence:}\\ 
          \makebox[4.5cm][c]{$\E^{u^1,u^2} \cap \E_{\beta^1}^{u^2} =
            \emptyset$}\\ $\forall \beta^1 \in \R^{\tau}$
        \end{Bcenter}}}
  \end{Bundle}}
\end{Bundle}
\begin{itemize}
\item[(b)] For $u^2 \in \mbox{dom}\, J^2$, there are three
  possibilities for $u^1 \in \mbox{dom} \, s^{u^2}$.  The two
  branches on the left lead to equivalence results, whereas the other
  branch leads to a nonequivalence result.  The sets $C^{u^2}$ and
  $T^{u^2}$ are defined in the next to last paragraph of Section 5.3.
\end{itemize}
\caption{{}Equivalence and nonequivalence of mixed and microcanonical
  ensembles.}
\label{fig:fixed-z-2}
\end{center}
\end{figure}

In this subsection we study, at the level of equilibrium macrostates, 
the equivalence and nonequivalence of the mixed
ensemble $\can(\hnone;\micro(\hntwo;P_n)^{\utwo,r})_{a_n \bone}$ 
and the microcanonical ensemble
$\micro(\hnone,\hntwo;P_n)^{\uone,u^2,r}$.  The parameters
$\beta^1$, $\uone$, and $u^2$ satisfy
$\beta^1 \in \R^\tau$, $\utwo \in \mbox{dom} \, J^2$, 
and $(\uone,u^2) \in \mbox{dom} \, J$, where
\[
J^2(u^2) \doteq \inf\{I(x) : x \in \X, \tilde{H}^2(x) = u^2\}
\]
and 
\[
J(\uone,\utwo) \doteq \inf\{I(x) : x \in \X, 
\thone(x) = \uone, \tilde{H}^2(x) = u^2\}.
\]
For any
$\uone$ and $\utwo$, $J^2(\utwo) \leq J(\uone,\utwo)$.  Hence, if $(\uone,\utwo)
\in \mbox{dom} \, J$, then $\utwo \in \mbox{dom} \, J^2$.
By a similar verification as in the paragraph after Proposition \ref{prop:ldpj},
the condition that $(\uone,\utwo) \in \mbox{dom} \, J$  
guarantees that both the mixed ensemble and the microcanonical
ensemble are well defined for all sufficiently large $n$.
The relationships between the sets of equilibrium macrostates
for the two ensembles follow immediately from Theorems \ref{thm:z},
\ref{thm:beta}, and \ref{thm:full} with minimal changes in proof.  
Hence we shall only
summarize them in Figure \ref{fig:fixed-z-2}.

By Theorem \ref{thm:ldpmixed}, for $(\bone,\utwo) \in
\R^\tau \times (\mbox{dom} \, J^2)$, with respect to
$\can(\hnone;$ $\micro(\hntwo;P_n)^{\utwo,r})_{a_n \bone}$
$\, Y_n$ satisfies the LDP with rate function
\be
\label{eqn:ibetaz22}
I_{\beta^1}^{u^2}(x) \doteq I(\{x\} \cap (\thtwo)^{-1}(\{\utwo\})) + \lan \beta^1,\tilde{H}^1(x) \ran - \psi_{\bone}^{\utwo},
\ee
where
\[
\psi_{\bone}^{\utwo} \doteq 
\inf\{I(y) + \lan \beta^1,\tilde{H}^1(y) \ran : y \in \X,
\tilde{H}^2(y) = u^2\}.
\]
In addition, for $(\uone,\utwo) \in \mbox{dom} \, J$, with respect to 
$\micro(\hnone,\hntwo;P_n)^{\uone,u^2,r}$
$\, Y_n$ satisfies the LDP with rate function
\[
I^{u^1,u^2}(x) \doteq  I(\{x\} \cap (\tilde{H}^1)^{-1}(\{u^1\}) \cap
(\tilde{H}^2)^{-1}(\{u^2\})) - J(u^1,u^2).
\] 
For $\bone \in \R^\tau$, $\utwo \in \mbox{dom} \, J^2$, and $(\uone,\utwo)
\in \mbox{dom} \, J$, we define the corresponding sets of equilibrium
macrostates
\beas
{\cal E}_{\beta^1}^{u^2} & \doteq & \{x \in \X :
I_{\beta^1}^{u^2}(x) = 0\} \\
& = & \{x \in \X : \tilde{H}^2(x) = \utwo,
I(x) + \lan \bone,\thone(x) \ran = \psi_{\beta^1}^{\utwo} \}
\eeas
and
\beas
\euoneutwo & \doteq & \{x \in \X: \iuoneutwo(x) = 0\} \\
& = & \{x \in \X : I(x) = J(\uone,\utwo), \thone(x) = \uone,
\thtwo(x) = \utwo\}.
\eeas
As the sets of points at which the corresponding rate functions
attain their minimum of 0,
the set $\eboneutwo$, for $\bone \in \R^\tau$ and 
$\utwo \in \mbox{dom} \, J^2$, and the set $\euoneutwo$, for
$(\uone,\utwo) \in \mbox{dom} \, J$, are nonempty and compact.
The purpose of this subsection is to record the relationships
between these sets.

Since for $(\uone,\utwo) \in \mbox{dom} \, J$ and each $n$ 
\[
\micro(H_n^1,H_n^2;P_n)^{u^1,u^2,r} \ \mbox{ and } \
\micro(H_n^1;\micro(H_n^2;P_n)^{\utwo,r})^{\uone,r}
\]
are equal, we can derive the relationships between $\eboneutwo$
and $\euoneutwo$
by applying the results of Section 4 to the canonical ensemble
and microcanonical ensemble
\[
\can(H_n^1;Q_n)_{a_n \beta^2} \ \mbox{ and } \
\micro(H_n^1;Q_n)^{u^1,r}, \mbox{ with } Q_n \doteq
\micro(H_n^2;P_n)^{\utwo,r}.
\]
To this end, we introduce the relevant thermodynamic functions.
By Theorem \ref{thm:microcan}, for $\utwo \in \mbox{dom} \, J^2$
the rate function in the LDP for $Y_n$
with respect to $\micro(\hntwo;P_n)^{\utwo,r}$ is
\[
I^{\utwo}(x) \doteq I(\{x\} \cap (\thtwo)^{-1}(\{\utwo\})) - J^2(\utwo).
\]
Hence by the Laplace principle, for $\utwo \in \mbox{dom} \, J^2$
the free energy with respect to the ensemble
$\can(\hnone;\micro(\hntwo;P_n)^{\utwo,r})_{a_n \bone}$ is given by
\bea
\nonumber 
\vphi^{\utwo}(\bone) & = & -\lim_{\ngi} \frac{1}{a_n} \log
\int_{\Omega^n} \exp[-a_n \lan \bone,\hnone \ran ] \, 
d\!\left(\micro\rule{0mm}{4mm}(\hntwo;P_n)^{\utwo,r}\right) 
\\ \label{eqn:vphiztwo}
& = & \inf_{x \in \X} \{I^{\utwo}(x) + \lan \bone,\thone(x) \ran\} 
 \\ & = &  \nonumber
\inf\{I(x) + \lan \bone,\thone(x) \ran : x \in \X, \thtwo(x) = \utwo\}
- J^2(\utwo).
\eea
The function $\vphi^{\utwo}$ is finite, concave, and continuous on $\R^\tau$.
For $\utwo \in \mbox{dom} \, J^2$ we define
\bea
J^{\utwo}(\uone) 
& \doteq & \inf\{I^{\utwo}(x) : x \in \X, \thone(x) = \uone\} \nonumber
\\ \label{eqn:jztwo} & = & \inf\{I(x): x \in \X, \thone(x) = \uone,
\thtwo(x) = \utwo\} - J^2(\utwo) \\ \nonumber
& = & J(\uone,\utwo) - J^2(\utwo).
\eea
With respect to 
$\micro(H_n^1;\micro(H_n^2;P_n)^{\utwo,r})^{\uone,r}$,
for $\utwo \in \mbox{dom} \, J^2$
the entropy is given by 
\be
\label{eqn:sztwo}
s^{\utwo}(\uone) \doteq -J^{\utwo}(\uone).
\ee
We have $\uone \in \mbox{dom} \, s^{\utwo}$ if and only if $(\uone,\utwo) \in
\mbox{dom} \, J$. 

As in Section \ref{section:equiv}, whether or not
$s^{\utwo}$ is concave on $\R^\tau$, its Legendre-Fenchel 
transform $(s^{\utwo})^*$ equals $\vphi^{\utwo}$.
If $\sutwo$ is concave on $\R^{\tau}$, then this formula
can be inverted to give $\sutwo= (\vphiutwo)^*$
for all $\uone \in \R^{\tau}$.

For $\utwo \in \mbox{dom} \, J^2$ the relationships between
$\E_{\beta^1}^{\utwo}$ and $\E^{\uone,u^2}$ are summarized in Figure
\ref{fig:fixed-z-2}.  These relationships depend on two sets that are
the analogues of the sets $C$ and $T$ defined in (\ref{eqn:C}) and
(\ref{eqn:T}).  For $\beta^1 \in \R^\tau$ we define
$C^{u^2}$ to be the set of $u^1 \in \R^{\tau}$ for which there exists
$\beta^1 \in \R^{\tau}$ such that 
\[
s^{u^2}(w) \leq s^{u^2}(u^1) +
\langle \beta^1, w-u^1 \rangle \mbox{ for all } w \in \R^{\tau}.
\]  
We also define $T^{u^2}$ 
to be the set of $u^1 \in \R^{\tau}$ for which there exists
$\beta^1 \in \R^{\tau}$ such that
\[
s^{u^2}(w) < s^{u^2}(u^1) +
\langle \beta^1, w-u^1 \rangle \mbox{ for all } w \not= u^1.
\]
As in Lemma \ref{lem:CtoC}, it can be shown that
$C^{u^2} = \Gamma^{u^2} \cap \mbox{dom}\, \partial (s^{u^2})^{**}$,
where
$\Gamma^{u^2} \doteq \{ u^1 \in \R^{\tau}: s^{u^2}(u^1)=
(s^{u^2})^{**}(u^1) \}$.

With Figure \ref{fig:fixed-z-2}, we complete our presentation 
of the equivalence 
and nonequivalence results for the mixed ensemble, the canonical
ensemble, and the microcanonical ensemble.

\renewcommand{\theequation}{\arabic{section}.\arabic{equation}}
\renewcommand{\thedefn}{\arabic{section}.\arabic{defn}}

\renewcommand{\theequation}{\arabic{section}.\arabic{subsection}.\arabic{equation}}
\renewcommand{\thedefn}{\arabic{section}.\arabic{subsection}.\arabic{defn}}

\renewcommand{\theequation}{\arabic{section}.\arabic{equation}}
\renewcommand{\thedefn}{\arabic{section}.\arabic{defn}}

\skp
\skp

\end{document}